\documentclass[11pt,letterpaper]{amsart}
\usepackage{lmodern}
\usepackage[T1]{fontenc}
\usepackage[latin9]{inputenc}
\usepackage{amsbsy}
\usepackage{amstext}
\usepackage{amsthm}
\usepackage{amssymb}
\usepackage{graphicx}
\usepackage{floatrow}
\usepackage[dvipsnames]{xcolor}
\usepackage{tikz}
\usepackage{tkz-euclide}
\usepackage{caption}
\usepackage{wrapfig}

\usepackage[makeroom]{cancel}
\usepackage{comment}

\usetikzlibrary{decorations.pathmorphing}
\definecolor{c1}{HTML}{4C79F5}
\definecolor{c2}{HTML}{4ED6DC}
\definecolor{c3}{HTML}{5FC178}

\makeatletter
\numberwithin{equation}{section}
\numberwithin{figure}{section}
\theoremstyle{plain}
\newtheorem{thm}{\protect\theoremname}
  \theoremstyle{remark}
  \newtheorem{rem}[thm]{\protect\remarkname}
  \theoremstyle{plain}
  \newtheorem{lem}[thm]{\protect\lemmaname}
  \theoremstyle{plain}
  \newtheorem{cor}[thm]{\protect\corollaryname}
\theoremstyle{plain}
  \newtheorem{prop}[thm]{\protect\propositionname}
  \theoremstyle{remark}
  \newtheorem*{rem*}{\protect\remarkname}
  \theoremstyle{definition}
  \newtheorem{defn}[thm]{\protect\definitionname}

\makeatother

\usepackage{babel}
  \providecommand{\corollaryname}{Corollary}
  \providecommand{\propositionname}{Proposition}
  \providecommand{\definitionname}{Definition}
  \providecommand{\lemmaname}{Lemma}
  \providecommand{\remarkname}{Remark}
\providecommand{\theoremname}{Theorem}

\begin{document}
\global\long\def\k{k}
\global\long\def\bs{\boldsymbol{\sigma}}
\global\long\def\ns{\nu_{\star}}
\global\long\def\bx{\mathbf{x}}
\global\long\def\by{\mathbf{y}}
\global\long\def\grad{\nabla_{{\rm sp}}}
\global\long\def\Hess{\nabla_{{\rm sp}}^{2}}
\global\long\def\ddq{\frac{d}{dR}}
\global\long\def\qs{q_{P}}
\global\long\def\Es{E_{\star}}
\global\long\def\nh{\boldsymbol{\hat{\mathbf{n}}}}
\global\long\def\BN{\mathbb{B}^{N}}
\global\long\def\SN{\mathbb{S}^{N-1}}
\global\long\def\SNq{\mathbb{S}^{N-1}(q)}
\global\long\def\US{S^{N-1}}
\global\long\def\USN{S^{N}}
\global\long\def\SNqs{\mathbb{S}^{N-1}(\qs)}
\global\long\def\bss{\bs_{\star}}
\global\long\def\nd{\nu^{(\delta)}}
\global\long\def\nz{\nu^{(0)}}
\global\long\def\cls{c_{LS}}
\global\long\def\qls{q_{LS}}
\global\long\def\dls{\delta_{LS}}
\global\long\def\Els{E_{LS}}
\global\long\def\P{\mathbb{P}}
\global\long\def\E{\mathbb{E}}
\global\long\def\etn{\rho_{N}}
\global\long\def\betn{\bar\rho_{N}}
\global\long\def\mn{m_{N}}
\global\long\def\bmn{\bar m_{N}}
\global\long\def\CDTN{Condition~A}
\global\long\def\spp{{\rm Supp}(\mu)}
\global\long\def\Fc{F^{\rm c}_{N,\beta}}
\global\long\def\V{{\rm Vol}}

\global\long\def\F{F_{N,\beta}}
\global\long\def\FF{\widetilde{F}_{N,\beta}}

\global\long\def\sppP{{\rm Supp}(\mu_\beta)}
\global\long\def\sppPprime{{\rm Supp}(\mu_P')}
\title{Free energy landscapes in spherical spin glasses}

\author{Eliran Subag}
\begin{abstract}
We introduce and analyze free energy landscapes defined by associating to any point inside the sphere a free energy calculated on a thin spherical band around it, using many orthogonal replicas.
This allows us to prove for general spherical models several results related to the Thouless-Anderson-Palmer (TAP) approach originally introduced in the 70s for the Sherrington-Kirkpatrick model. We establish a TAP representation for the free energy, valid for any overlap value which can be sampled as many times as we wish in an appropriate sense. We call such overlaps multi-samplable. 
The correction to the Hamiltonian in the TAP representation arises in our analysis as the free energy of a certain model on an overlap dependent band.
For the largest multi-samplable overlap it coincides with the Onsager reaction term from physics. For smaller multi-samplable overlaps the formula we obtain is new. We also derive the corresponding TAP equations for critical points.
We prove the above without appealing to the celebrated Parisi formula or the ultrametricity property. 

We prove that any overlap value in the support of the Parisi measure is multi-samplable. For generic models, we further show that the set of multi-samplable overlaps coincides with a certain set that arises in the characterization for the Parisi measure by Talagrand. Moreover, using the Parisi measure of the model on the sphere we express that of the aforementioned model on the band related to TAP correction. 
The ultrametric tree of pure states can be embedded in the interior of the sphere in a natural way. For this embedding, we show that the points on the tree uniformly maximize the free energies we define. From this we conclude that the Hamiltonian at each point on the tree is approximately maximal over the sphere of same radius, and that points on the tree approximately solve the TAP equations for critical points. 
 
\end{abstract}

\maketitle
\section{Introduction}
One of the most influential papers in the theory of spin glasses is the work  of Thouless, Anderson and Palmer (TAP) \cite{TAP}, 
 where they derived the 
famous mean field equations for the local magnetizations of the Sherrington-Kirkpatrick model \cite{SK75}. To analyze the thermodynamics, they proposed, one should investigate the solutions of those equations, to which they associated a certain free energy $F_{\rm TAP}(m)$, computed from a diagrammatic expansion of the partition function. 
Their approach was further developed in physics, see e.g. \cite{Bray1980,Cavagna2003,CrisantiSommersTAPpspin,DeDominicis1983,Gross1984,KurchanParisiVirasoro,Plefka},  with the general idea that the free energy $F_{N,\beta}$ associated to a mean-field Hamiltonian $H_N(\bs)$ (defined below in \eqref{eq:Hamiltonian}, \eqref{eq:F}) is approximated by
\begin{equation}
F_{N,\beta}\approx \max_{m}\Big\{ 
\frac{\beta}{N} H_N(m) + F_{\rm TAP}(m)
\Big\}, \label{eq:TAPoriginal}
\end{equation}
where the maximum is over solutions of the aforementioned equations which satisfy certain convergence conditions \cite{Plefka}.
It was also claimed that, in principle, one has to account for the multiplicity of solutions with the same free energy known as `complexity' \cite{TAP-SK4,Cavagna1998,Crisanti2003,Crisanti2005,Rieger1992a}.
If the complexity for the maximizer is zero, then the maximum can be taken over all points where the convergence conditions hold, and not only the solutions.
As observed already in the original paper  of Thouless, Anderson and Palmer \cite{TAP}, viewing the free energy they defined as a function on the whole space of possible magnetizations --- that is, the convex hull of the configuration space --- the solutions to their equations arise as critical points. One can therefore think of a \emph{free energy landscape} on the space of magnetizations, (some of) the maximizers of which should correspond to the so-called `states' of the Gibbs measure. 

In this paper we introduce an alternative, geometric approach to define and analyze free energy landscapes for spherical mixed $p$-spin  models,
which allows us to prove some of the main ideas from the TAP approach for the spherical models in the low temperature phase, and establish clear connections with the replica symmetry breaking theory of Parisi \cite{ParisiFormula,Parisi}. 
Our main results can be applied to any overlap value $q\in[0,1]$ which can be sampled as many times as we wish in an appropriate sense. We call such overlaps multi-samplable (Definition \ref{def:multi}). For each multi-samplable overlap we will identify the near-maximizers of the free energies we define as the near-maximizers of the Hamiltonian on a sphere of appropriate radius (Theorem \ref{thm:eqiuv} and \eqref{eq:equivalences}), derive a generalized TAP representation for the free energy (Theorem \ref{thm:FE}), and write the TAP equations for critical points (Section \ref{subsec:TAPeq}).
For any given multi-samplable overlap $q$, our representation is of the form in \eqref{eq:TAPoriginal} with the maximum taken over the sphere scaled by $\sqrt q$. Moreover, for the maximal such overlap the correction added to the Hamiltonian in our representation  coincides with the well-known Onsager reaction term from the physics literature  \cite{CrisantiSommersTAPpspin,KurchanParisiVirasoro} and the spherical analogue of the original TAP representation is recovered (Corollary \ref{cor:classicalTAP}). In general, however, we show that the correction in our representation corresponds to a free energy of another spherical model that depends on $q$ and can be different from the classical Onsager correction. 

The TAP representation was previously proved for some spherical models  \cite{BeliusKistlerTAP,geometryMixed,geometryGibbs} and for general models with Ising spins \cite{ChenPanchenkoTAP}, and the TAP equations were proved (in a different formulation than ours) for the Sherrington-Kirkpatrick model in high temperature phase \cite{BolthausenTAP,ChatterjeeTAP,TalagrandBookI} and for generic models with Ising spins in the low temperature phase \cite{AuffingerJagannthSpinDist,AuffingerJagannathTAP}. All the results cited here correspond to in our analysis to choosing the multi-samplable overlap value $q$ to be equal to the maximal point in the support of the Parisi measure. We will put those results into context and give more details below, after the precise statement of our representation for the free energy and the corresponding critical point equations.

All our results mentioned so far will follow mostly from first principles and concentration properties. In particular, their proofs are independent of Parisi's formula and the ultrametricity property. They are, however, closely related to both. First, any overlap in the support of the Parisi measure is multi-samplable, (Corollary \ref{cor:ParisiSupp}) and therefore all the aforementioned results apply to such overlaps. In fact, we shall see that the set of multi-samplable overlaps contains a set which in principle may be larger than the support --- defined through an optimality condition for the Parisi measure which was derived by Talagrand \cite{Talag}. Further,  for generic models this set is equal to the set of multi-samplable overlaps (Theorem \ref{thm:qsamp}). As mentioned above, we will express the correction in the TAP representation as the free energy of a certain spherical model. For multi-samplable $q$, we shall see that the Parisi distribution of this model coincides, up to scaling, with the restriction to $[q,1]$  of the Parisi distribution of the original model on the sphere. Moreover, the minimizer in the zero-temperature analogue of the Parisi distribution for the ground state energy on the sphere of radius $\sqrt{Nq}$ can be expressed using the restriction of the Parisi distribution of the original model on the sphere to $[0,q]$ (Proposition \ref{prop:Parisidistributions}).

It is well-known that for generic models the Gibbs measure splits into pure states that are arranged ultrametrically, i.e., in a tree structure according to the overlaps of their barycenters. See Section \ref{subsec:PS} for the precise statement.
This tree can be naturally embedded in the Euclidean ball of radius $\sqrt N$ (Corollary \ref{cor:embed}). We will see that, by construction, any point $m$ of the embedded tree maximizes the free energies we define. We will conclude from this that  the value of the Hamiltonian at  any point $m$ from the tree is approximately maximal over the sphere of radius $\|m\|$, uniformly over the tree (Corollary \ref{cor:Eontree}). And that each such $m$ approximately solves the critical point TAP equations (Corollary \ref{cor:VTAPsol}).

Our results can be used to relate the overlap distribution for independent samples from the Gibbs measure at two different temperatures to overlaps between near-maximizers of the Hamiltonian on spheres of two radii, and this can be used to deduce certain results on temperature chaos, see the arXiv version of this paper.

Finally, we remark that while the current paper focuses on spherical models, the main ideas we introduce here are more general and can be also applied in the setting of mixed models with Ising spins. This was recently done by Chen, Panchenko and the author in \cite{TAPChenPanchenkoSubag}, where analogues of many results from the current paper were derived. See also \cite{TAPIIChenPanchenkoSubag} for a zero-temperature analogue in the Ising case. We wish to emphasize, however, that there is a crucial difference in the Ising case:
the geometry of the band truly depends on the magnetization vector $m$, instead of only $\|m\|$ in the spherical case. This made the analysis of \cite{TAPChenPanchenkoSubag} and the computation of the TAP correction in particular, which in this case is a function of the empirical measure $\mu_m:=\frac1N \sum_{i\leq N}\delta_{|m_i|}$,
highly non-trivial and required new ideas and techniques.

\subsection{\label{subsec:1.1}The model and free energy landscapes}
Given a sequence $\gamma_{p}\geq 0$
 decaying exponentially fast 
in $p$, the \emph{Hamiltonian} of the spherical mixed $p$-spin model with
\emph{mixture} $\nu(x)=\sum_{p=1}^{\infty}\gamma_{p}^{2}x^{p}$ is the Gaussian process on $\SN$, the (Euclidean) sphere of radius $\sqrt N$ in dimension $N$, defined by 
\begin{equation}
H_{N}\left(\bs\right)=H_{N,\nu}\left(\bs\right):=\sum_{p=1}^{\infty}\frac{\gamma_{p}}{N^{\left(p-1\right)/2}}\sum_{i_{1},...,i_{p}=1}^{N}J_{i_{1},...,i_{p}}^{(p)}\sigma_{i_{1}}\cdots\sigma_{i_{p}},\label{eq:Hamiltonian}
\end{equation}
where $J_{i_{1},...,i_{p}}^{(p)}$ are i.i.d. real
standard Gaussian random variables and $\bs=(\sigma_{1},...,\sigma_{N})$. The covariance function of $H_N(\bs)$ is
\[
\E\big\{H_N(\bs) H_N(\bs') \big\} = N\nu(R(\bs,\bs')),
\]
where
$R(\bs,\bs'):=\frac1N\bs\cdot\bs'$ is the usual overlap function.  
The associated Gibbs measure on $\SN$ at inverse-temperature $\beta>0$  is  
\begin{equation}
\frac{dG_{N,\beta}}{d\bs}\left(\boldsymbol{\sigma}\right):=\frac{1}{Z_{N,\beta}}e^{\beta H_{N}\left(\boldsymbol{\sigma}\right)}, \label{eq:Gibbs}
\end{equation}
where  $d\bs$ denotes the normalized Haar measure and the partition function is defined as
\begin{equation*}
	Z_{N,\beta}:=\int_{\SN}e^{\beta H_{N}\left(\boldsymbol{\sigma}\right)}d\bs.
\end{equation*} 

The free energy is defined by
\begin{equation}
	\label{eq:F}
	F_{N,\beta}:=\frac1N\log Z_{N,\beta}.
\end{equation}
Its limit is given by the famous Parisi formula \cite{ParisiFormula,Parisi}, or its spherical version discovered by Crisanti and Sommers \cite{pSPSG}, which was first proved  by Talagrand \cite{Talag,Talag2} for models with even $p$-spin interactions, following the discovery of the interpolation scheme of Guerra \cite{GuerraBound}, and later generalized to spherical models with arbitrary interactions by Chen \cite{Chen}, based on the Aizenman-Sims-Starr scheme \cite{ASSscheme,AizenmanSimsStarrROST} and ultrametricity property proved by Panchenko \cite{ultramet}.
The formula states that 
\begin{equation}
\lim_{N\to\infty} \E F_{N,\beta} = \inf_{x\in D([0,1])} \mathcal P_\beta (x,\nu), 
\label{eq:parisi}
\end{equation} 
where $D([0,1])$ is the set of distribution functions $x$ on  $[0,1]$, i.e., non-decreasing right continuous $[0,1]\to[0,1]$ functions,  such that $x(\hat q)=1$ for some $\hat q<1$, and where
\begin{equation*}
	\mathcal P_\beta(x,\nu) =\frac{1}{2}\Big( 
	\beta^2\int_0^1x(q)\nu'(q)dq+\int_0^{\hat q}\frac{dq}{\int_q^1x(s)ds}+\log(1-\hat q) 
	\Big).
\end{equation*}  
Since $\mathcal P_\beta$ is strictly convex in $x$, there is a unique minimizer $x_\beta$ to \eqref{eq:parisi},  called the Parisi distribution. The corresponding Parisi measure will be denoted $\mu_\beta$. We remark that the existence of the limit of the free energy as in the left-hand side of \eqref{eq:parisi}, which will appear in our main results, was recently proved in \cite{FreeEnergyConvergence} without relying on the Parisi formula.

The space of magnetizations, on which we define the free energy landscapes, is the ball $\BN:=\{m:\,\|m\|<\sqrt N\}$. On $\BN$ we use the same formula \eqref{eq:Hamiltonian} to define $H_N(m)$. 
The first free energy we consider is defined by associating to each point $m$ a thin spherical band around it
\begin{equation}
{\rm Band}\left(m,\delta \right):=\left\{ \bs\in\SN:\,| R(\bs-m,m
)
|\leq\delta\|m\|/\sqrt N \right\} ,\label{eq:band}
\end{equation}
and  computing the free energy of the Hamiltonian centered by $H_N(m)$ over the band,
\begin{equation}
\F(m,\delta) := \frac1N\log
\int_{{\rm Band}(m , \delta)}e^{\beta (H_{N}(\bs)-H_N(m))}d\bs.
\label{eq:FEband}
\end{equation}

\begin{wrapfigure}{r}{5.2cm}
	\caption{\label{fig:1}Bands around points with $\|m\|=\sqrt{Nq}$. The two vectors in blue end at two points with $R(\bs_1,\bs_2)=R(m,m)$.
	}
	\begin{tikzpicture}[scale=0.4,domain=0:1]  	
	\node[inner sep=0pt] (1st) at (0,0)
	{\includegraphics[width=4cm]{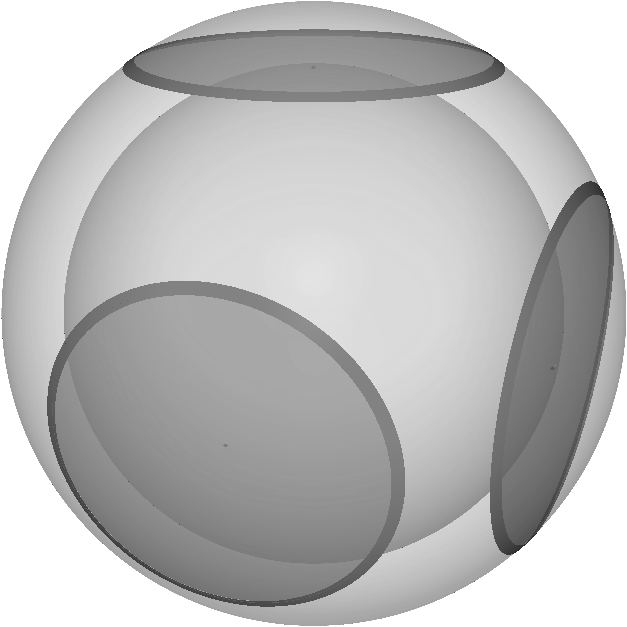}};
	\draw[color=black, line width = .3mm] (-.00,-.0) circle (5cm);
	\draw [-stealth, color=black, line width = .5mm,line cap=round] (0,0) -- (3.9,3.1);
	\draw [-stealth, color=black, line width = .5mm,line cap=round] (0,0) -- (-3.1,2.5);
	\node[below] at (2.15,1.6) {$\sqrt N$};
	\node[below] at (-3.1,1.6) {$\sqrt{Nq}$};
	\draw [-stealth, color=RoyalBlue, line width = .4mm,line cap=round] (-1.41,-2.1) -- (1.25,-2.1);
	\draw [-, color=RoyalBlue, line width = .4mm,line cap=round] (-.88,-2.1+0.56) -- (-.72,-2.1);
	\draw [-, color=RoyalBlue, line width = .4mm,line cap=round] (-1.51,-2.1+0.56) -- (-.88,-2.1+0.56);
	\draw [-stealth, color=RoyalBlue, line width = .4mm,line cap=round] (-1.41,-2.1) -- (-2.1,.4);
	\tkzDrawPoint[color=black, fill=orange, size=4](3.81,-.87) 
	\tkzDrawPoint[color=black, fill=orange, size=4](-.005,3.935) 
	\tkzDrawPoint[color=black, fill=orange, size=4](-1.41,-2.1) 
	\end{tikzpicture}
\end{wrapfigure}

\noindent Fixing some $\delta_N>0$ slowly going to $0$, we will denote $F_{N,\beta}(m):=F_{N,\beta}(m,\delta_N)$.
Note that since the thin band is approximately a sphere of co-dimension $1$, $\F(m)$ is the free energy of a spherical model with mixture that depends only on $q=\|m\|^2/N$ (see \eqref{eq:nuq_tilde}), up to a vanishing  in $N$ correction.

We may think of $\frac{\beta}{N}H_N(m)+\F(m)$ as a free energy landscape on $\BN$ which, loosely speaking, allows us to `scan'  regions of the sphere $\SN$ of different location and scale as we vary $m$.
We wish to understand the set of points whose bands have maximal free energy at  exponential scale,
\begin{equation}
\label{eq:Fapx}
\frac1N\log
\int_{{\rm Band}(m , \delta_N)}e^{\beta H_{N}(\bs)}d\bs =\frac{\beta}{N}H_N(m)+\F(m)\approx F_{N,\beta}.
\end{equation}

One of the main difficulties is that, although for any given non-random point  $\F(m)$ concentrates, the maximal fluctuations $\sup_{m\in \BN}|\F(m)-\E \F(m)|$ are typically of order $O(1)$,\footnote{We may write a Taylor expansion for $H_N(\bs)-H_N(m)$ around $m$. Its linear term, which acts in a sense as a `magnetic field' on the band, is $\langle\nabla H_N(m),\bs-m\rangle$. One expects that typically for points (outside of a neighborhood of the origin) where $\big|\frac{\|\nabla H_N(m)\|}{\E\|\nabla H_N(m)\|}-1\big|$ is of order $O(1)$, so is the deviation $|\F(m)-\E \F(m)|$. See more on this point in Footnote \ref{ft:maxfluct} below.
} so that $\F(m)$ is genuinely a random process. Another undesired property of $\F(m)$  is that the set of points where \eqref{eq:Fapx} holds is too large in a sense to allow one to gain any meaningful information.\footnote{To see that, imagine for a moment that instead of working with bands we work with spherical caps. If some small region has large Gibbs weight, then any cap containing it also heavy. In high dimensions, however, most of the volume of the cap concentrates close to its boundary, on a thin band.}
To deal with both issues simultaneously, we will define another free energy $\FF(m)$ with \emph{stronger concentration property} such that for any $m$,
\begin{equation}
\label{eq:ineq}
\frac{\beta}{N}H_N(m)+ \FF(m) \leq \frac{\beta}{N}H_N(m)+ \F (m) \leq F_{N,\beta},
\end{equation}
and for some special set of points both inequalities become \emph{approximate equalities.}
We define $\FF(m)$ by computing the free energy on the band using many replicas, constrained to be roughly orthogonal relative to $m$. Precisely, for $n$ replicas and overlap constraint $\rho>0$, we define
\begin{equation}
\label{eq:Ft}
\begin{aligned}
&\FF(m,n,\rho,\delta):  =\\&\quad\frac{1}{ Nn}\log\int_{B_{N}(m,n,\rho,\delta)}\exp\Big\{\beta\sum_{i=1}^{n}(H_{N}(\bs_{i})-H_{N}(m))\Big\}d\bs_{1}\cdots d\bs_{n},
\end{aligned}
\end{equation}
\begin{equation*}
\begin{aligned}
	&B_{N}(m,n,\rho,\delta): =\\
	&\quad\big\{(\bs_{1},\ldots,\bs_{n})\in({\rm Band}(m,\delta ))^{n}:\,|R(\bs_{i},\bs_{j})-R(m,m)|<\rho,\,\forall i\neq j\big\}.
\end{aligned}
\end{equation*}
Note that on a thin band ${\rm Band}(m,\delta_N)$,
\begin{equation}
\label{eq:orth}
R(\bs_i,\bs_j)\approx R(m,m) \iff R(\bs_i-m,\bs_j-m)\approx 0.
\end{equation}

Fix some sequences $n_N$ and $\rho_N$ going to $\infty$ and $0$ slowly, 
 and define
\begin{equation}
	\label{eq:Fhatbs}
	\FF(m):=\FF(m,n_N,\rho_N,\delta_N).
\end{equation}
Note that we can think of $\FF(m)$  as $\F(m)$ plus a penalty term for the conditional Gibbs probability of sampling many orthogonal points in the sense of \eqref{eq:orth}, 
\begin{equation}
\label{eq:1505-01}\begin{aligned}
&\FF(m)=\F(m)
\\
&+\frac1{ Nn_N}\log G_{N,\beta}^{\otimes n_N}\Big\{  \big|R( \bs_i,\bs_j) - R(m,m)   \big| < \rho_N,\, \forall i\neq j \,\Big| \, \bs_i\in  {\rm Band}(m, \delta_N) \Big\}.
\end{aligned}
\end{equation}
Also observe that $\FF(m,n,\rho,\delta)$ is the free energy corresponding  to the process $\sum_{i=1}^{n}(H_{N}(\bs_{i})-H_{N}(m))$ on $B_{N}(m,n,\rho,\delta)$ times $\frac1n$. As a consequence of the orthogonality \eqref{eq:orth}, uniformly over $B_{N}(m,n,\rho,\delta)$  the variance of this process is $O\big(N(n+n^2(\delta+\rho))\big)$. The $\frac1n$ factor which multiplies the free energy, compared to the $\sqrt n$ which comes from the standard deviation for large $N$ and small $\rho$ and $\delta$, will allow us to obtain the aforementioned concentration property of $\FF(m)$, which will play a crucial role in our analysis. We will prove the following proposition in Section \ref{sec:UnifConc}.
\begin{prop}[Uniform concentration]
	\label{prop:concExpec} Let $\beta>0$ and assume $\delta_N$, $\rho_N$ and $n_N$ tend to their limits slowly enough.\footnote{To be precise, here and in all other results in the paper, we mean that there exist some $\delta_N^0, \rho_N^0\to0$ and $n_N^0\to\infty$ which may depend on $\nu$ such that the statement holds for any $\delta_N\geq\delta_N^0, \rho_N\geq\rho_N^0$ and $n_N\leq n_N^0$ such that  $\delta_N, \rho_N\to0$ and $n_N\to\infty$.
	} Then, for any $t,c>0$, for large enough $N$,
	\begin{equation}
	\P\bigg(\,  \max_{\|m\|^2<N}\big|\, \FF(m)-\E \FF(m)\,\big|<t\, \bigg)>1-e^{-Nc}.\label{eq:concExpec}
	\end{equation}
\end{prop}

One of the motivations for the definition of $\FF(m)$ is that, as claimed above, we are guaranteed to have points with $\|m\|^2=Nq$ such that \eqref{eq:ineq} holds with approximate equalities instead of inequalities, provided that $q$ can be sampled as many times as we wish in the following sense.\footnote{In \eqref{eq:goodq}, the $\lim$ can be weakened to $\limsup$. More precisely, if $q$ satisfies \eqref{eq:goodq} with $\limsup$ for any $n$ and $\rho$, then it also satisfies it with $\lim$. See Remark \ref{rem:multi_limsup}.}
\begin{defn}\label{def:multi}
	Given $\beta>0$, an overlap value $q\in[0,1)$ is \emph{multi-samplable} 
	 if for any $n\geq1$ and $\rho>0$,
	\begin{equation}\label{eq:goodq}
		\lim_{N\to\infty}\frac1N\log \E G_{N,\beta}^{\otimes n}\Big(\,
		\big| R(\bs_i,\bs_j)-q\hspace{.025cm} \big|<\rho,\,\forall i\neq j\,
		\Big)=0.
	\end{equation}
\end{defn}

Note that \eqref{eq:goodq} is a lower bound on the Gibbs probability, since it is almost surely in $(0,1)$. 
We will see that if $q$ is multi-samplable, then w.h.p. there exists some point $m_\star$ such that $\|m_\star\|^2=Nq$ and
\[\frac{\beta}{N}H_N(m_\star)+ \FF(m_\star) = F_{N,\beta}+o(1).\] 
From the well-known concentration of $F_{N,\beta}$ (see e.g. \cite[Theorem 1.2]{PanchenkoBook})  and that of $\FF(m)$ in Proposition \ref{prop:concExpec}, 
\[\frac{\beta}{N}H_N(m_\star)+ \E\FF(m_\star) = \E F_{N,\beta}+o(1).\]
From concentration and \eqref{eq:ineq}, for \emph{any} $m=m_N$ with  $\|m\|^2=Nq$, the same holds with an inequality in the direction $\leq$.
Hence, the energy at $m_\star$ must be close to the ground-state energy\footnote{The existence of the limit in expectation can be deduced from the convergence of the free energy and the a.s. convergence  follows from the Borell-TIS inequality \cite{Borell,TIS}.} 
\begin{equation}
\Es(q):=\lim_{N\to\infty}\frac{1}{N}\E\max_{\|m\|^2=Nq}H_{N}(m)\stackrel{{\rm a.s.}}{=}\lim_{N\to\infty}\frac{1}{N}\max_{\|m\|^2=Nq}H_{N}(m),\label{eq:GS}
\end{equation}
using which we will prove the following theorem in Section \ref{sec:pfThm2}.

\begin{thm}\label{thm:eqiuv}Let $\beta>0$ and assume $\delta_N$, $\rho_N$ and $n_N$ tend to their limits slowly enough. Then $q\in[0,1)$ is multi-samplable if and only if 
	\begin{equation}\label{eq:E+F}
	\beta \Es(q)+ \lim_{N\to\infty}\E \FF(m) = \lim_{N\to\infty} \E F_{N,\beta},
	\end{equation}
	where $m=m_N$ is an arbitrary point such that $\|m\|^2=Nq$. 
\end{thm}

From the concentration of $F_{N,\beta}$ and of $\FF(m)$ in Proposition \ref{prop:concExpec}, if $q$ is multi-samplable, \eqref{eq:E+F} implies that
for any $t>0$, for some $c>0$ and large enough $N$,
\begin{equation}
\label{eq:temp}
\begin{aligned}
	\P\bigg(\sup_{\|m\|^2=Nq}\big|  \Delta_{N,\beta}(m)  \big| >t \bigg)<e^{-Nc},
	\end{aligned}
\end{equation}
where
\[
\Delta_{N,\beta}(m):=\Big|\frac{\beta}{N}H_N(m)+\FF(m) - F_{N,\beta} \Big|  - \Big| \frac{\beta}N H_N(m) - \beta\Es(q)\Big|.
\]
In particular, using \eqref{eq:ineq}, for any multi-samplable $q$ w.h.p. on $\sqrt q\cdot\SN = \big\{m: \|m\|^2=Nq\big\}$,
\begin{equation}
\label{eq:equivalences}
\frac{\beta}{N}H_N(m)+\FF(m)\approx
\frac{\beta}{N}H_N(m)+\F(m) \approx F_{N,\beta} \iff \frac1N H_N(m)\approx \Es(q).
\end{equation}
In words, the bands which approximately have the same free energy (in the sense of \eqref{eq:Fapx}) as the total free energy $F_{N,\beta}$ and for which the probability of sampling many orthogonal (as in \eqref{eq:orth}) points under the Gibbs measure is not exponentially small (see \eqref{eq:1505-01}) are exactly the bands around points which approximately maximize the energy $H_N(m)$ on $\sqrt q\cdot\SN$.

\subsection{\label{subsec:TAP}The TAP approach}
The thin band ${\rm Band}\left(m,\delta_N \right)$ is approximately a sphere in one dimension less than the original model, and we may roughly think of the restriction of $H_N(\bs)-H_N(m)$ to it as another spherical model. More precisely, given some $m$ with $\|m\|^2=Nq$, consider the transformation $\tilde \bs = (\bs-m)/\sqrt{1-q}$ and the process $\tilde H_N(\tilde \bs):=H_{N}(\bs)-H_N(m)$. Then uniformly in $\bs_1, \bs_2\in {\rm Band}\left(m,\delta_N \right)$, one can verify that
\begin{equation}
\label{eq:Htilde}
\frac1N\E\Big\{ \tilde H_{N}(\tilde\bs_{1}) \tilde H_{N}(\tilde \bs_{2})\Big\} =  \tilde\nu_{q}(R(\tilde \bs_1,\tilde \bs_2))+o(1),
\end{equation}
as $N\to\infty$, where we define the mixture
\begin{equation}
	\tilde\nu_{q}(t):=
	\nu((1-q)t+q)-\nu(q)
	=\sum_{k=1}^{\infty}a_{k}^2(q)t^{k}\label{eq:nuq_tilde}
\end{equation}
with $a_k^2(q):=\sum_{p\geq k}\gamma_p^2\binom{p}{k}(1-q)^{k}q^{p-k}$.
We also define the mixture 
\begin{equation}
	\nu_{q}(t):=
	\nu((1-q)t+q)-\nu(q)-\nu'(q)(1-q)t
	=\sum_{k=2}^{\infty}a_{k}^2(q)t^{k},
	\label{eq:nuq}
\end{equation}
obtained from $\tilde\nu_{q}(t)$ by removing the $1$-spin interaction. Finally, with $H_N^q(\bs)$ denoting the Hamiltonian corresponding to $\nu_q$, we define 
\[
F_\beta(q):=\lim_{N\to\infty}\frac1N\E \log\int_{\SN}e^{\beta H^q_{N}\left(\bs\right)}d\bs.
\]
\subsubsection{TAP representation for the free energy}
We prove the following proposition in Section \ref{sec:TAPcorrection}.

\begin{prop}
	\label{prop:FFlim}
	If $\delta_N, \rho_N\to0$ and $n_N\to\infty$ sufficiently slowly, then for any $t\in(0,1)$ and arbitrary $\bs=\bs_N\in\SN$,\footnote{\label{ft:maxfluct}One can also show that $\lim_{N\to\infty} \E\F(m)=\frac{1}{2}\log(1-q)+\bar F_{\beta}(q)$ where $\bar F_{\beta}(q)$ is defined similarly to $F_{\beta}(q)$, but using $\tilde \nu_q $ instead of $\nu_q$, and therefore
		$\E\FF(m)<\E\F(m)$ for any $m$ outside of a neighborhood of the origin, with the difference being bounded away from zero uniformly in $N$ large enough. For points as in \eqref{eq:equivalences}, $\F(m)\approx \FF(m)\approx \E\FF(m)$. In particular, for such points $\E \F(m)-\F(m)$ is of order $O(1)$.}
	\begin{equation}
		\label{eq:EFt}
		\lim_{N\to\infty} \sup_{q\in[0,t]}\Big| \E\FF(\sqrt q \bs) - \frac12\log(1-q)-F_\beta(q)\Big|=0.
	\end{equation}
\end{prop}
We therefore have the following.

\begin{thm}[TAP representation]
	\label{thm:FE} For any $\beta>0$,
	\begin{equation}
	\label{eq:TAP}
	\lim_{N\to\infty}\E F_{N,\beta} =  \beta\Es(q)+\frac{1}{2}\log(1-q)+F_{\beta}(q)
	\end{equation}	
		if and only if $q\in[0,1)$ is multi-samplable. Otherwise, \eqref{eq:TAP} holds with a strict inequality $>$.
\end{thm}
\begin{proof}
Let $q\in[0,1)$. Let $\delta_N$, $\rho_N$ and $n_N$ be some sequences such that the conclusion of Theorem \ref{thm:eqiuv} and Proposition \ref{prop:FFlim} hold (with $t>q$).
Then by substituting in \eqref{eq:E+F} the limit of $\E \FF(m)$ from \eqref{eq:EFt} it folows that \eqref{eq:TAP} holds if and only if $q$ is multi-samplable.

Let $q\in[0,1)$ be arbitrary and not necessarily multi-samplable. Let $m_\star=m_{N,\star}$ be the maximizer of $H_N(m)$ over the sphere $\sqrt{q}\SN$. Then, from \eqref{eq:ineq} and Proposition \ref{prop:concExpec}, with probability going to $1$,
\begin{equation*}
	\frac{\beta}{N}H_N(m_\star)+ \E\FF(\sqrt{q}\bs)  \leq F_{N,\beta},
\end{equation*}
where $\bs$ is an arbitrary (non-random) point in $\SN$. By the Borell-TIS inequality \cite{Borell,TIS} and concentration of the free energy \cite[Theorem 1.2]{PanchenkoBook}, the inequality also holds in expectation, up to a $o(1)$ error in the $N\to\infty$ limit. 
By Proposition \ref{prop:FFlim} it follows that \eqref{eq:TAP} holds with an inequality $\geq$.
\end{proof}

For the largest multi-samplable\footnote{Since the logarithmic term in \eqref{eq:TAP} tends to $-\infty$ as $q\to1$, it is easy to verify that for fixed $\beta$ values of $q$ close enough to $1$ are not multi-samplable. Thus, from Definition \ref{def:multi}, the set of multi-samplable overlaps is closed.}
\begin{equation}
\label{eq:qmax}
q_{\max}=q_{\max}(\beta):=\max\big\{q\in[0,1):\,q\mbox{ is multi-samplable} \big\},
\end{equation}  
we shall see that $\nu_q$ corresponds to a replica symmetric model and therefore $F_{\beta}(q)$ coincides with the well-known Onsager reaction term from the physics literature  \cite{CrisantiSommersTAPpspin,KurchanParisiVirasoro} as in \eqref{eq:RSTAP}. Note that \eqref{eq:RSTAP} holds if and only if the Hamiltonian $H_N^q(\bs)$ with mixture  $\nu_q$ exhibits replica symmetry at $\beta$ (i.e., the minimizer in the corresponding Parisi formula is $\delta_0$). 
\begin{cor}\label{cor:classicalTAP}
			 For $q=q_{\max}$,
			 \begin{equation} 
			 \label{eq:RSTAP}
			 F_{\beta}(q)=\frac12 \beta^2 \nu_q(1)
			 =\frac12 \beta^2 \Big(\nu(1)-\nu(q)-(1-q)\nu'(q)\Big).
			 \end{equation}
\end{cor}
We prove the corollary in Section \ref{sec:TAPpfs}. 
We will also prove that any overlap value $q$ in the support of the Parisi measure is multi-samplable and therefore satisfies the representation for the free energy above, and that the maximal point in the support, which we denote by $q_P$, satisfies \eqref{eq:RSTAP} as well  (see Section \ref{subsec:Parisi}). We note that while it seems natural to expect that typically $q_P=q_{\max}$, this is not always the case (see Remark \ref{rem:Sbeta}).

By the representation \eqref{eq:TAP}, as in \eqref{eq:TAPoriginal}, asymptotically the free energy is equal to the maximum of the normalized energy plus the deterministic function 
\[
	F_{\rm TAP}(m)=\frac{1}{2}\log(1-\|m\|^2/N)+F_{\beta}(\|m\|^2/N)
\] over the sphere of radius $\sqrt{qN}$, whenever $q$ is multi-samplable. For models with Ising spins, Chen and Panchenko \cite{ChenPanchenkoTAP} proved a TAP representation with the classical Onsgaer correction and overlap $q=q_P$, the rightmost point in the support of the Parisi measure. Namely, they showed that the free energy is given by the maximum of the normalized energy plus a certain deterministic function over the intersection of the sphere of radius $\sqrt{q_PN}$ and the full cube $[-1,1]^N$. They also proved that the representation still holds if one maximizes over points with norm larger than $\sqrt{q_PN}$ in the cube. The deterministic function added to the energy in their representation is the Ising analogue of $F_{\rm TAP}(m)$ above, if one substitutes the simplified expression from \eqref{eq:RSTAP}. In contrast to the spherical case, in the Ising case this function depends on the empirical measure $\mu_m:=\frac1N \sum_{i\leq N}\delta_{|m_i|}$, and not only $\|m\|$. 

In \cite{TAPChenPanchenkoSubag} Chen, Panchenko and the author prove analogues of several results from the current work for general mixed models with Ising spins and in particular of the TAP representation of Theorem \ref{thm:FE}. The latter applies to any $q$ in the support of the Parisi measure and not only $q_P$.

Coming back the spherical models, the TAP representation for the free energy was previously proved in some specific models only at $q=q_P$ and with the correction \eqref{eq:RSTAP}. For the 1-RSB pure $p$-spin models with $p\geq3$ and $\beta\gg 1$ it was proved by the author \cite{geometryGibbs}, for 1-RSB mixed models close to pure and $\beta\gg 1$ by Ben Arous, Zeitouni and the author \cite{geometryMixed}, and for the replica symmetric pure $2$-spin model with general $\beta>0$ and external field by Belius and Kistler \cite{BeliusKistlerTAP}. Both \cite{geometryGibbs} and \cite{geometryMixed} heavily build on the investigation of critical points by Auffinger, Ben Arous and {\v{C}}ern{\'y} \cite{A-BA-C}, Auffinger and Ben Arous \cite{ABA2}, the author \cite{2nd} and Zeitouni and the author \cite{pspinext}.

The logarithmic term in \eqref{eq:EFt} is easy to explain. It is an entropy term, equal to the limit of $\frac1N\log$ of the volume of ${\rm Band}(m,\delta_N)$. 
We will arrive at the
second term in \eqref{eq:EFt} as follows.
By a continuity argument, we will first relate $\E \FF(m)$ to a similar free energy on the sphere  
${\rm Band}(m,0)$ (of co-dimension $1$ in $\SN$). After mapping it to the sphere of radius $\sqrt N$, we will need to calculate a replicated free energy of the form
\begin{equation}
\label{eq:repF}
\frac{1}{ Nn}\E\log\int_{B_N(0,n,\rho,0)
}
e^{\beta\sum_{i=1}^{n}\tilde H^q_{N}(\bs_{i})}d\bs_{1}\cdots d\bs_{n},
\end{equation}
where $B_N(0,n,\rho,0)=
\{(\bs_1,\ldots,\bs_n)\in (\SN)^n:\, |R(\bs_i,\bs_j)|<\rho	\}$
and $\tilde H^q_{N}(\bs)$ is the Hamiltonian that corresponds to the mixture $\tilde\nu_q$. This Hamiltonian can be decomposed as $H^q_{N,\beta}(\bs)$, the Hamiltonian of $\nu_q$, plus an independent $1$-spin interaction term $\alpha_1(q) \sum_{i=1}^N J^{(1)}_i \sigma_i$.
The $1$-spin cannot simultaneously have significant contribution to the value of many orthogonal points (see Lemma \ref{lem:rem1sp}).
Thus, when we let $N\to\infty$, $n\to\infty$ and $\rho\to0$, we shall see that replacing $\tilde H^q_{N}(\bs)$ by $ H^q_{N}(\bs)$ in \eqref{eq:repF}  does not affect the limiting value.
Moreover, we will prove that if a Hamiltonian does not contain a $1$-spin interaction term, then the limiting replicated free energy as in \eqref{eq:repF} is equal to its usual free energy (see Lemma \ref{lem:0ismulti}). Namely, for $H^q_{N}(\bs)$ the limiting replicated free energy as in \eqref{eq:repF} is equal to  $F_\beta(q)$. 

\subsubsection{\label{subsec:TAPeq}TAP equations}
For large enough $s\in(0,1)$, it is easy to verify that uniformly in $m$ with $\|m\|^2\geq Ns$,  $\FF(m)$ is as negative as we wish w.h.p., since the volume of the corresponding band is small (see, e.g., the proof of Proposition \ref{prop:concExpec}). 
For any model, some $q\in[0,1)$ is multi-samplable, see Remark \ref{rem:exists_q}. Thus, from the inequality \eqref{eq:ineq}, the concentration of $F_{N,\beta}$ and of  $\FF(m)$ in Proposition \ref{prop:concExpec} and the uniform convergence  of $\E\FF(m)$ as in Proposition \ref{prop:FFlim}, with $q=\|m\|^2/N$,
\begin{equation}
\label{eq:limsupTAP}
\lim_{N\to\infty}\max_{\|m\|^2<N} \left(\frac{\beta}{N} H_N(m) +\frac{1}{2}\log(1-q)+F_{\beta}(q)\right) =  \lim_{N\to\infty} \E F_{N,\beta},
\end{equation}
in probability, and  w.h.p. the points $m$  for which 
\begin{equation}
\frac{\beta}{N}H_N(m)+\FF(m)\approx
\frac{\beta}{N}H_N(m)+\F(m) \approx F_{N,\beta},
\label{eq:Fapx2}
\end{equation}
are approximate global maximizers of 
\begin{equation}
\label{eq:fTAP}
\frac{\beta}{N} H_N(m) +\frac{1}{2}\log(1-q)+F_{\beta}(q).
\end{equation}

The TAP equations are critical point equations for the TAP free energy which is defined in physics \cite{CrisantiSommersTAPpspin,KurchanParisiVirasoro} (for the pure models) similarly to \eqref{eq:fTAP}, but using the Onsager correction \eqref{eq:RSTAP}  \emph{for all} $m$. For $q$ such that $\nu_q$ is not replica symmetric, neither the free energy defined using \eqref{eq:RSTAP} nor the corresponding  TAP equations have a clear meaning. Certain stability conditions, like Plefka's  condition \cite{Plefka}, were proposed in physics to identify the `physical' solutions of the TAP  equations and presumably get rid of ones which do not correspond to actual properties of the system. 

Using the general correction $F_{\beta}(q)$ we may write critical point equations which are meaningful for all $m$: from concentration, the critical points of \eqref{eq:fTAP} with largest values are also approximate maximizers of  $\frac{\beta}{N}H_N(m)+\FF(m)$,
and thus satisfy \eqref{eq:Fapx2}. 
Note that here we do not need to restrict to solutions which satisfy a condition like Plefka's. If one only restricts to critical points for which $q$ is such that $\nu_q(t)$ is replica symmetric this would essentially correspond to the usual Onsager correction, but in principle (i.e., if the model has more than 1-RSB) there are additional solutions which approximately maximize \eqref{eq:fTAP} and obtain the correct free energy $F_{N,\beta}$ in the sense of \eqref{eq:Fapx2}.

Let $\nabla_{\perp}H_N(m)$ denote the projection of the  gradient in $\mathbb{R}^N$ to the orthogonal space to $m\neq 0$ and $\frac{d}{dR}H_N(m):={\textstyle \frac{d}{dt}}\big|_{t=0}H_N(m+tm/\|m\|)$ denote the radial derivative. Then, denoting $q=q(m)=\|m\|^2/N$, since $\nabla_{\perp} q(m)=0$,   
 \begin{equation}
 \label{eq:TAPEqGrad}
 \nabla\left( \frac{\beta}{N} H_N(m) +\frac{1}{2}\log(1-q)+F_{\beta}(q)\right)=0
 \end{equation}
 if and only if
 \begin{equation}
 \label{eq:TAPeq}
 \nabla_{\perp}H_N(m) = 0 \mbox{\ \ \ and \ \ }  
 {
 	\frac{1}{\sqrt N} \frac{d}{dR}}H_N(m)=\frac{\sqrt q}{\beta}\Big(\frac{1}{1-q} -2\frac{d}{dq}F_\beta(q)\Big).
 \end{equation}
The derivatives of the free energy of a mixed model in each of the coefficients $\gamma_p$ can be be computed from the Parisi formula, see \cite{Talag}. Using a similar argument, (see also the proof of Corollary \ref{cor:VTAPsol}) one can show that
\begin{equation}
\label{qe:ddqFq}
\frac{d}{dq}F_\beta(q)=\frac12 \beta^2\int \Big(
\frac{d}{dq}\nu_q(1)-\frac{d}{dq}\nu_q(x)
\Big)d\mu_\beta^q,
\end{equation}
where $\mu_\beta^q$ is the Parisi measure that corresponds to the mixture $\nu_q$ and
\[
\frac{d}{dq}\nu_q(x) = \big(
\nu'((1-q)x+q)-\nu'(q)
\big)(1-x)-(1-q)\nu''(q)x.
\]

For multi-samplable $q$ we have  the following equivalent formulation of the TAP equations \eqref{eq:TAPEqGrad} and \eqref{eq:TAPeq}. 
\begin{cor}\label{cor:ddqF}For multi-samplable $q>0$, \eqref{eq:TAPeq} is equivalent to
	\begin{equation*}
	\nabla_{\perp}H_N(m) = 0 \mbox{\ \ \ and \ \ }  
	{
		\frac{1}{\sqrt N} \frac{d}{dR}}H_N(m)=2\sqrt q\frac{d}{dq}\Es(q).
	\end{equation*}
\end{cor}

\begin{proof}
Denote the function on the right-hand side of \eqref{eq:TAP} by $\Phi(q)$. If $q\in(0,1)$ is multi-samplable, then by Theorem \ref{thm:FE}, 
$$\Phi(q)=\lim_{N\to\infty}\E F_{N,\beta}=\max_{t\in[0,1)}\Phi(t).$$  In particular, $\Phi'(q)=0$.
Namely,
\begin{equation*}
	\beta\frac{d}{dq}\Es(q)-\frac{1}{2}\frac{1}{1-q}+\frac{d}{dq} F_{\beta}(q)=0,
\end{equation*}
and the corollary follows.
\end{proof}

The TAP equations were first proved for the Sherrington-Kirkpatrick model at high temperature by Talagrand \cite{TalagrandBookI}, and later using Stein's method by Chatterjee \cite{ChatterjeeTAP}. 
Bolthausen \cite{BolthausenTAP} introduced a recursive scheme that converges to solutions of the TAP equations for the Sherrington-Kirkpatrick model below the AT line.  Using this recursive scheme, he showed in \cite{BolthausenMorita} that the free energy converges to the replica symmetric expression at high temperature. Recently, Auffinger and Jagannath \cite{AuffingerJagannthSpinDist,AuffingerJagannathTAP} proved that for generic mixed models with Ising spins the TAP equations are satisfied within each pure state in an appropriate sense. Finally, in \cite{TAPChenPanchenkoSubag} the TAP equations were established by Chen, Panchenko and the author for any mixed model with Ising spins and $q$ in the support of the Parisi measure, in the form of critical point equations for the Hamiltonian plus a TAP correction which is defined as a replicated free energy on a band similarly to the current work.

\subsection*{Remarks on methods of proof} All results stated so far will mostly follow from first principles and concentration properties. In particular, we emphasize that their proofs do not use the Parisi formula and the ultrametricity property. 
We remark that for the discussion about the TAP equations in Section \ref{subsec:TAPeq}
we assume the derivatives of the limiting free energy and ground state exist, which is only known using the Parisi formula. The proof of the convergence of the free energy \cite{FreeEnergyConvergence}, however, does not require the Parisi formula as we mentioned above.

The main concentration property we use is the uniform concentration of Proposition \ref{prop:concExpec}. It will follow from the well-known concentration of free energies as Lipschitz functions of the i.i.d. Gaussian disorder coefficients, see e.g. \cite[Theorem 1.2]{PanchenkoBook}. 
Another concentration result will be used to show that for a Hamiltonian that does not contain  $1$-spin interaction, $q=0$ is multi-samplable and therefore the limiting replicated free energy coincides with the usual free energy, see the discussion after Corollary \ref{cor:classicalTAP}. Namely, we will use
the super concentration property proved by Chatterjee \cite{ChattBook,ChatterjeeDisChaos} which bounds the variance of the free energy, which we will combine with the main result of  Paouris and Valettas
\cite{paouris2018} to obtain exponential decay for the left tail of the free energy. We remark  that if one does not insist  not to use the Parisi formula and ultrametricity, then the fact that $q=0$ is multi-samplable can be proved more easily from the description of the overlap array by Ruelle probability cascades (see e.g. \cite{PanchenkoBook}), since $0$ is in the support of the Parisi measure.

\subsection{\label{subsec:Parisi}Connections to the Parisi formula} 
The results we have stated so far apply to any overlap $q$ which is  multi-samplable. One would like to know how  can one `get' such overlaps $q$? To answer this, we start with the following result about the \emph{Parisi measure} $\mu_\beta([0,q]):=x_\beta(q)$.
\begin{cor}
	\label{cor:ParisiSupp}	If $q\in{\rm supp}(\mu_\beta)$, then $q$ is multi-samplable.
\end{cor}
The corollary will follow from the more general Theorem \ref{thm:qsamp}, the proof of which uses explicitly the Parisi formula and the optimality condition of Talagrand \cite{Talag} stated in Proposition \ref{prop:optimalbeta} below. A more insightful argument for the corollary  was kindly communicated to the author by D. Panchenko, for which the author is thankful, and goes as follows.

Recall that the Hamiltonian $H_N(\bs)$, or the corresponding mixture $\nu(t)$, are called generic if
\begin{equation}
\label{eq:generic}
\sum_{{\rm odd\ }p}p^{-1}\mathbf 1\{\gamma_p\neq 0\}=\sum_{{\rm even\ }p}p^{-1}\mathbf 1\{\gamma_p\neq 0\}=\infty.
\end{equation}
For generic models, the Parisi distribution coincides with the limiting overlap distribution (see e.g. \cite[Section 3.7]{PanchenkoBook}),
\begin{equation}
\label{eq:limodist}
\lim_{N\to\infty}\E G_{N,\beta}^{\otimes 2}\{R(\bs_1,\bs_2)\in \cdot\}=\mu_\beta(\cdot).
\end{equation}
In fact, for generic models, the limiting distribution of the infinite array of overlaps $\{R(\bs_i,\bs_j)\}_{i,j\geq1}$ under $G_{N,\beta}^{\otimes \infty}$ is equal to that of the array of overlaps from the Ruelle probability cascade that corresponds to $\mu_\beta$ in the case of finite RSB, or can be approximated by a sequence of cascades in case of full-RSB (see p. 102 and Section 3.7 of \cite{PanchenkoBook}). From its definition, one sees that for a Ruelle probability cascade that corresponds to some measure $\mu$, for any fixed $k$, the probability to sample $k$ points with overlap $q$ in the support of $\mu$ is positive. Therefore, for generic models, the corollary above follows. This can be extended to non-generic models by approximating any such model $\nu$ by a sequence of generic models $\nu_n$ and showing, using properties of the Parisi formula, that if $q$ belongs to the support of the Parisi measure for $\nu$ then there exists a sequence $q_n\to q$ of overlaps, each in the support of the measure of the corresponding $\nu_n$. Then, from the fact that $q_n$ are multi-samplable (for $\nu_n$), one can also verify that $q$ is multi-samplable (for $\nu$).

To describe the more general condition for $q$ being multi-samplable mentioned above, we recall the optimality
condition for the Parisi distribution $x_\beta$ proved by Talagrand \cite[Proposition 2.1]{Talag}.
Given a distribution function $x:[0,1]\to[0,1]$ such that $x(\hat q)=1$ for some $\hat q<1$, define the functions
\begin{align}
\label{eq:Phi}\Phi(t)&=\beta^2 \nu'(t)-\int_0^t\frac{ds}{(\int_s^1x(r)dr)^2},\\
\label{eq:phi}\phi(s)&=\int_0^s\Phi(t)dt.
\end{align}
\begin{prop}[Talagrand \cite{Talag}] 
	\label{prop:optimalbeta}The Parisi distribution $x_\beta$ is the unique distribution $x$ such that for $\mu([0,t]):=x(t)$,
	\begin{equation}
	\label{eq:S}
	{\rm supp}(\mu)\subset \mathcal S:=\Big\{s\in[0,1]: \phi(s)=\sup_{t\in[0,1]}\phi(t)  \Big\}.
	\end{equation}
\end{prop}
We denote by $\mathcal{S}_{\beta}$ the set defined by \eqref{eq:S} associated to the Parisi distribution $x_\beta$.
We will prove the following sufficient, and necessary in the generic case, condition for multi-samplable overlaps.
\begin{thm}
	\label{thm:qsamp}
	For any mixture $\nu$ and  $\beta>0$, 
	\begin{equation}
	\label{eq:q_direct}
	q\in \mathcal S_\beta \implies q \mbox{ is multi-samplable.}
	\end{equation}
	If we assume  that $\nu$ is generic, also
	\begin{equation}
	\label{eq:q_converse}
	q\in \mathcal S_\beta \impliedby q \mbox{ is multi-samplable.}
	\end{equation}
\end{thm}

Chen and Sen \cite{ChenSen} and  Jagannath and Tobasco \cite{JagannathTobascoLowTemp} proved a zero temperature analogue of the Parisi formula for the ground-state energy
\begin{equation}
\label{eq:EsParisi}
\Es
:=\Es(1) = \inf_{(\alpha,c)}\mathcal P_{\infty}(\alpha,c,\nu),
\end{equation}
where the infimum is over all non-decreasing, right continuous and integrable functions $\alpha:[0,1)\to[0,\infty)$  and $c>0$, and where
\begin{equation}\label{eq:Pinfty}
\begin{aligned}
&\mathcal P_{\infty}(\alpha,c,\nu)=\\&\quad\frac12\Big(
\nu'(1)\big(\int_0^1\alpha(s)ds+c\big)-\int_0^1\nu''(q)\big(
\int_0^q\alpha(s)ds
\big)dq
+\int_0^1\frac{dq}{\int_q^1\alpha(s)ds+c}
\Big).
\end{aligned}
\end{equation} 
From convexity, there exists a unique minimizer to the above, which we will denote by $(\alpha_\infty,c_\infty)$.

Let $x_\beta^q$ denote the Parisi distribution corresponding to the mixture $\nu_q$ (see \eqref{eq:nuq}) at inverse-temperature $\beta$ and let $(\alpha_\infty^q,c_\infty^q)$ denote the optimizer of the zero temperature Parisi formula corresponding to the mixture $\hat\nu_q(t):=\nu(qt)$. Note that $\Es(q)$ is the ground-state energy corresponding to the mixture $\hat\nu_q(t)$. We will prove the following, which in particular applies to overlaps in the support of the Parisi distribution and, in the generic case, to multi-samplable overlaps.

\begin{prop}
	\label{prop:Parisidistributions}
	For any mixture $\nu$ and $\beta>0$, if $q\in\mathcal S_\beta$ then
	\begin{equation}
	\label{eq:xbetaq} x^q_\beta(t) = x_\beta(q+(1-q)t),
	\end{equation}
	and if we further assume that $q>0$ then
	\begin{equation}
	\label{eq:alphaq} \big(\alpha^q_\infty(t),c^q_\infty\big)=\Big(\beta x_\beta(qt),\frac{\beta}{q}\int_q^1x_\beta(s)ds\Big).
	\end{equation}
\end{prop}

Recall that by Corollary \ref{cor:classicalTAP}, for the largest multi-samplable $q$, the correction in the TAP representation is the usual Onsager reaction term. We denote the right-most point in the support of the Parisi measure by
\begin{equation*}
	q_P=q_{P,\beta}:=\max {\rm supp}(\mu_\beta).
\end{equation*}
From \eqref{eq:xbetaq}, if $q=q_P$, then $x_\beta^q\equiv 1$ which corresponds to the replica symmetric solution. In particular, we will conclude the following.

\begin{cor}
	\label{cor:topSuppParisi}
	If $q=q_P$, then $F_\beta(q)$ is given by \eqref{eq:RSTAP}.
\end{cor}

The results from this section will be proved in Section \ref{sec:Parisipfs}.
\begin{rem}[Non-equivalence of ${\rm supp}(\mu_\beta)$ and $\mathcal S_\beta$]\label{rem:Sbeta}
	Talagrand's \cite{Talag} extremality condition, stated in Proposition \ref{prop:optimalbeta}, implies that ${\rm supp}(\mu_\beta)\subset \mathcal S_\beta$. It is reasonable to expect that typically ${\rm supp}(\mu_\beta) = \mathcal S_\beta$. In this remark we discuss counterexamples where the containment is strict, and thus there are overlaps which can not be detected directly at the level of $\mu_\beta$ for which the TAP representation \eqref{eq:TAP} and other results from the previous sections hold.

	In \cite[Proposition 2.3]{Talag}, Talagrand showed that, assuming $\gamma_1=0$, $\mu_\beta=\delta_0$ if and only if
	\begin{equation}
	\label{eq:crit}
	\max_{s\in[0,1]} f_\beta(s)=f_\beta(0)=0, \mbox{\ \ \ where \ \ }f_\beta(s):= \beta^2\nu(s)+\log(1-s)+s.
	\end{equation} 
	This follows, as he noted, since for the distribution $x\equiv1$, $\phi(s)=f_\beta(s)$, and for the corresponding measure $\mu=\delta_0$, ${\rm supp}(\mu)=\{0\}\subset \mathcal S$ if and only if \eqref{eq:crit} holds. 
	
	Let $\beta_c$ be the critical inverse-temperature, i.e., the largest value of $\beta$ such that \eqref{eq:crit} holds. 
	Since $f_{\beta_c}'(0)=0$, if $f_{\beta_c}''(0)<0$, which is the case e.g. when $\gamma_2=0$, then from continuity in $\beta$,  the maximum in \eqref{eq:crit} is still equal to $0$ for $\beta=\beta_c$ and there exists another point $s>0$ such that $f_{\beta_c}(s)=0$ and therefore $\mu_{\beta_c}=\delta_0$ and
	\begin{equation}
	\label{eq:crit2}
	{\rm supp}(\mu_{\beta_c})=\{0\}\subsetneq \mathcal S_{\beta_c}.
	\end{equation}
	When $f_{\beta_c}''(0)=0$, and thus $\beta_c:=\nu''(0)^{-1/2}$, \eqref{eq:crit2} does not necessarily occur. For example, if $\nu''(s)^{-1/2}$ is concave on $(0,1)$ it is not difficult to verify that $\mu_{\beta_c}=\delta_0$ and $\mathcal S_{\beta_c}=\{0\}$.

	Another interesting example is
	the $2$-spin model, $\nu(s)=\frac12 s^2$. In this case, see \cite{multipleoverlap},
		\begin{equation*}
			\mu_\beta=
			\begin{cases}
				\delta_0 &\mbox{if }\beta\leq1\\
				\delta_{1-1/\beta} &\mbox{if }\beta>1
			\end{cases} 
			\mbox{\ \ \ \ \ and thus\ \ \ \ \ }
			\mathcal S_\beta=
			\begin{cases}
				\{0\} &\mbox{if }\beta\leq1\\
				[0,1-1/\beta] &\mbox{if }\beta>1
			\end{cases}
			,
		\end{equation*} 
	and the containment is strict  ${\rm supp}(\mu_{\beta})\subsetneq \mathcal S_{\beta}$ whenever $\beta>1$, while the model is replica symmetric.
\end{rem}

\subsection{\label{subsec:PS}Energies on the ultrametric tree}

While it allowed him to (non-rigorously) compute the free energy,  Parisi's celebrated replica symmetry breaking solution \cite{ParisiFormula,Parisi} did not fully explain the physical properties of the Sherrington-Kirkpatrick model. After its discovery, a rich theory was developed in the physics literature to complete the picture and relate the solution to the structure of the Gibbs measure. In \cite{ParisiOrderPar} Parisi made a connection between the order parameter from his solution and the overlap distribution under the Gibbs measure. A key ingredient in Parisi's solution was the ansatz that only ultrametric replica matrices should be considered when applying the so-called replica trick. Mezard et al. \cite{MPSTV2,MPSTV1} discovered that the Gibbs measure splits into ``pure states'' that are organized in an ultrametric, or hierarchical, structure according to their overlaps (see also the book \cite{MPVspinglass} by Mezard, Parisi and Virasoro). 

Two stability properties of the Gibbs measure were established several years later in breakthroughs by Ghirlanda and Guerra \cite{GhirlandaGuerra} and Aizenman and Contucci \cite{AizenmanContucci}.
Arguin and Aizenman \cite{ArguinAizenman} proved that, assuming that the overlap distribution is supported on a finite set in the thermodynamic limit, the Aizenman-Contucci stochastic stability implies ultrametricity. Panchenko \cite{PanchenkoUltrametricity_finite_support_II,PanchenkoUltrametricity_finite_support_I} proved that, under the same assumption, the Ghirlanda-Guerra identities imply ultrametricity.
The fact that the Ghirlanda-Guerra identities together with ultrametricity determine the distribution of the overlap array as a Ruelle probability cascade \cite{Ruelle} or a limit of such was well-known, see the works of Baffioni and Rosati \cite{BaffioniRosati}, Bovier and Kurkova \cite{BovKurk1} and Talagrand \cite{TalagrandBookI}. 
Finally, in his seminal work \cite{ultramet}, Panchenko proved that the  ultrametricity follows from the Ghirlanda-Guerra identities, without any further assumption, implying in particular asymptotic ultrametricity for arbitrary generic models.

Connections to ultrametrcity also appeared, for instance, in the work of Bolthausen and Sznitman on the coalescent process \cite{BolthausenSznitman}, the work of Aizenman, Sims and Starr  \cite{AizenmanSimsStarrROST} on random overlap structures and related works by Arguin \cite{Arguin2008} and Arguin and Chatterjee \cite{ArguinChatterjee}.

Most relevant to our results in this section is the work of  Jagannath \cite{JagannathApxUlt}, where he  established an approximate analogue of the ultrametric structure at large, finite dimension $N$.
Roughly  speaking, he showed that the Gibbs measure splits into pure states that are clustered according to the typical overlaps of samples from them. This result generalized an earlier construction of Talagrand \cite{TalagrandPstates} for pure states, which assumed that $\mu_\beta(q_P)>0$.

To relate our analysis to the ultrametric tree, we will define an embedding of the tree in $\BN$ by associating a point to each pure state or cluster of states.  
The leaves will be the magnetizations of the pure states (i.e., center of mass w.r.t. to $G_{N,\beta}$), and to associate a point each of the clusters, we will take the unweighted (and this is important) average of the magnetizations of pure states in it. The precise definition will be  given in the proof of the corollary below, whose derivation from \cite{JagannathApxUlt} is rather straightforward, see Section \ref{sec:Ultra}.
By full regular tree we mean a rooted tree such that all leaves have the same depth and all vertices which are not leaves have the same degree. Also,
for any two vertices $u$ and $v$, $u\wedge v$ denotes the least common ancestor, and we write $u\leq v$ if $v$ is an ancestor of $u$ or $v=u$.
\begin{cor}\label{cor:embed}
	Suppose that $\nu(s)$ is generic and $\beta>0$ is such that $|{\rm supp}(\mu_\beta)|>1$. For any positive $c_N\to 0$, there exist positive $\eta_N, \epsilon_N\to0$, integer $d_N\to\infty$ and, depending on whether $|{\rm supp}(\mu_\beta)|$ is infinite or finite, integer $r_N\to\infty$  or constant $r_N=r$, respectively, such that the following holds.
	
	There exits a (random) full regular tree with vertex set $V\subset \BN$, of degree $d_N$ and depth $r_N$, 
	with leaves $m_1,\ldots,m_{K_N}\in \sqrt{q_P}\cdot \SN$, where $K_N=d_N^{r_N}$, and disjoint sets $B_i\subset  \SN$,
	such that  with probability tending to $1$, for any $i,j\leq K_N$ and $u,v\in V$:
	\begin{enumerate}
		\item\label{enu:t1} $\displaystyle \vphantom{\Big(} G_{N,\beta}\big(\cup_{i\leq K_N} B_i\big) > 1-\epsilon_N$\ \, and \ \, $\displaystyle \min_{i\leq K_N}G_{N,\beta}(B_i) > c_N$.
		\item\label{enu:t12} $\displaystyle \vphantom{\Big(} m_i \leq v\implies B_i\subset  {\rm Band}(v,\eta_N)$.
		\item\label{enu:t2} $\displaystyle \vphantom{\Big(} |R(u,v)-R(u\wedge v,u\wedge v)|<\epsilon_N$.
		\item\label{enu:t3} $\displaystyle \vphantom{\bigg(} G_{N,\beta}^{\otimes 2}\Big(\big|R(\bs_1,\bs_2)-R(m_i\wedge m_j,m_i\wedge m_j)\big|<\epsilon_N\ \Big|\ \bs_1\in B_i,\,\bs_2\in B_j\Big)>1-\epsilon_N$.
	\end{enumerate}
	Moreover, if $\mu_\beta(\qs)=w>0$, then the sequence $G_{N,\beta}(B_{i})$ weakly converges to  a Poisson-Dirichlet distribution of parameter $1-w$, (with the convention that for $i>K_N$, we set $G_{N,\beta}(B_i)=0$) and if $\mu_\beta(\qs)=0$ then $\max_{i}G_{N,\beta}\big(B_i\big)\to0$ in probability.\footnote{It is expected that typically $\mu_\beta(q_P)>0$.  It was shown that the support of $\mu_\beta$ is $\{0\}$ or an interval $[0,q_P]$ for any $\beta$ if and only if $\nu''(t)^{-1/2}$ is concave on $(0,1]$, see Proposition 1 in \cite{GSFollowing} which collects several results from \cite{PaCheChaos,JagannathTobascoBdsCplxSph,Talag}. In this case, there is also a formula for $\mu_\beta(q_P)$ from which one can check whether it positive or not. In particular, one can easily find examples where $\mu_\beta(q_P)>0$.
	Talagrand gives in p. 345 of \cite{Talag} an example with full-RSB such that $\mu_\beta(q_P)=0$ for some specific $\beta$. } 
\end{cor}
We note that if $\mu_\beta(\qs)>0$, Jagannath \cite{JagannathApxUlt} also proves that the weights $G_{N,\beta}(B_{i})$ converge to a Ruelle probability cascade, when ordered appropriately. If $\mu_\beta(\qs)=0$, he proves that if one clumps together pure states that are descendants of the same inner vertex at some level determined by a self-overlap lower than $q_P$, one also obtains a sequence of weights converges to Poisson-Dirichlet or Ruelle proability cascade when ordered appropriately. These results about the weights also apply in the setting of the corollary above. However, stating them precisely involves certain technicalities that we wish to avoid. We refer the interested reader to \cite{JagannathApxUlt} for further details.

We remark that if $c_N$ goes to $0$ slowly enough, from \eqref{eq:limodist} and the construction in the proof of Corollary \ref{cor:embed} based on \cite{JagannathApxUlt}, it is not difficult to see that the following hold. First,
$Q=\{R(v,v): v\in V\}$ is close to ${\rm supp}(\mu_\beta)$ in the sense that
for any $\epsilon>0$, w.h.p.  any $q\in Q$ is within distance $\epsilon$ from some $q'\in {\rm supp}(\mu_\beta)$ and vice versa (i.e., their Hausdorff distance is small). Second, for all $v$ with the same depth in the tree $V$, the self-overlap $R(v,v)$ is roughly the same. Lastly, if $0\in {\rm supp}(\mu_\beta)$ then we can take the root vertex of $V$ to be the origin in $\mathbb R^N$. 

The point in constructing the embedding above is the following connection with the replicated free energy and the Hamiltonian, uniformly on the tree. 
\begin{cor}[Maximality of energies on the tree]\label{cor:Eontree}
	Assume the setting of Corollary \ref{cor:embed} and let $V$ be vertex set of the tree obtained in the corollary with some $c_N>0$. Suppose that $c_N,\delta_N,\rho_N\to0$ and $n_N\to\infty$ sufficiently slowly. Then,
	\begin{equation}
	\label{eq:Eontree1}
	\sup_{v\in V}\bigg|\frac{\beta}{N}H_N(v)+\FF(v) - F_{N,\beta}\bigg|\longrightarrow 0,\quad\mbox{in probability},
\end{equation} 
and therefore
	\begin{equation}
	\label{eq:Eontree2}
		\sup_{v\in V}\bigg| \frac{H_N(v)}{N} - \max_{\|m\|=\|v\|}\frac{H_N(m)}{N}   \bigg| \longrightarrow 0,\quad\mbox{in probability}.
	\end{equation} 
\end{cor}

Recall \eqref{eq:limsupTAP} and note that from the same argument used there, in the setting of Corollary \ref{cor:Eontree}, with $q=\|v\|^2/N$, in probability,
\begin{equation}
	\label{eq:TAPcorrV}
	\lim_{N\to\infty}\sup_{v \in V} \left| \frac{\beta}{N} H_N(v) +\frac{1}{2}\log(1-q)+F_{\beta}(q) -F_{N,\beta}\right| =0.
\end{equation}
That is, the points $v\in V$ on the tree are approximate maximizers of \eqref{eq:fTAP}. From this we will conclude that all points of the tree are approximate TAP solutions in the following sense. One should compare the bound to the typical size of $\|\nabla H_N(m)/N\|$ of order $O(N^{-1/2})$.
	
	\begin{cor}[Vertices are approximate TAP solutions] \label{cor:VTAPsol}In the setting of Corollary \ref{cor:Eontree}, 
		 \begin{equation}
		 \label{eq:TAPEqGradV}
		 \sqrt N\cdot \sup_{v\in V}\left\|\nabla\left( \frac{\beta}{N} H_N(v) +\frac{1}{2}\log(1-q)+F_{\beta}(q)\right)\right\| \longrightarrow 0,\quad\mbox{in probability}.
		 \end{equation}
	\end{cor}

For the barycenters of pure states, which are the leaves of the tree from the corollaries above, \eqref{eq:TAPcorrV} or its analogue in the Ising case was proved in the works \cite{BeliusKistlerTAP,geometryMixed,ChenPanchenkoTAP,geometryGibbs}  we mentioned above in connection to Theorem \ref{thm:FE}. Since the self-overlap of the centers of pure states is equal to the maximal point in the support of the Parisi measure, in this case $F_{\beta}(q)$ is given by \eqref{eq:RSTAP} (or its Ising analogue). In \cite{AuffingerJagannathTAP} it was proved that for models with Ising spins the barycenters of pure states satisfy the TAP equations for a different formulation than ours.

We will prove the three corollaries above is Section	
\ref{sec:Ultra}.  We remark that in the full-RSB case, one can conclude from them that there exist paths from the origin to $\SN$ on which the energy is consistently maximal in the sense of \eqref{eq:Eontree2}. This was the inspiration for the optimization algorithm designed in \cite{GSFollowing}, which outputs an approximate maximizer of the energy $H_N(\bs)$ on $\SN$ in polynomial time in $N$.

\noindent \textbf{Acknowledgments.}
I am grateful to G\'{e}rard Ben Arous and Ofer Zeitouni for helpful discussions and comments, to Wei-Kuo Chen for a useful conversation on the optimality criterion in the Parisi formula, and to Dmitry Panchenko for a careful reading and valuable comments which improved the paper and played a role in  identifying multi-samplability as a key property for the TAP approach, and  for communicating the argument that follows Corollary \ref{cor:ParisiSupp} which was used in an earlier version to remove the assumption of genericity from several main results.
I am also thankful to the anonymous referees whose comments helped improve and clarify this manuscript. 
This work was supported by the Simons Foundation, the Israel Science Foundation (Grant
Agreement No. 2055/21) and a research grant from the Center for Scientific Excellence at the Weizmann Institute of Science.

\section{\label{sec:TAPcorrection}Computation of the TAP correction: proof of Proposition \ref{prop:FFlim}}

Our first aim will be to prove \eqref{eq:EFt} for fixed $q$, without the supremum over $[0,t]$. This will be accomplished in \eqref{eq:corro2_2} and	\eqref{eq:corro2_3}. Afterwards we will see why this is enough for the statement as in Proposition \ref{prop:FFlim} with the supremum. 

There are two key steps to the proof. The first is Lemma \ref{lem:rem1sp} which will allow us to move from a free energy defined using the mixture $\tilde\nu_q$ to one defined with $\nu_q$, which does not include $1$-spin interaction. Those two free energies will be defined with asymptotically infinitely many, orthogonal replicas (see 	\eqref{eq:FtHq} and \eqref{eq:FHq}). The second key step is
Lemma \ref{lem:0ismulti} which will show that in the absence of $1$-spin interaction, the free energy with replicas is asymptotically equivalent to the standard free energy with one replica.  The proof of this lemma will use a bound on the variance of the free energy proved by Chatterjee \cite{ChatterjeeDisChaos,ChattBook} and a result of Paouris and Valettas \cite{paouris2018}.

If $\|m\|^2\leq N$, the zero-width band
\begin{equation*}
	{\rm Band}\left(m,0\right)=\big\{ \bs\in\SN: R(\bs-m,m)=0 \big\}
\end{equation*}
is a sphere.
In some parts of the proof it will sometimes be more convenient to work with the following free energy defined over it. Let
\begin{equation}
	\label{eq:Ftp}	
	\begin{aligned}
		&\FF^o(m,n,\rho)\\
		&\,\,:=\frac{1}{ Nn}\log\int_{B_{N}(m,n,\rho,0)}\exp\Big\{\beta\sum_{i=1}^{n}(H_{N}(\bs_{i})-H_{N}(m))\Big\}d\bs_{1}\cdots d\bs_{n},
	\end{aligned}
\end{equation} 
where $d\bs_i$ denotes integration w.r.t. the uniform measure on $	{\rm Band}\left(m,0\right)$.

To relate the free energy above to $\FF(m,n,\rho,\delta)$, we will need the following.

\begin{lem}
	\label{lem:Lip} Denote by $\mathcal E_N(L) $ the event that $H_N(\bs)$ has Lipschitz constant bounded by $\sqrt N L$ on $\BN$. Then for any $c'>0$, for large enough $L>0$,
	$\P\{ \mathcal E_N(L) \} > 1-e^{-Nc'}$.
\end{lem}

\begin{proof}
	Note that the Lipschitz constant of $H_N(\bs)$ is equal to 
	\[
	\sup_{m\in\BN,\, \mathbf{u}:\,\|\mathbf{u}\|=1}\nabla H_N(m)\cdot \mathbf{u}.
	\]
	Since this is the supremum of a Gaussian field with variance bounded by some constant $c_\nu$ depending only on $\nu$, by the Borell-TIS inequality \cite{Borell,TIS} the lemma will follow if we can show that the Lipschitz constant has expectation bounded by $\sqrt{N}c_\nu'$. This follows, e.g., from  \cite[Lemma 6.1]{TAPChenPanchenkoSubag}.
\end{proof}

\begin{lem}
	\label{lem:Fo}
	On the event that $H_N(\bs)$ has Lipschitz constant bounded by $\sqrt{N} L$ on $\BN$,
	\begin{equation}
		\label{eq:FtFo}
		\begin{aligned}
			&\FF(m,n,\rho,\delta)\leq 
			\FF^o(m,n,\rho+\sqrt{8\delta})
			+\frac1N\log{\rm Vol}({\rm Band}\left(m,\delta \right))+\beta L \sqrt{2\delta},\\
			&\FF(m,n,\rho+\sqrt{8\delta},\delta)\geq 
			\FF^o(m,n,\rho)
			+\frac1N\log{\rm Vol}({\rm Band}\left(m,\delta \right))-\beta L \sqrt{2\delta},
		\end{aligned}
	\end{equation}
	where ${\rm Vol}$ denotes the uniform measure on $\SN$.
\end{lem}

\begin{proof}

For $m=0$ the lemma is trivial.	Assume that $q\in(0,1)$ and let $m\in\sqrt q\cdot \SN$. For $\bs\in\SN\setminus\{\pm m/\sqrt{q}\}$ define the projections
\begin{equation*}
	\hat \bs = \bs - \frac1q R(\bs,m)m \mbox{\ \ \ and \ \ }\bs' = m+\sqrt{N(1-q)}\frac{\hat \bs}{\|\hat \bs\|}.
\end{equation*}

We wish to find a bound, as a function of $\delta$, for the distance between $\bs$ and $\bs'\in 	{\rm Band}\left(m,0\right)$, uniformly in $q$ and $\bs\in {\rm Band}(m,\delta)\setminus\{\pm m/\sqrt q\}$. Given some $\delta\in(0,1]$, it is easy to see that the supremum of the distance is obtained by choosing $\sqrt q+\delta =1 $ and letting $\bs$ go to $m/\sqrt q\in \SN$, in which case the distance goes to $\sqrt{N(1-q+\delta^2)}=\sqrt{2N\delta}$. 
From this we have that for any $\delta>0$ and  
$\bs\neq \pm m/\sqrt{q}$,
\begin{equation*}
	\bs \in {\rm Band}(m,\delta )\implies \|\bs'-\bs\|\leq \sqrt{2N\delta}.
\end{equation*}
Therefore,
\begin{equation*}
	\bs_1,\bs_2 \in {\rm Band}(m,\delta)\implies \big|R(\bs_1',\bs_2')-R(\bs_1,\bs_2)\big|\leq \sqrt{8\delta}.
\end{equation*}

From the above, for $\bs_i\in {\rm Band}(m,\delta )$ different from 
$\pm m/\sqrt{q}$, 
\begin{equation*}
	(\bs_{1},\ldots,\bs_{n})\in B_{N}(m,n,\rho,\delta)\implies (\bs_{1}',\ldots,\bs_{n}')\in B_{N}^o(m,n,\rho+\sqrt{8\delta}),
\end{equation*}
and 
\begin{equation*}
	(\bs_{1}',\ldots,\bs_{n}')\in B_{N}^o(m,n,\rho)\implies
	(\bs_{1},\ldots,\bs_{n})\in B_{N}(m,n,\rho+\sqrt{8\delta},\delta)	 .
\end{equation*}
Assuming the Lipschitz bound, \eqref{eq:FtFo} follow by integrating each of $\bs_i$ in the direction  $m/\|m\|$ on ${\rm Band}(m,\delta)$.
\end{proof}

Given a real function $f(\bs)$ on $\SN$, define 
\begin{equation}
	\label{eq:Ftf}
	\FF(f,n,\rho):=\frac{1}{ Nn}\log\int_{B_{N}(0,n,\rho)}\exp\Big\{\beta \sum_{i=1}^{n}f(\bs_{i})\Big\}d\bs_{1}\cdots d\bs_{n},
\end{equation}
where we define
\[
B_{N}(0,n,\rho):  =\big\{(\bs_{1},\ldots,\bs_{n})\in(\SN )^{n}:\,|R(\bs_{i},\bs_{j})|<\rho,\,\forall i\neq j\big\}.
\]

Let $m\in\SN$ with $\|m\|^2/N=q \in[0,1)$. For $\bs_1,\bs_2\in 	{\rm Band}\left(m,0\right)$, \eqref{eq:Htilde} holds without the $o(1)$ term. Hence, if we denote by $\tilde H_N^q(\bs)$ 
the Hamiltonian that corresponds to 
$\tilde\nu_q(t)$ defined in \eqref{eq:nuq_tilde}, then
(by a change of variables)
\begin{equation*}
	\FF^o(m,n,\rho) \overset{d}{=} \frac{N-1}{N} \widetilde F_{N-1,\beta}\Big({\sqrt{\tfrac{N}{N-1}} \tilde H_{N-1}^q(\bs)},n,\rho/(1-q)\Big).
\end{equation*}
Let $n_N\to\infty$ and $\rho_N'\to0$ be arbitrary.  Since the ground-state energy of $\tilde H_{N-1}^q(\bs)$ concentrates, for large $N$,
\begin{equation}
	\label{eq:FtHq}
	\E \FF^o(m,n_N,\rho_N') = \E \widetilde F_{N,\beta}\big( \tilde H_{N}^q(\bs),n_N,\rho_N'/(1-q)\big)+o(1).
\end{equation}
For $m=0$ this holds by definition without the $o(1)$.

\begin{lem}
	\label{lem:rem1sp}
	Let $\tilde \xi(t)=\sum_{p=1}^{\infty}\gamma_p^2t^p$  be a general mixture and let $\xi(t)=\sum_{p=2}^{\infty}\gamma_p^2t^p$, in which the $1$-spin interaction is removed. 
	Let $\tilde H_N(\bs)$ and $H_N(\bs)$ be the corresponding  Hamiltonians. Then, 
	\begin{equation}
		\label{eq:rem1sp}
		\big|\E \FF \big( \tilde H_{N}(\bs),n,\rho\big) - \E \FF \big(  H_{N}(\bs),n,\rho\big) \big|\leq \beta C\gamma_1 \sqrt{1/n+\rho },
	\end{equation}	
	for some universal constant $C>0$.
\end{lem}

\begin{proof}
	As follows from a covariance calculation, we can define $\tilde H_N(\bs)$ and $H_N(\bs)$ on the same probability space such that 
	\begin{equation}
		\label{eq:1spdecomp}
		\tilde H_N(\bs) = H_N(\bs) + \gamma_1 \bs \cdot \mathbf X,
	\end{equation}
	with $\mathbf X$ independent of $H_N(\bs)$,
	where $\mathbf X=(X_i)_{i=1}^N$ is a vector of i.i.d. Gaussian variables $X_i\sim N(0,1)$.
	By Cauchy-Schwarz,
	\begin{equation*}
		\begin{aligned}
			&\big| \FF\big( \tilde H_{N}(\bs),n,\rho\big) -  \FF\big(  H_{N}(\bs),n,\rho\big) \big|
			\leq \frac{\beta}{Nn}\sup \Big|
			\sum_{i=1}^{n} \tilde H_{N}(\bs_{i})
			-\sum_{i=1}^{n}  H_{N}(\bs_{i})	
			\Big|\\
			&=\frac{\beta\gamma_1}{Nn}\sup \Big| \sum_{i=1}^{n} \bs_i \cdot \mathbf X\Big|
			\leq \frac{\beta\gamma_1}{Nn} \sqrt{N(n+n(n-1)\rho )}\|\mathbf X\|,
		\end{aligned}
	\end{equation*}
	where the supremum is over $(\bs_1,\ldots,\bs_n)\in B_N(0,n,\rho)$.
	And \eqref{eq:rem1sp} follows form 
	Jensen's inequality.
\end{proof}

Let $q\in[0,1)$ and let $m=m_N$ be a sequence such that $\|m_N\|=Nq$.
Denote by $H_N^q(\bs)$ 
the Hamiltonian that corresponds to 
$\nu_q(t)$ from \eqref{eq:nuq}. Then, from the lemma and \eqref{eq:FtHq}, 
\begin{equation}
	\label{eq:FHq}
	\E \FF^o(m,n_N,\rho_N') = \E \FF\big( H_{N}^q(\bs),n_N,\rho_N'/(1-q)\big)+o(1).
\end{equation}

By definition, the right-hand side of \eqref{eq:FHq} is bounded from above by the same with $n_N$ replaced by $1$, which is just the usual free energy of $H_{N}^q(\bs)$ (without multiple replicas). This free energy is exactly how we defined $F_\beta(q)$, and therefore
\begin{equation}
	\label{eq:corro2}
	\limsup_{N\to\infty} \E\FF^o(m,n_N,\rho_N')\leq F_\beta(q).
\end{equation}

Let $\rho_N,\delta_N\to0$ and define $\rho_N'=\rho_N+\sqrt{8\delta_N}$. Assume that $\delta_N$ is sufficiently large (i.e., not exponentially small in $N$) so that independently of $q$,
\begin{equation}\label{eq:vol}
\lim_{N \to\infty}\frac1N\log{\rm Vol}({\rm Band}\left(m,\delta_N \right)) = \frac12\log(1-q).
\end{equation} 
From Lemmas \ref{lem:Lip} and \ref{lem:Fo} and the fact that both free energies in the latter lemma concentrate (see e.g. \cite[Theorem 1.2]{PanchenkoBook}), we also have that
\begin{equation}
	\label{eq:corro2_2}
	\limsup_{N\to\infty} \E\FF(m,n_N,\rho_N,\delta_N)\leq\frac12\log(1-q)+F_\beta(q).
\end{equation}

To prove a matching lower bound to \eqref{eq:corro2_2} for the liminf, we will need the following lemma.
\begin{lem}
	\label{lem:0ismulti}
	Suppose that the mixture $\nu(x)=\sum_{p\geq2}\gamma_{p}^{2}x^{p}$ does
	not contain a $1$-spin interaction (i.e., $\gamma_1=0$). If $n_N\to\infty$ and $\rho_N\to0$ sufficiently slowly, then
	\begin{equation}
		\label{eq:0ismulti}
		\lim_{N\to\infty}	\E \FF^o(0,n_N,\rho_N) =  \lim_{N\to\infty}\E F_{N,\beta},
	\end{equation}
	and the overlap $q=0$ is multi-samplable.
\end{lem}

We postpone its proof to Section \ref{subsec:pfFF4} below. The proof will use a bound on the variance of the free energy proved by Chatterjee \cite{ChatterjeeDisChaos,ChattBook} and a result of Paouris and Valettas \cite{paouris2018} that, combined with the latter, implies that the left tail of the free energy decays exponentially.

By applying the lemma to the Hamiltonian $H_N^q(\bs)$ with the mixture $\nu_q(t)$, assuming that 
$n_N\to\infty$ and $\rho'_N\to0$ sufficiently slowly, the limit of the free energy in the right-hand side of \eqref{eq:FHq} is equal to $F_\beta(q)$ and thus
\begin{equation}
	\label{eq:corro_3}
	\liminf_{N\to\infty} \E\FF^o(m,n_N,\rho_N')\geq F_\beta(q).
\end{equation}

Given some $\rho_N,\delta_N\to0$, for some smaller $\rho_N',\delta_N'$ we have that
\[
\FF(m,n,\rho_N,\delta_N)\geq \FF(m,n,\rho_N'+\sqrt{8\delta_n'},\delta_N').
\]
By assuming that  $\rho_N,\delta_N\to0$ sufficiently slow, we may assume the same about $\rho_N',\delta_N'$.
Hence, using Lemmas \ref{lem:Lip} and \ref{lem:Fo}, concentration and \eqref{eq:vol},
\begin{equation}
	\label{eq:corro2_3}
	\liminf_{N\to\infty} \E\FF(m,n_N,\rho_N,\delta_N)\geq\frac12\log(1-q)+F_\beta(q).
\end{equation}

\begin{lem}
	\label{lem:Fbcont}
	$F_\beta(q)$ is a continuous function of $q\in[0,1)$.
\end{lem}
\begin{proof}
	The Hamiltonians corresponding to $\nu_q(t)$ for two different values $q$ and $q'$ can be coupled by using the same disorder coefficients $J^{(p)}_{i_1,\ldots,i_p}$ in \eqref{eq:Hamiltonian} to define both. The difference of the two Hamiltonians is another mixed Hamiltonian, say with mixture $\sum_{p\geq2}a_{p,q,q'}^2t^p$. One can check that $\sum_{p\geq2}a_{p,q,q'}^2p^3\to0$ as we let $q-q'\to0$, uniformly in $q,q\in[0,1]$.
	Hence, by \cite[Lemma 6.1]{TAPChenPanchenkoSubag}, the expectation of the ground state of this difference Hamiltonian normalized by $N$ goes to $0$ as $q-q'\to0$ (where we use the fact that at the origin the Hamiltonian is equal to $0$).
	Using the Borell-TIS inequality \cite{Borell,TIS} and 
	the concentration of the free energies \cite[Theorem 1.2]{PanchenkoBook}, the lemma easily follows.
\end{proof}

We are now ready to prove Proposition \ref{prop:FFlim}. Let $\bs=\bs_N\in\SN$ be some arbitrary sequence. First note that for any $q, q'>0$ and $\Delta>0$, 
\begin{equation}
	\label{eq:FF'}
	\begin{aligned}
		&\mbox{if}\quad  |\sqrt q - \sqrt{q'}| \leq \Delta,\ \ \delta-\delta'\geq\Delta,\ \ \rho-\rho'\geq 2\Delta\\
		&\mbox{then}\quad \ \
		\FF(\sqrt q \bs,n,\rho,\delta) \geq \FF(\sqrt{q'}\bs,n,\rho',\delta'),
	\end{aligned}
\end{equation}
since $B_N(\sqrt{q}\bs,n,\rho,\delta)\supset B_N(\sqrt{q'}\bs,n,\rho',\delta')$.

Note that by a diagonalization argument, from \eqref{eq:corro2_3}, if  $\delta_N',\rho_N'\to0$ and $n_N\to\infty$ tend to their limits sufficiently slowly, then for arbitrary $\bs \in  \SN$,
\begin{equation}
	\label{eq:FtLB}
	\liminf_{N\to\infty} \inf_{q\in Q'_N}\E\FF(\sqrt q \bs ,n_N,\rho_N',\delta_N')- \frac12\log(1-q)-F_\beta(q)\geq 0,
\end{equation}
for some subset $Q'_N\subset [0,1)$ which increases with $N$ and such that the closure of $\cup_N Q'_N$ is $[0,1]$.
Define $\Delta_N=\max\{|\sqrt x-\sqrt q|:\, x\in[0,1], q\in Q_N'  \}$, $\delta_N= \delta_N'+\Delta_N$ and $\rho_N= \rho_N'+2\Delta_N$.

From \eqref{eq:FF'}, for any $t\in(0,1)$, using the fact $\log(1-q)/2+F_\beta(q)$ is uniformly continuous on $[0,t]$,
\begin{equation}
	\label{eq:liminfFt}
	\liminf_{N\to\infty} \inf_{q\in[0,t]}\Big(\E\FF(\sqrt q \bs ,n_N,\rho_N,\delta_N)-\frac12\log(1-q)-F_\beta(q)\Big)
\end{equation}
is bounded from below by \eqref{eq:FtLB}. 

Let $c>0$ be some arbitrary small number. Fix some $q_0\in[0,1)$. Note that from \eqref{eq:corro2_2},  with the same $n_N$ as above, if $\rho$ and $\delta$ are small enough, then
\begin{equation*}
	\limsup_{N\to\infty} \Big(\E\FF(\sqrt{q_0} \bs ,n_N,\rho,\delta)-\frac12\log(1-q_0)-F_\beta(q_0)\Big)<c.
\end{equation*}
From \eqref{eq:FF'}, for any $\Delta<\min\{\rho/4,\delta/2\}$,
\begin{equation*}
	\limsup_{N\to\infty}\sup_{q:|\sqrt q-\sqrt{q_0}|<\Delta} \E\FF(\sqrt{q} \bs ,n_N,\rho_N,\delta_N)
	\leq
	\limsup_{N\to\infty} \E\FF(\sqrt{q_0} \bs ,n_N,\rho,\delta).
\end{equation*}
Since $[0,t]$ is compact and $\log(1-q)/2+F_\beta(q)$ is uniformly continuous on it, 
\begin{equation*}
	\limsup_{N\to\infty} \sup_{q\in[0,t]}\Big(\E\FF(\sqrt q \bs ,n_N,\rho_N,\delta_N)-\frac12\log(1-q)-F_\beta(q)\Big)<2c.
\end{equation*}

Combining the above completes the proof of Proposition \ref{prop:FFlim}. It remains to prove Lemma \ref{lem:0ismulti}.\qed

\subsection{\label{subsec:pfFF4}Proof of Lemma \ref{lem:0ismulti}}
First, we note that the fact that  \eqref{eq:0ismulti} implies that $q=0$ is multi-samplable follows from the concentration of $F_{N,\beta}$  and of $\FF^o(0,n_N,\rho_N)=\FF(0,n_N,\rho_N,\delta)$ (which follow, e.g., by \cite[Theorem 1.2]{PanchenkoBook}; see also \eqref{eq:conc_1pt} below), since
\begin{equation}
\label{eq:FG}
\FF^o(0,n_N,\rho_N)-F_{N,\beta} = \frac1{ Nn_N}\log G_{N,\beta}^{\otimes n_N}\Big\{  \forall i\neq j :\,  \big|R( \bs_i,\bs_j)   \big| < \rho_N\, \Big\}.
\end{equation}

Since for any $n$ and $\rho$,
\[
\FF^o(0,n,\rho)\leq \FF^o(0,1,1)=F_{N,\beta},
\]
to prove \eqref{eq:0ismulti} it will be enough to show that for arbitrary $n$ and $\rho$,
\begin{equation*}
	\liminf_{N\to\infty}\E \FF^o(0,n,\rho)\geq  \lim_{N\to\infty}\E F_{N,\beta}.
\end{equation*}

Again, from concentration, it will be enough to show that for any $n\geq1$ and	 $\rho,t>0$,
\begin{equation}
\label{eq:Gbd1}
\lim_{N\to\infty}\P \Big(
\FF^o(0,n,\rho)  > \E F_{N,\beta}-t \Big)=1.
\end{equation}
Note that
$$\FF^o(0,n,\rho)\geq \frac{1}{Nn}\log\big(M_{N,\beta}(1,\rho)\cdots M_{N,\beta}(n,\rho) \big),$$ 
where $M_{N,\beta}(1,\rho):=Z_{N,\beta}$ and for $n\geq2$,
\begin{align*}
M_{N,\beta}(n,\rho)&:=	\inf_{\bs_1,\ldots,\bs_{n-1}\in\SN } 
Z_{N,\beta,\rho}(\bs_1,\ldots,\bs_{n-1}),
 \\
Z_{N,\beta,\rho}(\bs_1,\ldots,\bs_n)&:=
\int_{\big\{\bs\in\SN:\,\max_{i\leq n}|R( \bs,\bs_i) | < \rho  \big\}}e^{\beta H_{N}\left(\boldsymbol{\sigma}\right)}d\bs.
\end{align*} 

Hence, to prove \eqref{eq:Gbd1} it will be enough to show that
for any $n\geq1$ and	 $\rho,t>0$,
\begin{equation}
\label{eq:Gbd2}
\lim_{N\to\infty}\P \Big(	\frac{1}{N}\log\big(M_{N,\beta}(n,\rho) \big)   \leq  \E F_{N,\beta} -t \Big)=0.
\end{equation}

Let $n\geq1$ and $\epsilon>0$ be arbitrary. Denote
\[
D_N(n) = \big\{(\bs_1,\ldots,\bs_{n-1})\in (\SN)^{n-1}:\, R(\bs_i,\bs_j)=0,\, \forall i\neq j \, \big\},
\]
and let
\begin{equation}
\label{eq:packing}
\big\{ (\bs_1^{\ell},\ldots, \bs_{n-1}^{\ell})   \big\}_{\ell=1}^K\subset D_N(n)
\end{equation}
be a maximal subset such that
\begin{equation}
\label{eq:Aell}
A_\ell:=\Big\{  
(\bs_1,\ldots,\bs_{n-1})\in (\SN)^{n-1}:\, \max_{i\leq n-1} \|\bs_i-\bs_i^\ell\|_2\leq \sqrt N \epsilon/2\,
\Big\}
\end{equation}
are disjoint sets for any two  $\ell<\ell'\leq K$  (maximal in the sense that for any subset of $D_N^K(n)$ that strictly contains \eqref{eq:packing} there are two subsets $A_\ell$, $A_{\ell'}$ which intersect). Then, from maximality, for any $(\bs_1,\ldots,\bs_{n-1})\in D_N(n)$ there is some $\ell\leq K$ such that
\begin{equation}
\label{eq:ell1}
\max_{i\leq n-1} \|\bs_i-\bs_i^\ell\|_2\leq \sqrt N \epsilon.
\end{equation}
Note that since $A_\ell$ are disjoint, denoting by $\mbox{Vol}$ the uniform (product) measure on $(\SN)^{n-1}$,
\[
K\leq \mbox{Vol}\big((\SN)^{n-1} \big)/\mbox{Vol}\big(A_\ell \big)<e^{Na},
\]
for sufficiently large $a=a(\epsilon,n)>0$ independent of $N$.

Now, let $\hat \bs_1,\ldots,\hat \bs_{n-1}\in\SN$ be arbitrary points. Suppose that $(\bs_1,\ldots,\bs_{n-1})\in D_N(n)$ are some points whose span contains that of $\hat\bs_1,\ldots,\hat\bs_{n-1}$, and assume  $\ell\leq K$ satisfies \eqref{eq:ell1}. Then, for any $\bs\in\SN$ such that $|R(\bs,\bs_i^\ell)|<\eta$ for all $i$,
\[
|R(\bs,\bs_i)|=|R(\bs,\bs_i-\bs_i^\ell)+R(\bs,\bs_i^\ell)|\leq \epsilon+\eta.
\]
Hence, given $\rho$, if we choose some $\epsilon$ and $\eta$ so that $\rho\geq\epsilon+\eta$, then
\begin{equation*}
\big\{ \bs:\,   \big|R( \bs,\bs_i)   \big| < \rho,\, \forall i\leq n-1 \big\}
\supset \big\{ \bs:\,   \big|R( \bs,\bs_i^\ell)   \big| < \eta,\, \forall i\leq n-1 \big\}.
\end{equation*}

From a union bound, 
\begin{align*}
	&\P \Big(	\frac{1}{N}\log\big(M_{N,\beta}(n,\rho) \big)  \leq  \E F_{N,\beta}-t \Big)\\
	&\leq \P \Big(	\inf_{\ell\leq K} \frac{1}{N}\log\big(Z_{N,\beta,\eta}(\bs_1^{\ell},\ldots, \bs_{n-1}^{\ell})\big)       \leq\E F_{N,\beta} -t \Big)
	\\
	 &\leq e^{Na}\cdot
	\P \Big(	\frac{1}{N}\log\big(Z_{N,\beta,\eta}(\bs_1,\ldots,\bs_{n-1})\big)  \leq\E F_{N,\beta} -t \Big)
	,
\end{align*}
for arbitrary $(\bs_1,\ldots,\bs_{n-1})\in D_N(n)$.

Hence, to complete the proof it will be enough to show that for any $n\geq1$, $\eta,t,c>0$ and  $(\bs_1,\ldots,\bs_{n-1})\in D_N(n)$, for large $N$,
\begin{equation}
\label{eq:GBn1}
	\P \Big(	\frac{1}{N}\log\big(Z_{N,\beta,\eta}(\bs_1,\ldots,\bs_{n-1})\big)  \leq  \E F_{N,\beta} -t \Big)
	<e^{-Nc}.
\end{equation}

From the Lipschitz property of Lemma \ref{lem:Lip}, it is enough to prove that for arbitrary $n\geq1$, $t,c>0$ (which may change from \eqref{eq:GBn1}) and large $N$,
\begin{equation*}
\P \Big(	\frac{1}{N}\log\big(Z_{N,\beta}(\bs_1,\ldots,\bs_{n-1})\big)\leq \E F_{N,\beta}   -t \Big)
<e^{-Nc},
\end{equation*}
where
\begin{equation}
Z_{N,\beta}(\bs_1,\ldots,\bs_{n-1}):=
\int_{\big\{\bs\in\SN:\,R( \bs,\bs_i)=0,\,\forall i\leq n-1  \big\}}e^{\beta H_{N}\left(\boldsymbol{\sigma}\right)}d\bs,
\label{eq:Zorth}
\end{equation}
with $d\bs$ denoting integration w.r.t. the uniform measure over the sphere of co-dimension $n-1$.
Up to scaling of the Hamiltonian by $\sqrt{N/(N-n+1)}$, \eqref{eq:Zorth}  has the same distribution as $Z_{N-n+1,\beta}$.
Hence, since $\E F_{N,\beta}$ converges, the proof of Lemma \ref{lem:0ismulti} is completed by the following lemma. \qed

\begin{lem}
	\label{lem:Flb}
	For any mixture $\nu(s)=\sum_{p\geq2}\gamma_{p}^{2}s^{p}$ that does
	not contain a $1$-spin interaction (i.e., $\gamma_1=0$), for any $c,\,t>0$, for
	large enough $N$,
	\begin{equation}
	\label{eq:smalldev}
	\mathbb{P}(F_{N,\beta}<\mathbb{E}F_{N,\beta}-t)<e^{-cN}.
	\end{equation}
\end{lem}

\begin{proof}
If $\gamma_p>0$ for infinitely many values of $p$, with $H_N(\bs)$ defined by \eqref{eq:Hamiltonian}, let $H^{p_0}_N(\bs)$ be defined through the same formula with summation over $p\leq p_0$. By the chaining argument, see \cite[Theorem 1.3.3]{RFG}, one can verify that $\frac1N(H_N(\bs)-H^{p_0}_N(\bs))$ is as small as we wish in expectation and by the Borell-TIS inequality \cite{Borell,TIS} also w.h.p., if we take $k$ to be large enough. From this and the concentration of the free energy it follows that it is enough to prove the lemma assuming, as we will, that $\gamma_p>0$ for $p=p_0$ for some $p_0$.

From H\"{o}lder's inequality, $F_{N,\beta}$
is a convex function of the disorder coefficients $J_{i_{1},...,i_{p}}^{(p)}$.
Therefore, by the main result of Paouris and Valettas \cite{paouris2018}  to prove \eqref{eq:smalldev} it will be 
enough to prove that asymptotically $\text{Var}(F_{N,\beta})=o(N)$.
To prove the bound on the variance, we adapt Chatterjee's
proof of a similar bound for the Sherrington-Kirkpatrick
model in Theorem 6.3 of his book \cite{ChattBook}. The adaptation is straightforward and is given here
		for the sake of completeness.
		
		First we have that
		\begin{equation}
		\frac{\partial}{\partial J_{i_{1},...,i_{p}}^{(p)}}NF_{N,\beta}=\beta\gamma_{p}N^{-\frac{p-1}{2}}\langle\sigma_{i_{1}}\cdots\sigma_{i_{p}}\rangle,\label{eq:Fgrad}
		\end{equation}
		where $\langle\cdot\rangle$ denotes averaging under $G_{N,\beta}$.
		Similarly, letting $\bs^{1},\bs^{2},\ldots$ be an i.i.d. sequence from $G_{N,\beta}$
		and assuming $\ell_{1},\ldots,\ell_{k}$ are all different, 
		\begin{align*}
			& \frac{\partial}{\partial J_{i_{1},...,i_{p}}^{(p)}}\langle\sigma_{i_{1}^{1}}^{\ell_{1}}\cdots\sigma_{i_{s_{1}}^{1}}^{\ell_{1}}\cdots\sigma_{i_{1}^{k}}^{\ell_{k}}\cdots\sigma_{i_{s_{k}}^{k}}^{\ell_{k}}\rangle\\
			& \begin{aligned}=\beta\gamma_{p}N^{-\frac{p-1}{2}} & \Big[\sum_{\ell_{j}=1}^{k}\langle\sigma_{i_{1}^{1}}^{\ell_{1}}\cdots\sigma_{i_{s_{1}}^{1}}^{\ell_{1}}\cdots\sigma_{i_{1}^{k}}^{\ell_{k}}\cdots\sigma_{i_{s_{k}}^{k}}^{\ell_{k}}\sigma_{i_{1}}^{\ell_{j}}\cdots\sigma_{i_{p}}^{\ell_{j}}\rangle\\
				& -k\langle\sigma_{i_{1}^{1}}^{\ell_{1}}\cdots\sigma_{i_{s_{1}}^{1}}^{\ell_{1}}\cdots\sigma_{i_{1}^{k}}^{\ell_{k}}\cdots\sigma_{i_{s_{k}}^{k}}^{\ell_{k}}\sigma_{i_{1}}^{\ell_{k+1}}\cdots\sigma_{i_{p}}^{\ell_{k+1}}\rangle\Big].
				\end{aligned}
				\end{align*}
				Hence, by induction on $k$, for appropriate constants which satisfy
				$|c_{k}(\ell_{1},\ldots,\ell_{k})|\leq(k-1)!$ we have that 
				\begin{equation}
				\begin{aligned} & \frac{\partial}{\partial J_{i_{1}^{1},...,i_{p_{1}}^{1}}^{(p_{1})}}\cdots\frac{\partial}{\partial J_{i_{1}^{k},...,i_{p_{k}}^{k}}^{(p_{k})}}NF_{N,\beta}\\
				& =\beta^{k}N^{-\frac{\sum_{i=1}^{k}(p_{i}-1)}{2}}\prod_{i=1}^{k}\gamma_{p_{i}}\sum_{1\leq\ell_{1},\ldots,\ell_{k}\leq k}c_{k}(\ell_{1},\ldots,\ell_{k})\langle\sigma_{i_{1}^{1}}^{\ell_{1}}\cdots\sigma_{i_{p_{1}}^{1}}^{\ell_{1}}\cdots\sigma_{i_{1}^{k}}^{\ell_{k}}\cdots\sigma_{i_{p_{k}}^{k}}^{\ell_{k}}\rangle.
				\end{aligned}
				\label{eq:partialJF}
				\end{equation}
				
				Define 
				\[
				\Delta_{k}:=\sum\Big(\E\frac{\partial}{\partial J_{i_{1}^{1},...,i_{p_{1}}^{1}}^{(p_{1})}}\cdots\frac{\partial}{\partial J_{i_{1}^{k},...,i_{p_{k}}^{k}}^{(p_{k})}}NF_{N,\beta}\Big)^{2}
				\]
				where the sum is over all indices $(i_{1}^{j},...,i_{p_{j}}^{j})_{j=1}^{k}$
				with $2\leq p_{j}\leq p_{0}$ and $1\leq i_{l}^{j}\leq N$.

				Let $\tilde{\bs}^{1},\,\tilde{\bs}^{2},\ldots$ be i.i.d. samples
				from the Gibbs measure corresponding to an independent copy $\tilde{H}_{N}(\bs)$
				of the Hamiltonian $H_{N}(\bs)$. Then, 
				\begin{equation}
				\begin{aligned}
				& \Big(\sum_{1\leq\ell_{1},\ldots,\ell_{k}\leq k}\E \langle\sigma_{i_{1}^{1}}^{\ell_{1}}\cdots\sigma_{i_{p_{1}}^{1}}^{\ell_{1}}\cdots\sigma_{i_{1}^{k}}^{\ell_{k}}\cdots\sigma_{i_{p_{k}}^{k}}^{\ell_{k}}\rangle\Big)^2\\
				& \leq k^k \sum_{1\leq\ell_{1},\ldots,\ell_{k}\leq k}\big( \E \langle\sigma_{i_{1}^{1}}^{\ell_{1}}\cdots\sigma_{i_{p_{1}}^{1}}^{\ell_{1}}\cdots\sigma_{i_{1}^{k}}^{\ell_{k}}\cdots\sigma_{i_{p_{k}}^{k}}^{\ell_{k}}\rangle\big)^2\\
				& =k^k \sum_{1\leq\ell_{1},\ldots,\ell_{k}\leq k} \E\langle\sigma_{i_{1}^{1}}^{\ell_{1}}\tilde{\sigma}_{i_{1}^{1}}^{\ell_{1}}\cdots\sigma_{i_{p_{1}}^{1}}^{\ell_{1}}\tilde{\sigma}_{i_{p_{1}}^{1}}^{\ell_{1}}\cdots\sigma_{i_{1}^{k}}^{\ell_{k}}\tilde{\sigma}_{i_{1}^{k}}^{\ell_{k}} \cdots
				\sigma_{i_{p_{k}}^{k}}^{\ell_{k}}\tilde{\sigma}_{i_{p_{k}}^{k}}^{\ell_{k}}\rangle.
				\end{aligned}
				\label{eq:EJ}
				\end{equation}

				Combining
				(\ref{eq:EJ}) and (\ref{eq:partialJF}) we obtain that
				\begin{align*}
					\Delta_{k}  \leq&\sum_{2\leq p_{1},\ldots,p_{k}\leq p_{0}}\beta^{2k}N^{-\sum_{i=1}^{k}(p_{i}-1)}k^{3k}\\&\prod_{i=1}^{k}\gamma_{p_{i}}^{2} \sum_{1\leq\ell_{1},\ldots,\ell_{k}\leq k}\E\langle(\bs^{\ell_{1}}\cdot\tilde{\bs}^{\ell_{1}})^{p_{1}}\cdots(\bs^{\ell_{k}}\cdot\tilde{\bs}^{\ell_{k}})^{p_{k}}\rangle.
					\end{align*}
					By H\"{o}lder's inequality,
					\[
					\E\langle(\bs^{\ell_{1}}\cdot\tilde{\bs}^{\ell_{1}})^{p_{1}}\cdots(\bs^{\ell_{k}}\cdot\tilde{\bs}^{\ell_{k}})^{p_{k}}\rangle\leq\E\langle|\bs^{1}\cdot\tilde{\bs}^{1}|^{\sum_{i=1}^{k}p_{i}}\rangle.
					\]
					Using a co-area formula, for any $a\in(0,1)$ and large enough $N$,
					uniformly in $p\leq N^{a}$,
					\begin{align*}
						\E\langle|\bs^{1}\cdot\tilde{\bs}^{1}|^{p}\rangle & =N^{p/2}\E\langle|\sigma_{1}^{1}|^{p}\rangle=N^{p/2}\frac{1}{\sqrt{N}}\frac{\omega_{N-1}}{\omega_{N}}\int_{-\sqrt{N}}^{\sqrt{N}}|s|^{p}(1-s^{2}/N)^{\frac{N-3}{2}}ds\\
						& \leq2N^{p/2}\int_{-\sqrt{N}}^{\sqrt{N}}\frac{1}{\sqrt{2\pi}}|s|^{p}e^{-\frac{1}{2}s^{2}}ds\leq2N^{p/2}\E|X|^{p}\leq2N^{p/2}(p-1)!!,
						\end{align*}
						where $\omega_{N}=2\pi^{N/2}/\Gamma(N/2)$ is the area of the unit
						sphere in $\mathbb{R}^{N}$ and $X$ is a standard Gaussian variable.
						Hence, for large $N$ and any $k\leq N^{a}$,
						\[
						\Delta_{k}\leq2\left(k^{p_{0}+3}p_{0}^{p_{0}}\beta^{2}\sum_{p=2}^{p_{0}}\gamma_{p}^{2}N^{-\frac{p}{2}+1}\right)^{k}.
						\]
						
						Let $\nabla F_{N,\beta}=\big(\frac{\partial}{\partial J_{i_{1},...,i_{p}}^{(p)}}F_{N,\beta}\big)$
						be the array of the derivatives of $F_{N,\beta}$ w.r.t. all the disorder
						coefficients. Let $|\nabla F_{N,\beta}|$ denote its Frobenius norm
						and using (\ref{eq:Fgrad}) note that $N^{2}\E|\nabla F_{N,\beta}|^{2}\leq N\beta^{2}\sum_{p=2}^{p_{0}}\gamma_{p}^{2}$.
						By \cite[Eq. (6.3)]{ChattBook}, for any $d\geq1$,
						\[
						{\rm Var}(NF_{N,\beta})\leq\sum_{k=1}^{d-1}\frac{\Delta_{k}}{k!}+\frac{1}{d}N^{2}\E|\nabla F_{N,\beta}|^{2}.
						\]
						With the same choice of $d=c\log N/\log\log N$
						with small enough $c$ as in the proof of  \cite[Theorem 6.3]{ChattBook}, we obtain that for some large $C>0$,
						\[
						{\rm Var}(NF_{N,\beta})\leq CN\log\log N/\log N.\qedhere
						\]
						\end{proof}

For later use we will also prove the following zero-temperature analogue of Lemma \ref{lem:Flb}.
\begin{cor}
	\label{cor:GSlb}
	For any mixture $\nu(s)=\sum_{p\geq2}\gamma_{p}^{2}s^{p}$ that does
	not contain a $1$-spin interaction (i.e., $\gamma_1=0$), for any $c,\,t>0$, for
	large enough $N$,
	\begin{equation}
		\label{eq:GSsmalldev}
		\mathbb{P}\Big(\frac{1}{N}\max_{\bs\in\SN}H_{N}(\bs)<\Es(1)-t\Big)<e^{-cN}.
	\end{equation}
\end{cor}
\begin{proof}
	Denote $E_N:=\frac{1}{N}\max_{\bs\in\SN}H_{N}(\bs)$.
	Note that for any $\beta$, $\frac{1}{\beta}F_{N,\beta}\leq E_N$. If $H_N(\bs)$ is $\sqrt NL$-Lipschitz, then $H_N(\bs)\geq E_N-\sqrt \beta$ for any $\bs$ within distance $\sqrt{N\beta}/L$ from the maximizer of $H_N(\bs)$. From Lemma \ref{lem:Lip}, the Borell-TIS inequality and the concentration of the free energy,  one sees from this that $\lim_{N \to\infty}\frac{1}{\beta}\E F_{N,\beta} \to \Es(1) $ as $\beta\to\infty$. Hence, for large enough $\beta$ and $N$ the if the event in \eqref{eq:GSsmalldev} occurs then 
	\[		
		\frac{1}{\beta}F_{N,\beta}<\frac{1}{\beta}\E F_{N,\beta} -t/2.
	\]
	The corollary thus follow from Lemma \ref{lem:Flb}.
\end{proof}

\section{\label{sec:UnifConc}Uniform concentration of $\FF(m)$: proof of Proposition \ref{prop:concExpec}}

In this section we prove Proposition \ref{prop:concExpec}. We start with the following slightly weaker concentration property, uniform only over spheres. 

\begin{prop}
	\label{prop:concExpecq} 
	For any $\beta,t,c>0$ there exist $\delta_0$, $\rho_0$ and $n_0$ such that if $\delta\leq \delta_0$, $\rho\leq \rho_0$ and $n\geq n_0$, then for any $N$ and $q\in[0,1)$, 
	\begin{equation}
		\P\bigg(\,  \max_{\|m\|^2=Nq}\big|\, \FF(m,n,\rho,\delta)-\E \FF(m,n,\rho,\delta)\,\big|<t\, \bigg)>1-e^{-Nc}.\label{eq:concExpecq}
	\end{equation}
\end{prop}

\begin{proof}
Fix some $\beta,t,c>0$ and $q\in[0,1)$. 
Let $m$ such that $R(m,m)=q$ and let $(\bs_1,\ldots,\bs_{n})\in B_{N}(m,n,\rho,\delta)$,  and denote $R_{i,j}=R(\bs_i,\bs_j)$ and $R_i=R(\bs_i,m)$.
Then, for $i\neq j$, 
\[
\big| R_{i,j}-q\big| <\rho,\quad \big| R_i-q\big| \leq \delta \sqrt q.
\]
Therefore, since
\begin{align*}
	V&:=\frac1N \E \bigg\{ \Big(\sum_{i=1}^{n}(H_{N}(\bs_{i})-H_{N}(m))\Big)^2 \bigg\}
	\\&=\sum_{i=1}^{n}\sum_{j=1}^{n}\Big(
	\nu(R_{i,j})-\nu(q)-\big[\nu(R_i)-\nu(q)\big]-\big[\nu(R_j)-\nu(q)\big]
	\Big)
\end{align*}
and $\nu'(1)=\max_{t\in[-1,1]}|\nu'(t)|$, uniformly on $B_{N}(m,n,\rho,\delta)$,
\begin{align*}
	&\big| V -  n\big(\nu(1) - \nu(q)\big)  \big| <
	\nu'(1)n^2 \big( \rho + 2\delta
	\big).
\end{align*}

Note that $\FF(m,n,\rho,\delta)$ is equal to $n^{-1}$ times the free energy of the process $\sum_{i=1}^{n}(H_{N}(\bs_{i})-H_{N}(m))$ on $B_{N}(m,n,\rho,\delta)$, the covariance of which we have just estimated. 
Hence, from a well-known concentration result \cite[Theorem 1.2]{PanchenkoBook}, (the cited theorem concerns countable spaces, but can be easily extended by approximation to our setting with $B_{N}(m,n,\rho,\delta)$, say with prefactor $3$ instead of $2$, since the Gaussian process we consider is continuous) for a single vector $m$,
\begin{equation}
	\label{eq:conc_1pt}
	\begin{aligned}
		\P\bigg( \big|\, \FF(m,n,\rho,\delta)-\E \FF(m,n,\rho,\delta)\,\big|\geq t\, \bigg) &< 3 \exp \Big\{-\frac{(N n t)^2}{4\beta^2 NV}\Big\}
		\\& \leq 3 \exp \Big\{-\frac{N t^2C}{4\beta^2}\Big\},
	\end{aligned}
\end{equation}
where 
\begin{align*}
	n^2/V\geq C:=\Big(n^{-1}\nu(1)  + \nu'(1) \big( \rho + 2\delta
	\big) \Big)^{-1} .
\end{align*}

For $q=0$ this completes the proof, and from now on we will assume that $q\in(0,1)$.
Fix some $c'>c$ and let
$\mathcal E_N(L)$ be the event from Lemma \ref{lem:Lip}, with $L$ being large enough so that the bound there holds.  For some $\bs,\bs'\in \SN$, consider the integral
\begin{equation}
	\label{eq:integral}
	\int_{B_{N}(\sqrt q \bs,n,\rho,\delta)}\exp\Big\{\beta\sum_{i=1}^{n}(H_{N}(\bs_{i})-H_{N}(\sqrt q \bs))\Big\}d\bs_{1}\cdots d\bs_{n}
\end{equation}
and the same integral with $\bs'$. $\FF(\sqrt q \bs,n,\rho,\delta)$ and $\FF(\sqrt q \bs',n,\rho,\delta)$ are obtained from the latter by taking $\log$ and dividing  by $Nn$. 
Let $T$ be an orthogonal transformation such that $T\bs=\bs'$ and such that $T\bx=\bx$ for any $\bx$ which is orthogonal to $\bs$ and $\bs'$. Then the integral related to $\bs'$ is obtained from \eqref{eq:integral} if we replace in the exponent $\bs_i$ and $\bs$ by $T\bs_i$ and $T \bs$ but keep the same domain of integration. Note that for any $\bx\in\SN$, $\|T\bx - \bx\|\leq \|\bs- \bs'\|$.
From this it follows that on the event $\mathcal E_N(L)$, as a function on $\SN$, $\bs \mapsto \FF(\sqrt q \bs,n,\rho,\delta)$ has Lipschitz constant bounded by $2\beta L/\sqrt N$. 

Let $\epsilon:=t/4\beta L$
and let $M_N$ be some $\sqrt N \epsilon$-net of $\SN$, with $|M_N|<e^{N\tau}$ for some $\tau>0$. From \eqref{eq:conc_1pt} by a union bound,
\begin{equation}
	\label{eq:net}
	\begin{aligned}
		&\P\bigg( \max_{\bs\in M_N}\big|\, \FF(\sqrt q \bs ,n,\rho,\delta)-\E \FF(\sqrt q \bs ,n,\rho,\delta)\,\big|< \frac{t}{2}\, \bigg) \\&\quad>1-  
		\exp\big(N(\tau - t^2C/16\beta^2)\big).
	\end{aligned}
\end{equation}

Denote the event in \eqref{eq:net} by $\hat{\mathcal E}_N(t)$.
Then on $\hat{\mathcal E}_N(t) \cap \mathcal E_N(L)$ we have that
\begin{equation*}
	\begin{aligned} 
		&\max_{\bs\in\SN}\big|\, \FF(\sqrt q \bs ,n,\rho,\delta)-\E \FF(\sqrt q \bs ,n,\rho,\delta)\,\big|\\
		&\leq \max_{\bs\in M_N}\big|\, \FF(\sqrt q \bs ,n,\rho,\delta)-\E \FF(\sqrt q \bs ,n,\rho,\delta)\,\big|+\frac{t}{2} 
		<t,
	\end{aligned}
\end{equation*}
where we used the fact that $\E \FF(\sqrt q \bs ,n,\rho,\delta)$ is constant for $\bs\in\SN$.
Since $$\P\big\{ \hat{\mathcal E}_N(t) \cap \mathcal E_N(L) \big\}>1 -e^{-Nc'} -\exp\big(N(\tau - t^2C/16\beta^2)\big),
$$ by choosing small enough $\rho$ and $\delta$ and large enough $n$, so that $C$ is sufficiently large, the proof is completed.
\end{proof}

We now proceed to the proof of Proposition \ref{prop:concExpec}. Here we will use results from Section \ref{sec:TAPcorrection}. 
Let $\delta_N',\rho_N'\to0$ and $n_N\to\infty$ be sequences tending to their limits sufficiently slowly so that the conclusion of Proposition \ref{prop:FFlim} holds. 
Let $k_N\geq 1$ be some sequence of integers such that $k_N\to\infty$ and $k_N/N\to0$. Define $Q_N=\{(i2^{-k_N})^2: i=0,1,\ldots,2^{k_N}-1 \}$, $\delta_N=\delta_N'+2^{-k_N}$, $\rho_N=\rho_N'+2^{-k_N+1}$ and $\delta_N''=\delta_N+2^{-k_N}$, $\rho_N''=\rho_N+2^{-k_N+1}$. 

Note that for any $q\in[0,1)$ there exists $q_0\in Q_N$ such that $|\sqrt q-\sqrt{q_0}|<2^{-k_N}$, and
from \eqref{eq:FF'}, for any $\bs \in \SN$,
\begin{equation}
\label{eq:FFsandwich}
\FF\big(\sqrt{q_0} \bs,n_N,\rho_N',\delta_N'\big) \leq\FF\big(\sqrt{q} \bs,n_N,\rho_N,\delta_N\big)\leq \FF\big(\sqrt{q_0} \bs,n_N,\rho_N'',\delta_N''\big).
\end{equation}

Fix some $t,c>0$. From Proposition \ref{prop:concExpecq} and a union bound over $q\in Q_N$, for large $N$, 
\begin{equation*}
\P\bigg(\,  \max_{q\in Q_N}\max_{\bs\in\SN}\big|\, \FF(\sqrt q \bs,n_N,\rho_N',\delta_N')-\E \FF(\sqrt q \bs,n_N,\rho_N',\delta_N')\,\big|<\frac{t}{2}\, \bigg)>1-e^{-3Nc},
\end{equation*}
and the same holds with $\rho_N''$ and $\delta_N''$.

Let $s\in(0,1)$ be an arbitrary number. Note that by Lemma \ref{lem:Fbcont} and \eqref{eq:EFt}, for large enough $N$, for any $q\in[0,s]$ with $q_0\in Q_N$ as above, the expectations of all three free energies in \eqref{eq:FFsandwich} are within distance less than $t/2$ from each other.
Hence, for large $N$,
\begin{equation*}
	\P\bigg(\,  \max_{q\in [0,s]}\max_{\bs\in\SN}\big|\, \FF(\sqrt q \bs,n_N,\rho_N,\delta_N)-\E \FF(\sqrt q \bs,n_N,\rho_N,\delta_N)\,\big|< t \, \bigg)>1-e^{-2Nc}.
\end{equation*}

Lastly, using Lemma \ref{lem:Lip}, choose $L>0$ large enough so that $H_N(\bs)$ has Lipschitz constant bounded by $\sqrt NL$ with probability at least $1-e^{-2Nc}$. Also, choose $s$ close enough to $1$ so that $\sqrt{(1-s-\epsilon)}<t/2\beta L$, for some small $\epsilon>0$. Then, for large $N$, on the event that the Lipschitz bound holds, for any $\bs \in\SN$ and $q\in[s,1]$,
\[
\sup_{\bs' \in {\rm Band}(\sqrt q\bs,\delta_N)}\big|
 H_N(\bs') - H_N(\sqrt q \bs)
\big|< Nt/2\beta.
\]
From the definition of $\FF(\bs,n,\rho,\delta))$ and since it concentrates around its expectation (e.g., as in Proposition \ref{prop:concExpecq}), on the same event,
\[
\pushQED{\qed} 
\max_{q\in [s,1]}\max_{\bs\in\SN}\big|\, \FF(\sqrt q \bs,n_N,\rho_N,\delta_N)-\E \FF(\sqrt q \bs,n_N,\rho_N,\delta_N)\,\big|< t.
\qedhere
\popQED
\]

\section{\label{sec:pfThm2}Characterization of multi-samplable overlaps: proof of Theorem \ref{thm:eqiuv}}
In this section we prove Theorem \ref{thm:eqiuv}. 
Let $q\in[0,1)$ and  $m^\star\in \sqrt q\cdot \SN$ such that $H_N(m^\star)=\max_{\|m\|^2=Nq}H_N(m)$. Let $m=m_N\in\sqrt q\cdot \SN$ be some arbitrary, non-random sequence. From the Borell-TIS inequality \cite{Borell,TIS},
 the concentration of $F_{N,\beta}$  \cite[Theorem 1.2]{PanchenkoBook} and Proposition \ref{prop:concExpecq}, for any
 $\epsilon>0$, for some $c>0$, with probability at least $1-e^{-Nc}$,
 \begin{equation}
 \label{eq:conc3}
 \big|N^{-1} H_N(m^\star) -\Es(q) \big|,\ \big| \FF(m^\star)-\E \FF(m) \big|,\ \big| F_{N,\beta}-\E F_{N,\beta} \big|<\epsilon.
\end{equation}
Here we used the fact that, by symmetry, $\E \FF(m)$ does not depend on the choice of $m$.
Hence, if $q$ satisfies \eqref{eq:E+F}, for some $\epsilon_N\to0$, with probability at least $1-\epsilon_N$,
\begin{equation}
\label{eq:FFineq}
\frac{\beta}{N}H_N(m^\star)   +\FF(m^\star) -F_{N,\beta} >-\epsilon_N.
\end{equation}
From the definition of $\FF(m)$ the left-hand side of \eqref{eq:FFineq} is equal to 
\begin{equation*}
\begin{aligned}
\frac1{ Nn_N}\log G_{N,\beta}^{\otimes n_N}\Big\{\, \forall i\leq n_N:\,&\big|R( \bs_i,m^\star) - q  \big| \leq \delta_N\sqrt q,\\  \forall i< j \leq n_N: \, &\big|R( \bs_i,\bs_j) - q  \big| < \rho_N\,  \Big\}.
\end{aligned}
\end{equation*}  

Note that
\begin{equation*}
		Y_N(n)=\frac1{ N}\log G_{N,\beta}^{\otimes n}\Big\{\, \forall i< j \leq n: \, \big|R( \bs_i,\bs_j) - q  \big| < \rho_N\,  \Big\},
\end{equation*} 
is subadditive in $n$. Thus, denoting $n_N'=n\lfloor n_N/n \rfloor$, for fixed $n$ and large enough $N$, with probability at least $1-\epsilon_N$,
\[
\frac{1}{n}Y_N(n)\geq \frac{1}{n_N'}Y_N(n_N')\geq \frac{1}{n_N'}Y_N(n_N) \geq \frac{2}{n_N}Y_N(n_N)>-2\epsilon_N.
\] 
It clearly follows that $q$ is multi-samplable.

Next, suppose $q\in[0,1)$ 
is  a multi-samplable overlap. If $n_N\to\infty$ and $\rho_N\to0$ slowly enough,
\begin{equation}
\label{eq:t1}
\lim_{N\to\infty}\frac1N\log \E G_{N,\beta}^{\otimes 2n_N}\Big(\,\forall i< j\leq 2n_N:\,
\big| R(\bs_i,\bs_j)-q\hspace{.025cm} \big|<\rho_N\,
\Big)=0.
\end{equation}
Hence, for any $t>0$,
\begin{equation*}
	\lim_{N\to\infty}\frac1N\log \P \left(\frac1N \log G_{N,\beta}^{\otimes 2n_N}\Big(\,\forall i< j\leq 2n_N:\,
	\big| R(\bs_i,\bs_j)-q\hspace{.025cm} \big|<\rho_N\,
	\Big)>-t\right)=0,
\end{equation*}
and for an appropriate sequence $t_N\to0$,
\begin{equation*}
	\lim_{N\to\infty}\frac1N\log \P \left( \frac1N \log  G_{N,\beta}^{\otimes 2n_N}\Big(\,\forall i< j\leq 2n_N:\,
	\big| R(\bs_i,\bs_j)-q\hspace{.025cm} \big|<\rho_N\,
	\Big) > -t_N \right)=0.
\end{equation*}

By conditioning on the first $n_N$ points, $\bs_1,\ldots,\bs_{n_N}$, we conclude that, with probability at least $e^{-Nc_N}$ for some sequence $c_N\to0$, there exist $\bar \bs_1,\ldots,\bar \bs_{n_N}\in \SN$ such that for any $i<j \leq n_N$, $\big|R(\bar  \bs_i, \bar \bs_j) - q   \big| < \rho_N$, and such that 
\begin{align*}
	&\frac1{N}\log G_{N,\beta}^{\otimes 2n_N}\Big\{ \forall i\neq j\leq 2n_N,\, \big|R( \bs_i,\bs_j) - q   \big| < \etn  \,\Big|\, \forall i\leq n_N,\, \bs_i = \bar \bs_i \Big\}\\
	&:= \frac1{N}\log \int_T \exp\Big\{\beta\sum_{i=n_N+1}^{2n_N}H_N(\bs_i)\Big\}d\bs_{n_N+1}\cdots d\bs_{2n_N}-n_NF_{N,\beta}
	>-t_N,
\end{align*}
where $T$ is the set of points $(\bs_{n_N+1},\ldots, \bs_{2n_N})\in (\SN)^{n_N}$ such that  $|R( \bs_i,\bs_j) - q  |<\etn$ and $|R( \bs_i,\bar \bs_k) - q  |<\etn$, for any $k\leq n_N$ and $n_N+1\leq i<j\leq 2n_N$.

Setting $\bar m=n_N^{-1}\sum_{i\leq{n_N}}\bar \bs_i$, for any $(\bs_{n_N+1},\ldots, \bs_{2n_N})\in T$ and  $n_N+1\leq i\leq 2n_N$,
\[
\big|R(\bs_i,\bar m )-q\,\big|<\etn \ \ \ \mbox{and}\ \ \ \big|R(\bar m,\bar m )-q\,\big|<\frac{n_N-1}{n_N}\etn+\frac{1}{n_N}.
\]
Hence,  defining $m_\star:=0$ if $q=0$ and  $m_\star:=\sqrt{Nq}\bar  m/\|\bar m\|$ otherwise, for $\delta_N$ going to $0$ sufficiently slowly, $T\subset({\rm Band}(m_\star,\delta_N))^{n_N}$ and
therefore, with probability at least $e^{-Nc_N}$,
\begin{equation*}
\begin{aligned}	
&\frac1{N}\log \int_{B_N(m_\star,n_N,\rho_N,\delta_N)} \exp\Big\{\beta\sum_{i=n_N+1}^{2n_N}H_N(\bs_i)\Big\}d\bs_{n_N+1}\cdots d\bs_{2n_N}-n_NF_{N,\beta}\\
&=n_N\Big(
\beta \frac{H_N(m_\star)}{N}+\FF(m_\star) - F_{N,\beta}
\Big)>-t_N.
\end{aligned}
\end{equation*}

By the same concentration results we used to justify \eqref{eq:conc3}, for any
$\epsilon>0$, for some $c>0$, with probability at least $1-e^{-Nc}$,
\begin{equation}
	\label{eq:conc3_2}
	\big|N^{-1} H_N(m^\star) -\Es(q) \big|,\ \big| \FF(m_\star)-\E \FF(m) \big|,\ \big| F_{N,\beta}-\E F_{N,\beta} \big|<\epsilon.
\end{equation}
Note that by definition, $H(m^\star)\geq H_N(m_\star)$. Combining these facts with the last inequality we conclude that for large $N$,  with probability at least $e^{-Nc_N}-e^{-Nc}>0$ and therefore deterministically,
\[
\beta \Es(q)+ \E \FF(m) - \E F_{N,\beta}
>-\frac{t_N}{n_N}-3\epsilon.
\]
Since $\epsilon$ is arbitrary,
\begin{equation}
\label{eq:limineq2}
\beta \Es(q)+ \lim_{N\to\infty}\E \FF(m) \geq \lim_{N\to\infty} \E F_{N,\beta}.
\end{equation}
Here we assume that $\delta_N,\rho_N$ and $n_N$, with which $\FF(m)$ is defined (see \eqref{eq:Fhatbs}),  tend to their limits slowly enough so that by Proposition \ref{prop:FFlim} the limit of  $\E \FF(m)$ exists.
From the inequality \eqref{eq:ineq} and the concentration results above, 
\begin{equation}
\label{eq:limineq}
\beta \Es(q)+ \lim_{N\to\infty}\E \FF(m) \leq \lim_{N\to\infty} \E F_{N,\beta}.
\end{equation}
(Note that for this inequality $q$ does not have to be multi-samplable.) 
This implies \eqref{eq:E+F} and completes the proof. \qed

\begin{rem}
	\label{rem:multi_limsup}
	Assume that for some $q\in[0,1)$, Definition \ref{def:multi} holds only with $\limsup$ in \eqref{eq:goodq}. Then, using the same argument in the proof above, \eqref{eq:E+F} holds on some subsequence $N_i$. However, all limits in \eqref{eq:E+F} exist, using Proposition \ref{prop:FFlim}. Hence, \eqref{eq:E+F} holds also without taking a subsequence. From the first part of the proof of Theorem \ref{thm:eqiuv}, it then follows that \eqref{eq:goodq} also holds with $\lim$ instead of $\limsup$.
\end{rem}

\begin{rem}
	\label{rem:exists_q}
	Let $n$ and $k$ be two natural numbers and consider a partition of $[-1,1]$ into $k$ intervals of equal length. Note that 
	by Ramsey's theorem, for large enough $K=K(n)$, if we take $K$ samples from the Gibbs measure, then $n$ of them fall into the same interval and therefore
	\eqref{eq:goodq} holds with some $q\in[-1,1]$ and $\rho=1/k$, on some subsequence $N_i$. Hence, from compactness,  there exists $q\in[-1,1]$ which satisfies \eqref{eq:goodq}  with arbitrary $n$ and $\rho$ on some other subsequence $N_i$. Since any matrix $\{(\bs_i,\bs_j)\}_{i,j=1}^n$ of overlaps with arbitrary $\bs_i\in\SN$ must be positive semi-definite,  $q\in[0,1]$. It is also easy to verify  that $q<1$. By Remark \ref{rem:multi_limsup}  this $q$ is multi-samplable. We showed that for any model there exits $q\in[0,1)$ which is multi-samplable. This, of course, also follows from Corollary \ref{cor:ParisiSupp}, but the argument here is basic and does not rely on Parisi's formula.
\end{rem}

\section{\label{sec:TAPpfs}The Onsgaer correction at $q_{\max}$: proof of Corollary \ref{cor:classicalTAP}}
We split the proof of the corollary into three steps stated in the following auxiliary lemmas.

\begin{lem}
	\label{lem:nesting}
	Given $\beta>0$, if $q$ is a multi-samplable overlap w.r.t. $\nu(t)$ and $q'$ is a multi-samplable overlap w.r.t. $\nu_q(t)$, then $\hat q:=q+(1-q)q'$ is a multi-samplable overlap w.r.t. $\nu(t)$.
\end{lem}

\begin{lem}
	\label{lem:qmax0}
	Let $\beta>0$ and let $H_N(\bs)$ be an arbitrary spherical Hamiltonian. If $0$ is the only multi-samplable overlap, then for any $\epsilon>0$, 
	\begin{equation}\label{eq:RSoverlap}
		\lim_{N\to\infty}\E G_{N,\beta}^{\otimes 2}\left\{  
		|R(\bs_1,\bs_2)|>\epsilon 
		\right\} =0.
	\end{equation}
\end{lem}

\begin{lem}
	\label{lem:RSFE}
	Let $\beta>0$ and let  $H_N(\bs)$ be the Hamiltonian corresponding to some mixture $\nu(t)$. If \eqref{eq:RSoverlap} holds for any $\epsilon>0$, then $\lim_{N\to\infty}\E F_{N,\beta} = \frac12\beta^2\nu(1)$.
\end{lem}

Note that by Lemmas \ref{lem:0ismulti} and \ref{lem:nesting}, for $q=q_{\max}$, the set of multi-samplable overlaps for the mixture $\nu_q(t)$ is $\{0\}$. Corollary \ref{cor:classicalTAP} therefore follows from Lemmas \ref{lem:qmax0} and \ref{lem:RSFE}. It remains to prove the lemmas above which we will in the following subsections.

\subsection{\label{subsec:aux1}Proof of Lemma \ref{lem:nesting}}
Suppose that $q$ and $q'$ are as in the lemma and define $\hat q:=q+(1-q)q'$. 
From Theorem 
\ref{thm:FE}, 
\begin{align}
\label{eq:TAP1}
\lim_{N\to\infty}\E F_{N,\beta} &=  \beta\Es(q)+\frac{1}{2}\log(1-q)+F_{\beta}(q),\\
\label{eq:TAP2}
\lim_{N\to\infty}\E F^q_{N,\beta} &=  \beta\Es^q(q')+\frac{1}{2}\log(1-q')+F^q_{\beta}(q'),
\end{align}	
where in the first line all terms are defined using the mixture $\nu$, and the superscript $q$ in the second line is used to denote the same but working with the mixture $\nu_q$.
What we need to prove is that
\begin{equation}
\label{eq:TAP3}
\lim_{N\to\infty}\E F_{N,\beta} \leq  \beta\Es(\hat q)+\frac{1}{2}\log(1-\hat q)+F_{\beta}(\hat q),
\end{equation}	
since  by Theorem 
\ref{thm:FE} we know that the opposite inequality holds.

For a polynomial function $f(t)$ define 
\[
\mathcal{A}_q(f)(t)
:=
f((1-q)t+q)-f(q)-f'(q)(1-q)t
,
\] so that $\nu_q = \mathcal{A}_q(\nu)$ (see 
\eqref{eq:nuq}). The free energy $F^q_{\beta}(q')$ corresponds to the mixture $\mathcal A_{q'}(\nu_q)= \mathcal A_{q'} \circ \mathcal A_{q} (\nu)$ and the free energy $F_{\beta}(\hat q)$ corresponds to $\nu_{\hat q} = \mathcal A_{\hat q}(\nu)$.
It is straightforward to verify that those two mixture are in fact equal, and thus so are the two free energies. By definition,  the left-hand side of \eqref{eq:TAP2} is equal to $F_{\beta}(q)$. 
Note that the sum 
of the two logarithmic terms in \eqref{eq:TAP1} and \eqref{eq:TAP2} is equal to the logarithmic term in \eqref{eq:TAP3}.
Thus, to prove \eqref{eq:TAP3} we need to show that
\begin{equation}
\label{eq:TAP4}
\Es(q)+\Es^q(q')\leq
\Es(\hat q).
\end{equation}	

We will show that for any $\epsilon>0$,
\begin{equation}
\label{eq:TAP5}
\lim_{N \to\infty}\P \bigg( \forall m_\star\in \sqrt q\cdot \SN:\  \max_{\substack{\|m\|^2=N\hat q \\ R(m-m_\star,m_\star)=0}} \frac{1}{N} 
	\Big(H_N(m)-H_N(m_\star)\Big)\geq \Es^q(q')-\epsilon \bigg)=1.
\end{equation}
By applying the inequality above with a point $m_\star$
such that  $H_N(m_\star)=\max_{\|m\|^2=Nq}H_N(m)$, using the Borell-TIS inequality \cite{Borell,TIS}, \eqref{eq:TAP4} follows.

Since $H_N(\bs)$ has Lipschitz constant bounded by $\sqrt NL$ for some large $L>0$ with probability going to $1$ (see Lemma \ref{lem:Lip}) and since $\epsilon>0$ is arbitrary, we may prove \eqref{eq:TAP5} by a union bound and a net argument as in the proof of Theorem \ref{thm:eqiuv} if we can show that for an arbitrary non-random $m_\star\in \sqrt q\cdot \SN$, and any $C>0$,
\begin{equation}
\label{eq:TAP6}
\limsup_{N\to\infty}\frac1N\log\P \bigg( \max_{\substack{\|m\|^2=N\hat q \\ R(m-m_\star,m_\star)=0}} \frac{1}{N} 
\Big(H_N(m)-H_N(m_\star)\Big)< \Es^q(q')-\epsilon \bigg)<-C.
\end{equation}

By computing the covariance and scaling the radius of the sphere over which we maximize  (see the discussion around \eqref{eq:Htilde}),  the probability in \eqref{eq:TAP6} is equal to 
\begin{equation}
\label{eq:TAP7}
\P \bigg( \max_{\bs\in \mathbb S^{N-2}} \frac{1}{N}\sqrt{\frac{N}{N-1}} 
\hat H_{N-1}(\bs) < \check \Es(1)-\epsilon \bigg),
\end{equation}
where $\hat H_{N}(\bs)$ is the Hamiltonian that corresponds to the mixture $\tilde \nu_q(q't)$ and $\check \Es(1)$ is the limiting ground-state energy, defined  similarly to \eqref{eq:GS}, that corresponds to the mixture $\nu_q(q't)$. 
By applying the decomposition in \eqref{eq:1spdecomp} to $\hat H_{N}(\bs)$ we may write 
\begin{equation*}
\hat H_N(\bs) = \check  H_N(\bs) + \gamma \bs \cdot \mathbf X,
\end{equation*}
where $\mathbf X$ is a vector of i.i.d. standard Gaussian variables, $\check  H_N(\bs)$ is the Hamiltonian which corresponds to the mixture  $\nu_q(q't)$, they are independent, and $\gamma^2 = \tilde \nu_q(q')- \nu_q(q')$.
To upper bound the probability in \eqref{eq:TAP7}, we may replace $\hat H_{N-1}(\bs)$ by $\check  H_{N-1}(\bs)$ and restrict $\bs$ to the subset of points that are orthogonal to $\mathbf X$, which is a sphere of co-dimension $1$. Hence, to prove \eqref{eq:TAP6}, we need to show that for any $\epsilon,C>0$,
\begin{equation}
\label{eq:TAP8}
\limsup_{N\to\infty}\frac1N\log\P \bigg( \max_{\bs\in  \SN} \frac{1}{N} 
\check  H_{N}(\bs) < \check  \Es(1)-\epsilon \bigg)<-C.
\end{equation}
Since the  mixture $\nu_q(q't)$ does not contain a $1$-spin interaction, \eqref{eq:TAP8} follows from Corollary \ref{cor:GSlb}.\qed

\subsection{\label{subsec:aux2} Proof of Lemma \ref{lem:qmax0}}
Define ${\rm Cap}(\bs_0,t)=\{\bs\in\SN:\, R(\bs,\bs_0)\geq t\}$ and 
\begin{equation*}
	q:=\sup\Big\{s\geq 0:\, \limsup_{N\to\infty} \frac1N\log \E \sup_{\bs \in \SN}G_{N,\beta}\big( {\rm Cap}(\bs,\sqrt s) \big) =0   \Big\}.
\end{equation*}
We will show that $q$ is multi-samplable. This will imply that $q=0$, from which Lemma \ref{lem:qmax0} obviously follows by conditioning on the first sample from the Gibbs measure.

Note that for any $\epsilon>0$, for some $t>0$ and $t_N>0$ going to $0$, on some subsequence $N_k$, 
\begin{equation}
\label{eq:event1}
\begin{aligned}
\lim_{N\to\infty}\frac1N \log\P \Big\{ \exists \bs_\star \in \SN&:\, G_{N,\beta}\big( {\rm Cap}(\bs_\star,\sqrt q-\epsilon) \big)>e^{-Nt_N},\\
\forall  \bs \in \SN&:\, G_{N,\beta}\big( {\rm Cap}(\bs,\sqrt q+\epsilon) \big)<e^{-Nt}
\Big\}=0.
\end{aligned}
\end{equation}

Assume that the event in \eqref{eq:event1} (which depends on the randomness in the disorder only) occurs.  Then the probability to sample a point $\bs_1\in {\rm Cap}(\bs_\star,\sqrt q-\epsilon)\setminus {\rm Cap}(\bs_\star,\sqrt q+\epsilon)$ from $G_{N,\beta}$ is at least $e^{-Nt_N}/2$ for large $N$. For such point,
write $\bs_1=c_\star \bs_\star + c_1 \hat \bs_1$, where $\hat\bs_1\in\SN$,  $R(\bs_\star,\hat\bs_1)=0$, $|c_\star-\sqrt q|<\epsilon$. Let $\eta,\rho>0$ be small numbers such that 
\begin{equation}
2\epsilon+\eta^2<\eta(\rho-3\epsilon),
\label{eq:rhoeta}
\end{equation}
and define $\boldsymbol \tau_{\pm}\in\SN$ by $\boldsymbol \tau_{\pm}:= \sqrt{1-\eta^2}\bs_\star \pm \eta \hat\bs_1$.

Still assuming the event in \eqref{eq:event1}, conditional on $\bs_1$, under $G_{N,\beta}$ the probability to sample 
an independent point $\bs_2$ such that $\bs_2\in {\rm Cap}(\bs_\star,\sqrt q-\epsilon)\setminus {\rm Cap}(\bs_\star,\sqrt q+\epsilon)$  and  $\bs_2\notin {\rm Cap}(\boldsymbol\tau_{\pm},\sqrt q+\epsilon)$ is at least $e^{-Nt_N}/2$ for large $N$. Assume that $\bs_2$ is such a point, and assume towards contradiction that $|R(\bs_2,\bs_1)-q|>\rho$. Then, assuming $\epsilon$ is small,  $|c_\star R(\bs_2,\bs_\star)-q|<3\epsilon$, and 
since 
\[
R(\bs_2,\bs_1)=c_\star R(\bs_2,\bs_\star) + c_1 R(\bs_2,\hat \bs_1),
\]
we have that $|R(\bs_2,\hat \bs_1)|>\rho-3\epsilon$. Therefore, assuming $\eta$ and $\epsilon$ are small,
\[
R(\bs_2,\boldsymbol\tau_{\pm})=\sqrt{1-\eta^2}R(\bs_2,\bs_\star) \pm\eta R(\bs_2,\hat \bs_1)\geq \sqrt q-\eta^2-\epsilon \pm\eta R(\bs_2,\hat \bs_1),
\]
and, using \eqref{eq:rhoeta}, either $R(\bs_2,\boldsymbol\tau_{+})$ or $R(\bs_2,\boldsymbol\tau_{-})$ is larger than $\sqrt q+\epsilon$. Since this is a contradiction, we conclude that $|R(\bs_2,\bs_1)-q|\leq \rho$.

On the event in \eqref{eq:event1}, given independent samples $\bs_1,\ldots,\bs_n$ from $G_{N,\beta}$ such that $\bs_1,\ldots,\bs_n\in {\rm Cap}(\bs_\star,q-\epsilon)\setminus {\rm Cap}(\bs_\star,q+\epsilon)$ and $|R(\bs_i,\bs_j)-q|\leq \rho$ for $i\neq j$, the conditional probability of sampling a new point $\bs_{n+1}$ from $G_{N,\beta}$
such that $\bs_{n+1}\in {\rm Cap}(\bs_\star,q-\epsilon)\setminus {\rm Cap}(\bs_\star,q+\epsilon)$ and $\bs_{n+1}\notin {\rm Cap}(\boldsymbol\tau^i_{\pm},q+\epsilon)$ for $i\leq n$, where $\boldsymbol\tau^i_{\pm}$ is defined similarly to $\boldsymbol\tau_{\pm}$
but using $\bs_i$, is at least  
$e^{-Nt_N}/2$, assuming that $N$ is large.

By applying the argument above to each sample separately, we also have that
for samples as in the last paragraph, $|R(\bs_{n+1},\bs_{i})-q|\leq \rho$ for $i\leq n$. Note that by taking $\epsilon$ to be small enough, we can choose $\rho$ to be as small as we wish.

From this it follows that \eqref{eq:goodq} holds with arbitrary $n$ and $\rho$ on the subsequence $N_k$, on which \eqref{eq:event1} holds.
Finally, note that by Remark \ref{rem:multi_limsup}, this is enough to  conclude that $q$ is multi-samplable. This completes the proof of the lemma.\qed

\subsection{\label{subsec:aux3} Proof of Lemma \ref{lem:RSFE}}

	The upper bound $\limsup_{N\to\infty}\E F_{N,\beta}\leq\frac{1}{2}\beta^{2}\nu(1)$
	follows from Jensen's inequality. To prove the matching lower bound we will use a simple cavity and a second moment arguments to show that 
	\begin{equation}
	\liminf_{N\to\infty}\Big(\E \log Z_{N+1,\beta}-\E\log Z_{N,\beta}\Big)\geq \frac{1}{2}\beta^{2}\nu(1).\label{eq:ZN+1}
	\end{equation}
	
	We reparameterize $H_{N}(\bs)$ as a random process on   the
	unit sphere $\US:=\{\bs\in\mathbb{R}^{N}:\,\|\bs\|=1\}$,
	\[
	h_{N}(\bs):=H_{N}(\sqrt{N}\bs)=\sum_{p=1}^{\infty}\gamma_{p}h_{N,p}(\bs),
	\]
	where
	\begin{equation}
	\label{eq:hNp}
	h_{N,p}(\bs):=\sqrt{N}\sum_{i_{1},\ldots,i_{p}=1}^{N}J_{i_{1},\ldots,i_{p}}^{(p)}\sigma_{i_{1}}\cdots\sigma_{i_{p}}.
	\end{equation}
	
	In this section, we use the notation $\boldsymbol{\rho}\in\USN$
	and $(\bs,\varepsilon)\in\US\times[-1,1]$ and assume they are related to each other
	by $\boldsymbol{\rho}=(\sqrt{1-\varepsilon^{2}}\bs,\varepsilon)$. Here, $\varepsilon$ plays the role of  the so-called cavity coordinate and the decomposition we will use is well-known.
	Let $\hat	h_{N+1}(\boldsymbol{\rho})$ have the same distribution as $h_{N+1}(\boldsymbol{\rho})$ and write 
	\[
\hat	h_{N+1}(\boldsymbol{\rho})=\hat h_{N+1}^{(1)}(\boldsymbol{\rho})+\hat h_{N+1}^{(2)}(\boldsymbol{\rho}),
	\]
	where, expressing $\hat h_{N+1}(\boldsymbol{\rho})$ as a polynomial in the coordinates of $\boldsymbol{\rho}$ using \eqref{eq:hNp},  $\hat h_{N+1}^{(1)}(\boldsymbol{\rho})$ is defined as the sum of all terms which  do not contain the last coordinate, and $\hat h_{N+1}^{(2)}(\boldsymbol{\rho})$ consists of all other terms.
	
	Now, define 
	\[
	h_{N+1}^{(1)}(\boldsymbol{\rho})= \sqrt{\tfrac{N}{N+1}} \hat h_{N+1}^{(1)}(\boldsymbol{\rho})
	\mbox{\ \ \ and \ \ }
	h_{N+1}^{(2)}(\boldsymbol{\rho})= \sqrt{\tfrac{N}{N+1}} \hat h_{N+1}^{(2)}(\boldsymbol{\rho}) + \sqrt{\tfrac{1}{N+1}}  h'_{N+1}(\boldsymbol{\rho})
	\]
	where $\hat h'_{N+1}(\boldsymbol{\rho})$ is an independent copy of $h_{N+1}(\boldsymbol{\rho})$. Then, we can define all processes on the same probability space  so that 
	\[
	h_{N+1}(\boldsymbol{\rho})=h_{N+1}^{(1)}(\boldsymbol{\rho})+h_{N+1}^{(2)}(\boldsymbol{\rho}),
	\]
	where the Gaussian fields on the right-hand side are independent. By an abuse of notation we may write
	\begin{equation}
	h_{N+1}^{(1)}(\boldsymbol{\rho})=h_{N+1}^{(1)}(\bs,\varepsilon):=\sum_{p=1}^{\infty}\gamma_{p}(1-\varepsilon^{2})^{p/2}h_{N,p}(\bs).\label{eq:Htilde1}
	\end{equation}

	We note that
	\[
	\E\left\{ h_{N+1}^{(2)}(\boldsymbol{\rho})h_{N+1}^{(2)}(\boldsymbol{\rho}')\right\} =(N+1)\nu(\boldsymbol{\rho}\cdot\boldsymbol{\rho}')-N\nu\big((1-\varepsilon^{2})^{\frac{1}{2}}(1-\varepsilon'^{2})^{\frac{1}{2}} \bs \cdot \bs' \big).
	\]
	By a Taylor approximation one obtains that for small $t>0$, whenever
	$\varepsilon,\,\varepsilon'\in(-t,t)$, 
	\begin{equation}
	\left|\E\left\{ h_{N+1}^{(2)}(\boldsymbol{\rho})h_{N+1}^{(2)}(\boldsymbol{\rho}')\right\} -\nu\left( \bs \cdot \bs'\right)-N\varepsilon\varepsilon'\nu'\left(\bs \cdot \bs'\right)\right|\leq4Nt^{4}\nu''(1)+4t^{2}\nu'(1).\label{eq:varh2}
	\end{equation}

	Fix $C>0$ and define $B_{N}:=\{\boldsymbol{\rho}\in\USN:\,|\varepsilon|<N^{-1/2}C\}$.
	Let $\E_{N}$ denote the conditional expectation given the process
	$h_{N+1}^{(1)}(\boldsymbol{\rho})$. 
\begin{lem}
\label{lem:RSFElem1}
	In the above setting, for any $\delta>0$ with probability going to $1$ as $N\to\infty$,
\begin{equation}
\begin{aligned}\int_{B_{N}}e^{\beta h_{N+1}(\boldsymbol{\rho})}d\boldsymbol{\rho} & \geq(1-\delta)\E_{N}\left\{ \int_{B_{N}}e^{\beta h_{N+1}(\boldsymbol{\rho})}d\boldsymbol{\rho}\right\} \\
& =(1-\delta)\int_{B_{N}}e^{\beta h_{N+1}^{(1)}(\boldsymbol{\rho})+\frac{1}{2}\beta^{2}\text{Var}(h_{N+1}^{(2)}(\boldsymbol{\rho}))}d\boldsymbol{\rho}.
\end{aligned}
\label{eq:hN}
\end{equation}
\end{lem}
\begin{proof}
	First note that the lemma will follow from Chebyshev's inequality
	if we show that for any $\delta>0$, with probability going to $1$,
	\[
	\E_{N}\bigg\{\left(\int_{B_{N}}e^{\beta h_{N+1}(\boldsymbol{\rho})}d\boldsymbol{\rho}\right)^{2}\bigg\}/\bigg(\E_{N}\left\{ \int_{B_{N}}e^{\beta h_{N+1}(\boldsymbol{\rho})}d\boldsymbol{\rho}\right\} \bigg)^{2}\leq1+\delta.
	\]
	Denoting $f_{N+1}(\boldsymbol{\rho})=\beta h_{N+1}^{(1)}(\boldsymbol{\rho})+\frac12 \beta^2\text{Var}(h_{N+1}^{(2)}(\boldsymbol{\rho}))$, note that the left-hand side above is equal to
	\[
	\frac{\int_{ B_N\times B_N} 
	\exp\Big\{  
	f_{N+1}(\boldsymbol{\rho})+f_{N+1}(\boldsymbol{\rho}')+\beta^2 \text{Cov}\big(h_{N+1}^{(2)}(\boldsymbol{\rho}),h_{N+1}^{(2)}(\boldsymbol{\rho}')\big)
	\Big\}d\boldsymbol{\rho}d\boldsymbol{\rho}'}
	{\int_{ B_N\times B_N} 
	\exp\Big\{  
	f_{N+1}(\boldsymbol{\rho})+f_{N+1}(\boldsymbol{\rho}')
	\Big\}d\boldsymbol{\rho}d\boldsymbol{\rho}'}
	\]
	
	If we define
	\[
	\kappa(\alpha):=
	\sup_{{|\bs \cdot \bs'|\leq\alpha,\
			|\varepsilon|,|\varepsilon'|\leq N^{-1/2}C
		}
	}\text{Cov}\big(h_{N+1}^{(2)}(\boldsymbol{\rho}),h_{N+1}^{(2)}(\boldsymbol{\rho}')\big),
	\]
	then $\kappa(1)\leq C'$ for some constant $C'>0$ and $\lim_{\alpha\to0^{+}}\kappa(\alpha)=0$.
Hence, denoting 
\[
A_{N}(I)=
\int_{(\boldsymbol{\rho},\boldsymbol{\rho}')\in B_{N}\times B_{N}:\
	\bs \cdot \bs' \in I} 
\exp\Big\{  
f_{N+1}(\boldsymbol{\rho})+f_{N+1}(\boldsymbol{\rho}')+\beta^2 \text{Cov}\big(h_{N+1}^{(2)}(\boldsymbol{\rho}),h_{N+1}^{(2)}(\boldsymbol{\rho}')\big)
\Big\}d\boldsymbol{\rho}d\boldsymbol{\rho}',
\]
it is sufficient to show that for any $\delta,\,\alpha>0$,
with probability going to $1$ as $N\to\infty$,
\begin{equation}
	\frac{A_{N}([-1,1]\setminus[-\alpha,\alpha])}{A_{N}([-\alpha,\alpha])}\leq\delta.\label{eq:AN}
\end{equation}
							
Using Dudley's entropy bound (see e.g. \cite[Theorem 1.3.3]{RFG}), it
is standard to show that $\E\max_{\bs}h_{N,p}(\bs)\leq NK\sqrt{p}$
for some universal constant $K>0$.\footnote{This can be seen using the fact that the metric $\big(\E\big\{ (h_{N,p}(\bs)-h_{N,p}(\bs'))^2\big\}\big)^{1/2}$ associated
to $h_{N,p}(\bs)$ is bounded by $\sqrt{2Np(1-\langle \bs,\bs' \rangle)}$  and a simple analysis of the covering number of the
sphere.} It follows that, for any $N$, 
\[
	\E\sup_{\boldsymbol{\rho}\in B_{N}}\left(h_{N+1}^{(1)}(\bs,\varepsilon)-h_{N}(\bs)\right)<c
\]
for some $c>0$. From the Borell-TIS inequality, for
large $N$ and any $t>0$,
\begin{equation}
	\P\left\{ \sup_{\boldsymbol{\rho}\in B_{N}}|h_{N+1}^{(1)}(\bs,\varepsilon)-h_{N}(\bs)|>c+t\right\} \leq ce^{-c't^{2}N},\label{eq:BorelhN}
\end{equation}
for appropriate $c,\,c'>0$. Hence, since $\text{Var}(h_{N+1}^{(2)}(\boldsymbol{\rho}))$
is bounded uniformly in $N$, (\ref{eq:AN}) would follow (with probability
going to $1$) if we can prove it with $A_{N}(I)$ replaced by
\[
A_N'(I):=
	\int_{{(\bs,\bs')\in\US\times\US:\
	 \bs \cdot \bs' \in I
}
}e^{\beta h_{N}(\bs)+\beta h_{N}(\bs')}d\bs d\bs'.
\]
This, however, is equivalent to our assumption that $\lim_{N\to\infty}\E G_{N,\beta}^{\otimes 2}\left\{  
|R(\bs_1,\bs_2)|>\alpha
\right\} =0$,
and the proof is completed.
\end{proof}
\begin{lem}
\label{lem:RSFElem2}
In the above setting, for any $\delta>0$ with probability going to $1$ as $N\to\infty$,
there exists a (random) subset $D_{N}\subset\US$ such that 
\begin{equation}
\label{eq:DN1}
	\int_{D_{N}}e^{\beta h_{N}(\bs)}d\bs  \geq(1-\delta)\int_{\US}e^{\beta h_{N}(\bs)}d\bs,
\end{equation}
and
\begin{equation}
\sup_{\bs\in D_{N}}\sum_{p=1}^{\infty}p\gamma_{p}h_{N,p}(\bs)  \leq(1+\delta)N\beta\nu'(1).
\label{eq:DN2}
\end{equation}
\end{lem}
\begin{proof}
Let $G_{N,\beta}$ and $F_{N,\beta}$ be the Gibbs measure and free
energy associated to $h_{N}(\bs)$ (which coincide with those of $H_{N}(\bs)$,
up to reparameterization of the measure). Denote 
\[
	g_{N}(\bs)=\sum_{p=1}^{\infty}p\gamma_{p}h_{N,p}(\bs)\text{\, and \,}g_{N,t}(\bs)=\sum_{p=1}^{\infty}\sqrt{2tp+t^{2}p^{2}}\gamma_{p}h_{N,p}'(\bs),
\]
where $h_{N,p}'(\bs)$ are independent copies of $h_{N,p}(\bs)$ and
note that $h_{N}(\bs)+tg_{N}(\bs)\overset{d}{=}h_{N}(\bs)+g_{N,t}(\bs)$.
											
Assume towards contradiction that for some $\delta>0$, for $N$ as
large as we wish,
\[
	\P\left\{ G_{N,\beta}\left(\left\{ \bs:\,g_{N}(\bs)>(1+\delta)N\beta\nu'(1)\right\} \right)>\delta\right\} >\delta.
\]
Then, from the concentration of the free energy \cite[Theorem 1.2]{PanchenkoBook}, on some
subsequence $N_{i}$ going to $\infty$, for any $t>0$,
\[
	\frac{1}{N}\E\log\int e^{\beta(h_{N}(\bs)+tg_{N}(\bs))}d\bs>(1+o(1))\big(\E F_{N,\beta}+(1+\delta)t\beta^{2}\nu'(1)\big).
\]
													
On the other hand, by applying Jensen's inequality conditional on
$h_{N}(\bs)$,
\[
	\frac{1}{N}\E\log\int e^{\beta(h_{N}(\bs)+g_{N,t}(\bs))}d\bs<(1+o(1))\Big(\E F_{N,\beta}+\frac{1}{2}\beta^{2}\sum_{p}(2tp+t^{2}p^{2})\gamma_{p}^{2}\Big).
\]
By taking $t$ small enough, since $\nu'(1)=\sum_p p\gamma_p^2$, we arrive at a contradiction. 
\end{proof}

By a similar argument to the one we used for (\ref{eq:BorelhN}),  for appropriate constants $c,\,c'>0$,
\[
\P\left\{ \sup_{\boldsymbol{\rho}\in B_{N}}\Big|h_{N+1}^{(1)}(\bs,\varepsilon)-\big[h_{N}(\bs)-\frac{1}{2}\varepsilon^{2}\sum_{p=1}^{\infty}p\gamma_{p}h_{N,p}(\bs)\big]\Big|>\frac{c+t}{N}\right\} \leq ce^{-c't^{2}N}.
\]

From (\ref{eq:varh2}), 
\begin{equation}
\sup_{\boldsymbol{\rho}\in B_{N}}\left|\text{Var}(h_{N+1}^{(2)}(\boldsymbol{\rho}))-\nu(1)-N\varepsilon^{2}\nu'(1)\right|=O(N^{-1}).\label{eq:Varh2}
\end{equation}

Combining the above with Lemma \ref{lem:RSFElem1}, we have that for any $\delta>0$ with probability going to $1$, 
\begin{equation}
	\int_{B_{N}}e^{\beta h_{N+1}(\boldsymbol{\rho})}d\boldsymbol{\rho}  \geq (1-\delta)\int_{B_{N}}e^{\beta\big(h_{N}(\bs)-\frac{1}{2}\varepsilon^{2}\sum_{p=1}^{\infty}p\gamma_{p}h_{N,p}(\bs)\big) +\frac{1}{2}\beta^{2}\big(\nu(1)+N\varepsilon^{2}\nu'(1)\big)}d\boldsymbol{\rho}.
	\label{eq:hN}
\end{equation}

For any smooth function $f$ on the sphere, by a coarea formula, we have that
\begin{equation}
\label{eq:coarea}
\begin{aligned}
	\int_{B_{N}}f(\boldsymbol{\rho})d\boldsymbol{\rho} & =\frac{\omega_{N-1}}{\omega_{N}}\int_{\US}\int_{-N^{-1/2}C}^{N^{-1/2}C}(1-\varepsilon^{2})^{\frac{N-3}{2}}f(\bs,\varepsilon)d\varepsilon d\bs\\
	& =(1+o_{N}(1))\sqrt{\frac{N}{2\pi}}\int_{\US}\int_{-N^{-1/2}C}^{N^{-1/2}C}e^{-\frac{N}{2}\varepsilon^{2}}f(\bs,\varepsilon)d\varepsilon d\bs,
\end{aligned}
\end{equation}
where  $\omega_{N-1}=2\pi^{N/2}/\Gamma(N/2)$ is the surface area
of the unit sphere in $\mathbb R^{N}$, where $\frac{\omega_{N-1}}{\omega_{N}}=\sqrt{N/2\pi}(1+o_{N}(1))$, and where  by an abuse of notation $f(\bs,\varepsilon):=f(\boldsymbol{\rho})$, assuming as before that
$\boldsymbol{\rho}=(\sqrt{1-\varepsilon^{2}}\bs,\varepsilon)$. 

Hence, from \eqref{eq:hN} and Lemma \ref{lem:RSFElem2}, for any $\delta>0$ with probability going to $1$, 
\begin{equation*}
\begin{aligned}
\int_{B_{N}}e^{\beta h_{N+1}(\boldsymbol{\rho})}d\boldsymbol{\rho} & \geq
(1-\delta)\sqrt{\frac{N}{2\pi}}\int_{D_N}\int_{-N^{-1/2}C}^{N^{-1/2}C}e^{-\frac{N}{2}\varepsilon^{2}}e^{\beta h_{N}(\bs) +\frac{1}{2}\beta^{2}\big(\nu(1)-N\delta \varepsilon^{2}\nu'(1)\big)}d\varepsilon d\bs\\
&\geq (1-\delta)e^{\frac{1}{2}\beta^{2}\big(\nu(1)-\delta C^2 \nu'(1)\big)}\sqrt{\frac{N}{2\pi}} \int_{-N^{-1/2}C}^{N^{-1/2}C}e^{-\frac{N}{2}\varepsilon^{2}}d\varepsilon
 \int_{D_N} e^{\beta h_{N}(\bs)} d\bs\\
 &\geq(1-\delta)^2 e^{\frac{1}{2}\beta^{2}\big(\nu(1)-\delta C^2 \nu'(1)\big)}  
 \P \{|W|\leq C\}
 \int_{\US} e^{\beta h_{N}(\bs)} d\bs,
\end{aligned}
\end{equation*}
where  $W$ is a standard normal variable, and the $(1+o_{N}(1))$ factor from \eqref{eq:coarea} was absorbed in $(1-\delta)$.
For  arbitrary $\delta'>0$, by taking $C$ to be large enough and then $\delta$ small enough, 
\begin{equation}
\label{eq:logdiff}
\log \int_{B_{N}}e^{\beta h_{N+1}(\boldsymbol{\rho})}d\boldsymbol{\rho}  
-\log \int_{\US} e^{\beta h_{N}(\bs)} d\bs\geq  \frac{1}{2}\beta^{2}\nu(1)  -\delta'.
\end{equation}

By conditioning on the process $h_{N+1}^{(1)}(\boldsymbol{\rho})$ and then using the concentration of the free energy \cite[Theorem 1.2]{PanchenkoBook}, denoting $A:=\beta^{2}(\nu(1)+C^{2}\nu'(1))$,
\begin{equation*}
	\P\bigg\{\left|\,\log\int_{B_{N}}e^{\beta h_{N+1}(\boldsymbol{\rho})}d\boldsymbol{\rho}-\log\int_{B_{N}}e^{\beta h_{N+1}^{(1)}(\boldsymbol{\rho})}d\boldsymbol{\rho}\,\right|>A+t\bigg\}\leq2e^{-\frac{t^{2}}{4A}}.
\end{equation*}
From (\ref{eq:BorelhN}) and the fact that $\V(B_{N})/\V(\US)$
is bounded from below uniformly in $N$, for some constant $c,c'>0$,
\begin{equation*}
	\P\bigg\{\left|\,\log\int_{B_{N}}e^{\beta h_{N+1}(\boldsymbol{\rho})}d\boldsymbol{\rho}-\log\int_{\US}e^{\beta h_{N}(\bs)}d\boldsymbol{\rho}\,\right|\leq A+c+2t\bigg\}\leq 1-2e^{-\frac{t^{2}}{4A}}-ce^{-c't^{2}N}.
\end{equation*}

The uniform integrability property of the last inequality combined with  \eqref{eq:logdiff}, which holds with probability going to $1$, imply
\eqref{eq:ZN+1}, which completes the proof
of Lemma \ref{lem:RSFE}. \qed

\section{\label{sec:Parisipfs}Relations to Parisi's formula: proofs of results from Section \ref{subsec:Parisi}}
In this section we will prove Theorem \ref{thm:qsamp} and Proposition \ref{prop:Parisidistributions}.
Corollary 
\ref{cor:ParisiSupp}
is an immediate consequence of Proposition \ref{prop:optimalbeta} and Theorem \ref{thm:qsamp}.  Corollary \ref{cor:topSuppParisi}  also follows quickly, and to explain why we start with the following remark. Here, and throughout the section, we use the notation from Section \ref{subsec:Parisi}.
\begin{rem}
	\label{rem:gamma1}
	If $\nu(s)=\sum_{p\geq1}\gamma_ps^p$ is a mixture such that $\gamma_1>0$, then for any distribution $x$, $\phi'(0)=\Phi(0)>1$ and therefore $0\notin \mathcal S_\beta$ and $0\notin {\rm supp}(\mu_\beta)$, by Proposition \ref{prop:optimalbeta}.
\end{rem}
If $q_P=0$, then from the remark $\gamma_1=0$. In this case,  $\nu_{q_P}(s)=\nu(s)$ by definition, and the Parisi distribution corresponding to $\nu_{q_P}(s)$ is $x_\beta^{q_P}\equiv 1$. If $q_P>0$, then by Propositions \ref{prop:optimalbeta} and \ref{prop:Parisidistributions}, 
again $x_\beta^{q_P}\equiv 1$. From this, Corollary \ref{cor:topSuppParisi} follows by plugging $x_\beta^{q_P}\equiv 1$ in the 
Parisi formula.

Before attending to the proofs, we state
an optimality condition for the zero temperature Parisi formula analogous  to  Proposition \ref{prop:optimalbeta}, which was proved by Chen and Sen \cite[Theorem 2]{ChenSen}.
Given non-decreasing, right continuous and integrable function $\alpha:[0,1)\to[0,\infty)$  and $c>0$ define
\begin{align}
	\label{eq:Psi}
	\Psi(t)&=\nu'(t)-\int_0^t\frac{ds}{(\int_s^1\alpha(r)dr+c)^2},\\
	\label{eq:psi}
	\psi(s)&=\int_s^1\Psi(t)dt,
\end{align}
and let $\eta$ be the measure on $[0,1)$ defined by $\eta([0,s])=\alpha(s)$ for any $s\in[0,1)$.
\begin{thm}[Chen and Sen \cite{ChenSen}]
	\label{thm:optimalinfty}
	The minimizer $(\alpha_\infty,c_\infty)$ of \eqref{eq:EsParisi} is the unique pair $(\alpha,c)$ such that $\Psi(1)=0$, $\min_{s\in[0,1]}\psi(s)= 0$ 
	 and 
	\begin{equation}
	\label{eq:Sinf}
	{\rm supp}(\eta)\subset \mathcal S:=\big\{s\in[0,1]: \psi(s) = 0 \big\}.
	\end{equation}
\end{thm}

\subsection{\label{subsec:pfParisidist}Proof of Proposition \ref{prop:Parisidistributions}}
If $q=0\in\mathcal S_\beta$, then by Remark \ref{rem:gamma1}, $\gamma_1=0$ and thus $\nu_q=\nu$ and \eqref{eq:xbetaq} follows. We will assume from now on that $q>0$.
To prove \eqref{eq:xbetaq}, we will show that $x_\beta(q+(1-q)s)$ satisfies the characterizing condition from Proposition \ref{prop:optimalbeta}.
Denote by $\Phi_\beta$ and $\phi_\beta$ the functions defined by \eqref{eq:Phi} and \eqref{eq:phi} that correspond to the Parisi distribution $x_\beta$.
We have that
\begin{align*}
	\Phi_q(t)&:= \beta^2 \nu_q'(t)-\int_0^t\frac{ds}{(\int_s^1  x_\beta(q+(1-q)r) dr)^2}\\
	&=\beta^2(1-q)\Big(\nu'(q+(1-q)t)-\nu'(q) \Big)-(1-q)^2\int_0^t\frac{ds}{(\int_{q+(1-q)s}^1x_\beta(r)dr)^2}\\
	&=(1-q)\Big[\beta^2\Big(\nu'(q+(1-q)t)-\nu'(q) \Big)-\int_q^{q+(1-q)t}\frac{ds}{(\int_s^1x_\beta(r)dr)^2}\Big]\\
	&=(1-q)\big(\Phi_\beta(q+(1-q)t)-\Phi_\beta(q) \big).
\end{align*}
Since $0<q\in\mathcal S_\beta$ and $\phi_{\beta}'(q)=\Phi_{\beta}(q)$, $\Phi_\beta(q)=0$. Thus,
\begin{align}\label{eq:phiqbeta1}
	\phi_q(s)&:=\int_0^s\Phi_q(t)dt=(1-q)\int_0^s\Phi_\beta(q+(1-q)t)dt\\
	&=\int_q^{q+(1-q)s}\Phi_\beta(t)dt=\phi_\beta(q+(1-q)s)-\phi_\beta(q).
\end{align}
Since $q\in\mathcal S_\beta$, from Proposition \ref{prop:optimalbeta}, $\phi_\beta(q)=\max_{s\in[0,1]}\phi_\beta(s)$ and  $0=\phi_q(0)=\max_{s\in[0,1]}\phi_q(s)$.
If $0<s$ belongs to the support of the measure corresponding to $x_\beta(q+(1-q)s)$,  then clearly $q+(1-q)s$ belongs to the support of the measure corresponding to $x_\beta$ and, from Proposition \ref{prop:optimalbeta}, $q+(1-q)s\in\mathcal S_\beta$. Therefore, $\phi_q(s)=\max_{t\in[0,1]}\phi_q(t)$. We proved that $x_\beta(q+(1-q)s)$ satisfies the condition in Proposition \ref{prop:optimalbeta}, which proves \eqref{eq:xbetaq}.

Let $\alpha(t)$ and $c$ be the function and constant in the right-hand side of \eqref{eq:alphaq}. To prove \eqref{eq:alphaq} we will show that $(\alpha(t),c)$ satisfies the characterizing condition from Theorem \ref{thm:optimalinfty}. Note that
\[
\int_s^1\beta x_\beta(qr)dr=\frac{\beta}{q}\int_{qs}^{q}x_\beta(r)dr.
\]
Thus,
\begin{align*}
	\Psi_q(t)&:=\hat\nu_q'(t)-\int_0^t\frac{ds}{(\int_s^1\alpha(r)dr+c)^2},\\
	&=q\nu'(qt)-\int_0^t\frac{ds}{(\frac{\beta}{q}\int_{qs}^1x_\beta(r)dr)^2},\\
	&=q\nu'(qt)-\frac{q}{\beta^2}\int_0^{qt}\frac{ds}{(\int_{s}^1x_\beta(r)dr)^2}=\frac{q}{\beta^2}   \Phi_\beta(qt)   ,
\end{align*}
and
\begin{align}\label{eq:phiqbeta2}
	\psi_q(s):=\int_s^1\Psi_q(t)dt=\frac{q}{\beta^2}  \int_s^1   \Phi_\beta(qt)    dt=\frac{1}{\beta^2}  \int_{qs}^q   \Phi_\beta(t)    dt=\frac{1}{\beta^2}\big(   
	\phi_\beta(q)-\phi_\beta(qs)
	\big).
\end{align}
Similarly to the proof of the first part, from the above, the fact that $(\alpha(t),c)$ satisfies the  condition from Theorem \ref{thm:optimalinfty} can be easily verified from our assumption that $q\in\mathcal S_\beta$ and Proposition \ref{prop:optimalbeta}. This proves \eqref{eq:alphaq}. \qed

\subsection{\label{subsec:pfqsamp}Proof of Theorem \ref{thm:qsamp}} 

We will break the proof into several lemmas which we will prove in the subsections below.
First, 
for $q>0$, the direct part \eqref{eq:q_direct} of Theorem \ref{thm:qsamp} will follow from Proposition \ref{prop:Parisidistributions} and the following lemma.

\begin{lem}
	\label{lem:dist_to_ms}
	For any mixture $\nu$, $\beta>0$ and $q>0$, if $x^q_\beta(t)$ and $\big(\alpha^q_\infty(t),c^q_\infty\big)$ are as in \eqref{eq:xbetaq} and \eqref{eq:alphaq} of Proposition \ref{prop:Parisidistributions}, then $q$ is multi-samplable.
\end{lem}
 
Still for $q>0$, the converse part \eqref{eq:q_converse} of Theorem \ref{thm:qsamp} will follow from the following two lemmas.
\begin{lem}
	\label{lem:ms_to_dist}
	Let $\nu$ be a generic mixture and let $\beta>0$. If $q>0$ is multi-samplable, then $x^q_\beta(t)$ and $\big(\alpha^q_\infty(t),c^q_\infty\big)$ are as in  \eqref{eq:xbetaq} and \eqref{eq:alphaq} of Proposition \ref{prop:Parisidistributions}.
\end{lem}	
\begin{lem}
	\label{lem:dist_to_Sbeta}
	For any mixture $\nu$, $\beta>0$ and $q>0$, if $x^q_\beta(t)$ and $\big(\alpha^q_\infty(t),c^q_\infty\big)$ are as in  \eqref{eq:xbetaq} and \eqref{eq:alphaq} of Proposition \ref{prop:Parisidistributions}, then $q\in S_\beta$.
\end{lem}

Finally, we treat the case $q=0$. Assume for a moment that $\gamma_1>0$ and recall that by Remark \ref{rem:gamma1}, $0\notin \mathcal S_\beta$.
For $q=0$, $\nu_q(s)=\sum_{p\geq2}\gamma_ps^p$ while $\nu(s)=\sum_{p\geq1}\gamma_ps^p$. From this it is easy to see that, when $\gamma_1>0$, $\lim_{N\to\infty}\E F_{N,\beta}>F_\beta(q)$ and $\Es(q)=0$, so that by Theorem \ref{thm:FE}, $q=0$ is not multi-samplable.

What remains is to deal with the case where  $q=0$ and $\gamma_1=0$. In this case, $\nu_q(s)=\nu(s)$ and by Theorem \ref{thm:FE}, $q=0$ is multi-samplable. Finally, the following lemma shows that $q=0$ is in $\mathcal S_\beta$.

\begin{lem}
	\label{lem:0insupp}
	Let $\beta>0$ and suppose that $\nu(t)=\sum_{p\geq 1}\gamma_p^2 t^p$ is a mixture such that $\gamma_1=0$. If $\gamma_p>0$ only for $p=2$, then $0$ belongs to $\mathcal S_\beta$. If $\gamma_p>0$ for some $p\geq 3$, then $0$ belongs to the support of the Parisi measure $\mu_\beta$ and thus in particular also to $\mathcal S_\beta$.
\end{lem}

\subsubsection{Proof of Lemma \ref{lem:dist_to_ms}}

Let $q>0$ and assume that $x^q_\beta(t)$ and $\big(\alpha^q_\infty(t),c^q_\infty\big)$ are as in \eqref{eq:xbetaq} and \eqref{eq:alphaq} of Proposition \ref{prop:Parisidistributions}. To prove that $q$ is multi-samplable, by Theorem \ref{thm:FE} it is enough to show that
	\begin{equation}
	\label{eq:F1}
		\beta \Es(q)+ \frac{1}{2}\log(1-q)+F_{\beta}(q) = \lim_{N\to\infty} \E F_{N,\beta}.
	\end{equation}
	 
	 Recall that $\Es(q)$ is the ground-state energy corresponding to the mixture $\hat\nu_q(t):=\nu(qt)$ and $F_{\beta}(q)$ is the  free energy corresponding to $\nu_q(t)$ defined in \eqref{eq:nuq}. Using  the Parisi formula at positive \eqref{eq:parisi} and zero temperature \eqref{eq:EsParisi}, to prove \eqref{eq:F1} we need to show that
	\begin{equation}
	\label{eq:F2}
		\beta \mathcal P_{\infty}(\alpha^q_\infty,c^q_\infty,\hat\nu_q)+ \frac{1}{2}\log(1-q)+\mathcal P_\beta(x_\beta^q,\nu_q) = \mathcal P_\beta(x_\beta,\nu).
	\end{equation}
	
	For $\hat q<1$ large enough so that $x_\beta^q(\hat q)=1$, denoting $\bar q= q+(1-q)\hat q$,
	\begin{align*}
		&\mathcal P_\beta(x_\beta^q,\nu_q)+\frac{1}{2}\log(1-q)\\
		&=\frac12\Big( 
		\beta^2\int_0^1x_\beta(q+(1-q)t)
		(1-q)\big[
		\nu'((1-q)t+q)-\nu'(q)
		\big]dt\\
		&+\int_0^{\hat q}\frac{dt}{\int_t^1x_\beta(q+(1-q)s)ds}+\log(1-\hat q) +\log(1-q)
		\Big) \\
		&=\frac12\Big( 
		\beta^2\int_q^1x_\beta(t)
		\big[
		\nu'(t)-\nu'(q)
		\big]dt+(1-q)\int_0^{\hat q}\frac{dt}{\int_{q+(1-q)t}^1x_\beta(s)ds}
		+\log(1-\bar q)
		\Big) \\
		&=\frac12\Big( 
		\beta^2\int_q^1x_\beta(t)
		\big[
		\nu'(t)-\nu'(q)
		\big]dt
		+\int_{q}^{\bar q}\frac{dt}{\int_{t}^1x_\beta(s)ds}
		+\log(1-\bar q)
		\Big). 
	\end{align*}

Since
	\begin{equation*}
		\int_t^1\alpha_\infty^q(s)ds+c_\infty^q = \frac{\beta}{q}\int_{qt}^1x_\beta(s)ds,
	\end{equation*}
	we have that
	\begin{align*}
		&\mathcal P_{\infty}(\alpha^q_\infty,c^q_\infty,\hat\nu_q)\\
		&=\frac12\Big(
		q\nu'(q)
		\frac{\beta}{q}\int_{0}^1x_\beta(t)dt
		-\int_0^1q^2\nu''(qt)\big(
		\int_0^t
		\beta x_\beta(qs)ds
		\big)dt
		+\frac{q}{\beta}\int_0^1\frac{dt}{\int_{qt}^1x_\beta(s)ds}
		\Big)\\
		&=\frac12\Big(
		\beta\nu'(q)
		\int_{0}^1x_\beta(t)dt
		-\int_0^1q\nu''(qt)\big(
		\int_0^{qt}
		\beta x_\beta(s)ds
		\big)dt
		+\frac{1}{\beta}\int_0^q\frac{dt}{\int_{t}^1x_\beta(s)ds}
		\Big)\\
		&=\frac12\Big(
		\beta\nu'(q)
		\int_{0}^1x_\beta(t)dt
		-\beta\int_0^q \nu''(t)\big(
		\int_0^{t}
		 x_\beta(s)ds
		\big)dt
		+\frac{1}{\beta}\int_0^q\frac{dt}{\int_{t}^1x_\beta(s)ds}
		\Big).
	\end{align*}
	
	By definition, for $\bar q<1$ large enough so that $x_\beta(\bar q)=1$,
\begin{equation*}
	\mathcal P_\beta(x_\beta,\nu) =\frac12\Big(\beta^2 
	\int_0^1x_\beta(t)\nu'(t)dt+\int_0^{\bar q}\frac{dt}{\int_t^1x_\beta(s)ds}+\log(1-\bar q) 
	\Big).
\end{equation*}	
	
	By integration by parts,
	\begin{align*}
		\int_0^qx_\beta(t)
		\big[
		\nu'(t)-\nu'(q)
		\big]dt+\int_0^q\nu''(t)\big( 
		\int_0^tx_\beta(s)ds
		\big)dt=0.
	\end{align*}

Combining all the above, yields \eqref{eq:F2} which completes the proof. \qed
\subsubsection{Proof of Lemma \ref{lem:ms_to_dist}}

Let $\nu(t)=\sum_{p\geq 1} \gamma_p^2t^p$ be a generic mixture and let $\beta>0$. Suppose that $q>0$ is multi-samplable.
By Theorem \ref{thm:FE},
 \begin{equation}
 \label{eq:F3}
 \beta \Es(q)+ \frac{1}{2}\log(1-q)+F_{\beta}(q) = \lim_{N\to\infty} \E F_{N,\beta},
 \end{equation}
 and for any other choice of the coefficients $\gamma_p$, \eqref{eq:F3} holds with an inequality $\leq$. Therefore, assuming the derivatives exist,
 \begin{equation}
 \label{eq:F4}
\beta \frac{d}{d\gamma_p} \Es(q)+ \frac{d}{d\gamma_p}F_{\beta}(q) = \frac{d}{d\gamma_p}\lim_{N\to\infty} \E F_{N,\beta},
 \end{equation}
 where we view the three terms above as functions of the coefficients $\gamma_1,\gamma_2,\ldots$, with $\beta$ and $q$ fixed.

 As can be seen from \eqref{eq:nuq}, the Hamiltonian that corresponds to $\nu_q(t)$ can be written as $\sum_{p\geq2}\gamma_p H_{N,p}^{q}(\bs)$ for some independent mixed models $H_{N,p}^{q}(\bs)$ which depend on $q$. Hence, by H\"{o}lder's inequality, $F_\beta(q)$ is a convex function of $\gamma_p$, for any $p$. Since $\Es(q)$ is the ground state energy corresponding to $\nu(qt)$, it is easy to verify that it is a convex function of $\gamma_p$, for any $p$.
 
 In the proof of Theorem 1.2 of \cite{Talag}, Talagrand showed that the derivative of $\lim_{N\to\infty} \E F_{N,\beta}$ exists and is equal to
 \begin{equation}
 \label{eq:Fderivative}
 	\frac{d}{d\gamma_p}\lim_{N\to\infty} \E F_{N,\beta} = \beta^2 \gamma_p \Big(1-\int s^p d\mu_\beta(s)\Big),
 \end{equation}
 where $\mu_\beta$ is the measure corresponding to the Parisi distribution $x_\beta$. 
 Using the same argument as in \cite[Theorem 1.2]{Talag} and the convexity of $\Es(q)$ mentioned above it is straightforward to show that its derivative exists and is equal to 
 \begin{equation*}
 	\frac{d}{d\gamma_p} \Es(q) = \gamma_p p q^p \Big(
 	\int_0^1s^{p-1}\alpha_\infty^q(s)ds+c_\infty^q
 	\Big).
 \end{equation*}

 Again, from the argument in \cite[Theorem 1.2]{Talag} and the convexity of $F_\beta(q)$ mentioned above, one can calculate its derivative. One obtains the following, which can also be formally derived from \eqref{eq:Fderivative}  and the chain rule since 
 $F_\beta(q)$ is the free energy corresponding to the mixture $\nu_q(t)=\sum_{k=2}^{\infty}a_{k}^2(q)t^{k}$, 
  \begin{align*}
  \frac{d}{d\gamma_p} F_{\beta}(q) &= \beta^2 
  \sum_{k=2}^{\infty}a_{k}(q)\frac{d}{d\gamma_p}\big(a_{k}(q)\big)\Big(1-\int s^k d\mu^q_\beta(s)\Big)\\
  & = \beta^2 
  \sum_{k=2}^{\infty}\sum_{p= k}^\infty\gamma_p\binom{p}{k}(1-q)^{k}q^{p-k}\Big(1-\int s^k d\mu^q_\beta(s)\Big),
  \end{align*}
  where $\mu^q_\beta$ is the measure corresponding to the Parisi distribution $x_\beta^q$ of the mixture $\nu_q(t)$. 
  
Suppose that $u_p\in\mathbb R$ is a sequence with only finitely many non-zero elements, and such that if $u_p\neq 0$, then necessarily $\gamma_p> 0$. If we let  $f(s)=\sum_{p\geq1}u_p s^p$ be the corresponding polynomial,  then  
 \begin{align*}
 \sum_{p:\,u_p\neq 0}\frac{u_p}{\gamma_p} \frac{d}{d\gamma_p}\lim_{N\to\infty} \E F_{N,\beta} &= \beta^2 \int \Big(f(1)-f(s) \Big)d\mu_\beta(s),\\
 \sum_{p:\,u_p\neq 0}\frac{u_p}{\gamma_p}\frac{d}{d\gamma_p} F_{\beta}(q) &= \beta^2 \int \Big(f_q(1)-f_q(s)\Big) d\mu^q_\beta(s),\\
 \sum_{p:\,u_p\neq 0}\frac{u_p}{\gamma_p}\frac{d}{d\gamma_p} \Es(q) &= q\int f'(qs)\alpha_\infty^q(s)ds  +qf'(q)c_\infty^q,
 \end{align*}
 where, similarly to the definition of $\nu_q(t)$,
\begin{equation*}
f_q(s)=f((1-q)s+q)-f(q)-(1-q)f'(q)s.
\end{equation*}

If $\mu$ is a measure on $[0,1]$ with cumulative distribution $x(s)=\mu([0,s])$, then
\begin{align*}
\int_{[0,1]} \big(f(1)-f(s)\big)d\mu(s) & = \int_{[0,1]} \int_{[s,1]}f'(t)dtd\mu(s) \\
&=\int_{[0,1]} \int_{[0,t]} f'(t)d\mu(s)dt=\int_{[0,1]} x(t)f'(t)dt.
\end{align*}

Thus, from \eqref{eq:F4}, 
\begin{equation*}
\begin{aligned}
\beta q\int_0^1 f'(qs)\alpha_\infty^q(s)ds  +\beta  qf'(q)c_\infty^q
+
\beta^2 \int_0^1 f_q'(s)x_\beta^q(s)ds 
= 
\beta^2 \int_0^1 f'(s)x_\beta(s)ds.
\end{aligned}
\end{equation*}

Let $g(s)$ be some continuous real function on $[0,1]$ with $g(0)=0$. From the  M\"{u}ntz-Sz\'{a}sz Theorem and the assumption that $\nu$ is generic, there exists a sequence of functions $f_n(s)=\sum_{p\geq 1}u_{n,p} s^p$ as above, such that $f_n'(s)$ converges in supremum norm to $g(s)$. Also using the fact that 
\begin{equation*}
	f'_q(s)=(1-q)f'((1-q)s+q)-(1-q)f'(q),
\end{equation*}
we obtain that
\begin{equation*}
\begin{aligned}
&\beta \int_0^1 g(s)x_\beta(s)ds= q\int_0^1 g(qs)\alpha_\infty^q(s)ds  
+  qg(q)c_\infty^q\\
&\quad +
\beta (1-q)\int_0^1 g((1-q)s+q)x_\beta^q(s)ds
-\beta (1-q)g(q)\int_0^1 x_\beta^q(s)ds.
\end{aligned}
\end{equation*}
By a change of variables,
\begin{equation*}
\begin{aligned}
&\beta \int_0^1 g(s)x_\beta(s)ds= \int_0^q g(s)\alpha_\infty^q(s/q)ds 
+  qg(q)c_\infty^q\\
&\quad +
\beta \int_q^1 g(s)x_\beta^q\Big(\frac{s-q}{1-q}\Big)ds
-\beta (1-q)g(q)\int_0^1 x_\beta^q(s)ds.
\end{aligned}
\end{equation*}

By considering functions $g(s)$ such that $g(s)=0$ for all $t\leq q$, we have that $x_\infty^q(s)$ is as in \eqref{eq:xbetaq}.
By considering functions $g(s)$ such that $g(s)=0$ for all $t\geq q$, we have that $\alpha_\infty^q(s)$ is as in \eqref{eq:alphaq}.  Finally, substituting the latter two in the above equation (or by approximating the indicator of the singleton $\{q\}$ by a sequence of functions $g_n(s)$), we obtain that 
\[
qg(q)c_\infty^q-\beta (1-q)g(q)\int_0^1 x_\beta^q(s)ds=0,
\]
and thus $c_\infty^q(s)$ is as in \eqref{eq:alphaq}. \qed

\subsubsection{Proof of Lemma \ref{lem:0insupp}}
The case where $\nu(t)=\gamma_2^2t^2$ was discussed in Remark \ref{rem:Sbeta}. Suppose that $\gamma_p>0$ for some $p\geq 3$. Assume towards contradiction that $q_0:=\min {\rm supp}(\mu_\beta)>0$. Let $\Phi_\beta$ and $\phi_\beta$ be the functions \eqref{eq:Phi} and \eqref{eq:phi} corresponding to the Parisi distribution $x_\beta$. For any $s<q_0$, $x_\beta(s)=0$ and therefore for any $s\leq q_0$,
$$\Phi_\beta(s)=\beta^2\nu'(s)-Cs \mbox{\ \ \ where \ \ } C=\big(\int_{q_0}^1x(r)dr\big)^{-2}.$$

By Proposition \ref{prop:optimalbeta}, $q_0\in \mathcal S_\beta$. Hence, $\phi_\beta'(q_0)=\Phi_\beta(q_0)=0$ and
$\phi_\beta''(q_0)=\Phi_\beta'(q_0)\leq 0$. Namely,
\begin{align*}
	\beta^2\nu'(q_0)/q_0&=\beta^2\sum_{p\geq2} \gamma_p^2 p q_0^{p-2}=C,\\
	\beta^2\nu''(q_0)&=\beta^2\sum_{p\geq2} \gamma_p^2 p(p-1) q_0^{p-2}\leq C.
\end{align*}

This is, of course, a contradiction since $p-1>1$ for $p\geq3$, which  completes the proof. \qed

\subsubsection{Proof of Lemma \ref{lem:dist_to_Sbeta}}

Let $q>0$ and assume that $x^q_\beta(t)$ and $\big(\alpha^q_\infty(t),c^q_\infty\big)$ are as in  \eqref{eq:xbetaq} and \eqref{eq:alphaq}. 
We denote by $\phi_\beta(t)$ and by $\phi_q(t)$  the functions defined by \eqref{eq:phi} corresponding to the mixture $\nu(s)$ and its Parisi distribution $x_\beta$ and corresponding to the mixture $\nu_q(s)$ and its Parisi distribution $x_\beta^q$, respectively, and denote by  $\psi_q(t)$  the function defined by \eqref{eq:psi} corresponding to the mixture $\nu(qs)$ and $\big(\alpha^q_\infty(t),c^q_\infty\big)$

In \eqref{eq:phiqbeta1} and \eqref{eq:phiqbeta2}, we saw that
\begin{align*}
	\phi_q(s)&=\phi_\beta(q+(1-q)s)-\phi_\beta(q),\\
	\psi_q(s)&=\frac{1}{\beta^2}\big(   
	\phi_\beta(q)-\phi_\beta(qs)
	\big).
\end{align*}
By Theorem \ref{thm:optimalinfty}, for any $s\in[0,1]$,  $\psi_q(s)\geq \psi_q(1)=0$.
By Proposition \ref{prop:optimalbeta} and Lemma \ref{lem:0insupp} , for any $s\in[0,1]$, and $\phi_q(s)\leq \phi_q(0)=0$. It clearly follows that 
\begin{equation}
\label{eq:phiopt}
\forall s\in[0,1]:\ \ \phi(q)\geq\phi(s),
\end{equation}
and therefore $q\in S_\beta$.\qed

\section{\label{sec:Ultra}The ultrametric tree: proofs of Corollaries \ref{cor:embed}, \ref{cor:Eontree} and \ref{cor:VTAPsol}}

\subsection{Proof of Corollary \ref{cor:embed}}

Let $\nu(s)$ be a generic mixture and $\beta>0$ be such that $|{\rm supp}(\mu_\beta)|>1$. 
The proof will be based on the main results of Jagannath \cite{JagannathApxUlt}. Those results apply for generic spherical models, see Section 2 in \cite{JagannathApxUlt}. 
Genericity is required for two reasons. First, it implies the Ghirlanda-Guerra identities, and therefore ultrametricity \cite{ultramet}. Second, (see \cite[Theorem 1.2]{Talag})
\begin{equation}
\label{eq:odist}
\lim_{N\to\infty}\E G_{N,\beta}^{\otimes 2}\{R(\bs_1,\bs_2)\in \cdot\}=\mu_\beta.
\end{equation}

Similarly to the notation of \cite{JagannathApxUlt}, let
\[
\tau_{d,r}:=\Big(\cup_{k=1}^r \{1,\ldots,d\}^k  \Big)\cup \{
\varnothing \},
\]
and assume a tree structure on $\tau_{d,r}$ with  edges between $\varnothing $ and each vertex $(i_1)$ and edges between each two vertices $(i_1,\ldots,i_k)$ and $(i_1,\ldots,i_k,i_{k+1})$. 
Using the terminology of \cite{JagannathApxUlt}, suppose that the sequence
\begin{equation}
\label{eq:qseq}
0=q_0<q_1<\cdots<q_r<q_{r+1}=1
\end{equation}
is $\mu_\beta$-admissible. Namely,  assume that $\mu_\beta(\{q_i\})=0$ for $1\leq i\leq r$ and $\mu_\beta([q_i,q_{i+1}])>0$ for $0\leq i\leq r$. 
The next corollary follows from \cite[Theorems 1.5, 1.7]{JagannathApxUlt}.
\begin{cor}[Jagannath \cite{JagannathApxUlt}]
	There exist a sequence $d_N\to\infty$ and a (random) collection of subsets $\big(C_{\alpha,N}\big)_{\alpha\in\tau_{d_N,r}}$  of $\SN$ such that the following holds. Let $\ell_1,\ldots,\ell_{K_N}$, $K_N=d_N^r$, be an enumeration of the leaves of $\tau_{d_N,r}$ and set $B_i:=B_{i,N}:=C_{\ell_i,N}$ for the corresponding subsets.  Assume that $B_i$ are ordered so that $G_{N,\beta}(B_i)\geq B_{N,\beta}(B_{i+1})$. Then:
	\begin{enumerate}
		\item \label{cor:J1} If $\alpha$ is a descendant of $\alpha'$ then $C_{\alpha,N}\subset C_{\alpha',N}$, and if neither of $\alpha$ and $\alpha'$ is a descendant of the other then  $C_{\alpha,N}\cap C_{\alpha',N}=\varnothing$.
		\item \label{cor:J2} The weights $\big(G_{N,\beta}(C_{\alpha,N})\big)_{\alpha\in\tau_{d_N,r}}$ 
		converge to a Ruelle probability cascade with parameters $\mu_\beta([0,q_{1})),\ldots,\mu_\beta([0,q_{r})),1$.
		\item \label{cor:J3} The weights $\big(G_{N,\beta}(B_{i})\big)_{i\leq K_N}$ 
		converge to a Poisson-Dirichlet process of parameter $1-\mu_\beta([q_{r},1])$.
		\item \label{cor:J4} For some sequence $t_N\to0$, with probability going to $1$, for any $\ell_i$ and $\ell_j$, if $|\ell_i\wedge \ell_j|=n$, then 
		\begin{equation}
		\label{eq:clusteroverlap1}
		G_{N,\beta}^{\otimes2}\Big( \bs_1\in B_i,\ \bs_2\in B_j,\ R(\bs_1,\bs_2)\notin (q_n-t_N,q_{n+1}+t_N)  \Big)\leq t_N.
		\end{equation}
	\end{enumerate}
\end{cor}

We make several remarks. By Point \eqref{cor:J1}, $B_i$ are disjoint sets. In Point \eqref{cor:J2}, for the convergence it is assumed  that the collection of weights is extended to an infinite tree by adding zeros.
Similarly, in  Point \eqref{cor:J3} the sequence is extended by zeros to an infinite sequence. Point \eqref{cor:J3} follows from Point \eqref{cor:J2}. Lastly, in the notation of \cite{JagannathApxUlt}, \eqref{eq:clusteroverlap1} holds with $t_N=\max\{a_N,2b_N\}$.

From the convergence to a Ruelle probability cascade in Point \eqref{cor:J2} we have the following. Let $c_N>0$ be an arbitrary sequence going to $0$. By decreasing $d_N\to\infty$ if needed, keeping the children with largest weights $G_{N,\beta}(C_{\alpha,N})$ at each vertex of the tree $\tau_{d_N,r}$, we have that the corollary above still holds, and in addition we can to it the following point:
\begin{enumerate}
	\item[(5)] \label{cor:J5} With probability going to $1$,
	\begin{equation*}
	\min_{i\leq K_N}G_{N,\beta}(B_i)>c_N	.
	\end{equation*}
\end{enumerate}

Denote $S={\rm supp}(\mu_\beta)$.
If $q_n\in S$ set $q_n^+:=q_n$, otherwise define $q_n^+:=\max\{q:(q_n,q)\cap S=\varnothing \}$. Similarly, if $q_{n+1}\in S$ let $q_{n+1}^-:=q_{n+1}$ and otherwise define $q_{n+1}^-:=\min\{q:(q,q_{n+1})\cap S=\varnothing \}$. Note that by \eqref{eq:odist}, Point \eqref{cor:J4} also holds if we replace $q_n$ and $q_{n+1}$ by $q_n^+$ and $q_{n+1}^-$, and increase $t_N$ if needed.
Let us also assume that $c_N$ is large enough so that $t_N/c_N^2\to0$. Then, again by increasing $t_N\to0$ if needed, the corollary holds also if we replace Point \eqref{cor:J4} by the following:
\begin{enumerate}
\item[(4')] \label{cor:J4'} For some sequence $t_N\to0$, with probability going to $1$, for any $\ell_i$ and $\ell_j$, if $|\ell_i\wedge \ell_j|=n$, then 
\begin{equation*}
G_{N,\beta}^{\otimes2}\Big( R(\bs_1,\bs_2)\notin (q_n^+-t_N,q_{n+1}^-+t_N) \, \Big| \, \bs_1\in B_i, \bs_2\in B_j \Big)\leq t_N.
\end{equation*}
\end{enumerate}

We now allow $q_r=q_{r,N}$ to depend on $N$, keeping all other $q_i$, $i\leq r-1$, the same. Specifically, consider some sequence such that $q_{r,N}\nearrow q_P:=\max S$ . Then by a diagonalization argument,  repeating each element of $q_{r,N}$ as many times as needed and working for each $N$ with the corresponding sets $B_i$, from properties of the Poisson-Dirichlet distribution one can conclude the following.
\begin{cor}
	\label{cor:qrtoqP}
	Assume  $q_r=q_{r,N}$ is as described above and let $c_N>0$ be an arbitrary sequence going to $0$. 
	Then, there exist a sequence $d_N\to\infty$ and a (random) collection of disjoint subsets of $\SN$, $B_i:=B_{i,N}$, $i\leq K_N=d_N^r$,  an enumeration $\ell_1,\ldots,\ell_{K_N}$ of the leaves of $\tau_{d_N,r}$, and a sequence $t_N>0$ going to $0$ such that:
	\begin{enumerate}
		\item \label{enum:qrtoqP1}Point \eqref{enu:t1} of Corollary \ref{cor:embed} and the statement about $G_{N,\beta}(B_i)$ after Point \eqref{enu:t3} of Corollary \ref{cor:embed} hold.
		\item \label{enum:qrtoqP2}For any $n\leq r-1$, Point (4') above holds.
		\item \label{enum:qrtoqP3}For $n=r$, with probability going to $1$, for any $i\leq k_N$, \begin{equation*}
			G_{N,\beta}^{\otimes2}\Big( R(\bs_1,\bs_2)\notin (q_{r,N}^+-t_N,q_P+t_N) \, \Big| \, \bs_1,\bs_2\in B_i \Big)\leq t_N.
		\end{equation*} 
	\end{enumerate}
\end{cor}

Assume the setting of the corollary and for the rest of the proof assume the event that the inequalities in Points \eqref{enum:qrtoqP2} and \eqref{enum:qrtoqP3} occurs.  Denote $\theta_N:=q_P-q_{r,N}^+$
and let $E_{N,\beta}$ and $E_{N,\beta}^{\otimes2}$  denote expectation w.r.t. $G_{N,\beta}$ and $G_{N,\beta}^{\otimes2}$ and define the magnetizations
\[
m_i:= \frac{1}{G_{N,\beta}(B_i)}\int_{B_i} \bs dG_{N,\beta}(\bs).	
\]
 
From Point \eqref{enum:qrtoqP3} of Corollary \ref{cor:qrtoqP}, by Jensen's inequality, for all $i\leq K_N$,
\begin{equation}
\label{eq:clusteroverlapE}
|R(m_i, m_i)-q_P|\leq E_{N,\beta}^{\otimes2}\Big\{ \big| R(\bs_1,\bs_2)-\qs\big| \, \Big| \, \bs_1 ,\bs_2\in B_i \Big\}\leq \theta_N   +3t_N.
\end{equation}
Hence, from  Proposition 7.6 of \cite{TalagrandPstates}, for any $i,j\leq K_N$ and $\bx \in \mathbb R^N$ with $\|\bx\|^2\leq N$, 
\begin{align}
	&E_{N,\beta}^{\otimes2}\Big\{ \big|R(\bs_1,\bs_2)-R( m_i, m_j)\big| \, \Big| \, \bs_1\in B_i, \bs_2\in B_j \Big\}\leq8\sqrt{\theta_N   +3t_N} ,\\
	&E_{N,\beta}\Big\{ \big|R(\bs,\bx)-R( m_i,\bx)\big| \, \Big| \, \bs\in B_i \Big\}\leq 8\sqrt{\theta_N   +3t_N} .	
\end{align}

By Markov's inequality, for $\epsilon_N=8(\theta_N   +3t_N)^{1/4}$,
\begin{align}
	\label{eq:clusteroverlap5}
	&G_{N,\beta}^{\otimes2}\Big\{ \big|R(\bs_1,\bs_2)-R( m_i,m_j)\big| > \epsilon_N \, \Big| \, \bs_1\in B_i, \bs_2\in B_j \Big\}\leq\epsilon_N ,\\
	\label{eq:clusteroverlap51}
	&G_{N,\beta}\Big\{ \big|R(\bs,\bx)-R(m_i,\bx)\big| > \epsilon_N\, \Big| \, \bs\in B_i \Big\}\leq \epsilon_N .	
\end{align}
For the rest of the proof $\epsilon_N>0$ may increase from line to line, while we keep it going to $0$.
From Point \eqref{enum:qrtoqP2} of Corollary \ref{cor:qrtoqP} and \eqref{eq:clusteroverlapE},  for large $N$ and any $\ell_i$ and $\ell_j$, 
\begin{equation}
\label{eq:mioverlap}
	 |\ell_i\wedge \ell_j|=n \implies
	  R( m_i,m_j)\in (q_n^+-\epsilon_N,q_{n+1}^-+\epsilon_N) .
\end{equation} 
In particular, for $i=j$, since $q_{r+1}^-=q_P$,
\begin{equation}
\label{eq:mioverlap2}
|R( m_i,m_i)-q_P|< q_P-q_{r,N}+\epsilon_N\longrightarrow 0.
\end{equation} 

Recall that the magnetizations $m_i$ are related to the leaves $\ell_i$ of $\tau_{d_N,r}$ by the enumeration from Corollary \ref{cor:qrtoqP}. For each $\alpha\in\tau_{d_N,r}$ which is not a leaf, we associate the vector $v_\alpha=|I(\alpha)|^{-1}\sum_{i\in I(\alpha)} m_i$, where we define $I(\alpha)$ as the set of all indices corresponding to leaves $\ell_i$ that are descendants of $\alpha$. Of course, for $\alpha=\ell_i$, setting  $I(\ell_i)=\{i\}$, the latter definition coincides with $m_i$.
Assume the tree structure on the set of vertices $\{v_\alpha:\,\alpha \in\tau_{d_N,r}\}$ induced by that of $\tau_{d_N,r}$.

For any two vertices $\alpha,\alpha'\in \tau_{d_N,r}$,
\begin{equation*}
\begin{aligned}
R(v_\alpha,v_{\alpha'})&=\frac{1}{|I(\alpha)||I(\alpha')|}\sum_{i\in I(\alpha),\,j\in I(\alpha')}R(m_i,m_j),\\
R(v_{\alpha\wedge\alpha'},v_{\alpha\wedge\alpha'})&=\frac{1}{|I(\alpha\wedge\alpha')|^2}\sum_{i,j\in I(\alpha\wedge\alpha')}R(m_i,m_j).
\end{aligned}
\end{equation*} 
In both sums, $\ell_i\wedge \ell_j=\alpha\wedge \alpha'$ for all but a proportion of $1/d_N$ at most of the pairs $i,j$. 
Thus, from \eqref{eq:mioverlap},
\begin{equation}
\label{eq:voverlap}
	|R(v_\alpha,v_{\alpha'})-R(v_{\alpha\wedge\alpha'},v_{\alpha\wedge\alpha'})|<\epsilon_N+\eta,
\end{equation}
where we define
\begin{equation}
\label{eq:kappa}
\eta:= \max_{n=0,\ldots,r-1} \big(q_{n+1}^--q_{n}^+\big).
\end{equation}

By plugging this back to \eqref{eq:clusteroverlap5} and \eqref{eq:clusteroverlap51},  on the event we restricted to, we obtain that
\begin{equation}
\label{eq:overlapkappa}
	G_{N,\beta}^{\otimes2}\Big\{ \big|R(\bs_1,\bs_2)-R( m_i\wedge m_j,m_i\wedge m_j)\big| > \epsilon_N+\eta \, \Big| \, \bs_1\in B_i, \bs_2\in B_j \Big\}\leq\epsilon_N .	
\end{equation}
and that for any $v_\alpha$ which is an ancestor of $m_i$
\begin{equation}
\label{eq:Bkappa}
	G_{N,\beta}\Big\{ \big|R(\bs,v_\alpha)-R(v_\alpha,v_\alpha)\big| > \epsilon_N+\eta \, \Big| \, \bs\in B_i \Big\}\leq \epsilon_N .	
\end{equation}

Since $\tau_{d_N,r}$ is of finite depth, from \eqref{eq:Bkappa}, if we redefine $B_i$ to be its intersection with all the bands ${\rm Band}(v_\alpha,\delta)$ such that $v_\alpha$ is an ancestor of $m_i$ and $\delta=\eta+\epsilon_N$, and if we increase $\epsilon_N$, then \eqref{eq:overlapkappa} still holds and the statements we made about the Gibbs weights $G_{N,\beta}(B_i)$ still remain true, where we may need to change $c_N$ say to $c_N/2$. 
Similarly, from \eqref{eq:mioverlap2}, if we replace $m_i$ by $\sqrt{Nq}\cdot m_i/\|m_i\|$, then \eqref{eq:voverlap} and \eqref{eq:overlapkappa} still hold.

When the support $S$ is finite, this completes the proof of Corollary \ref{cor:embed} by choosing an admissible sequence \eqref{eq:qseq} such that $\eta=0$ (see \eqref{eq:kappa}). If $S$ is infinite, 
we proved that Corollary \ref{cor:embed} holds with constant $r_N=r$, if we replace $\eta_N$ and $\epsilon_N$ by $\epsilon_N+\eta$. Note that by using more overlaps in \eqref{eq:kappa}, we can make $\eta=\eta(r,S)$ go to $0$ as $r\to\infty$.
Thus, in the case where $S$ is infinite, Corollary \ref{cor:embed} follows by a standard diagonalization-type argument, taking $r=r_N$ to $\infty$ slowly enough. \qed

\subsection{Proof of Corollary \ref{cor:Eontree}}

Assume the setting of Corollary \ref{cor:embed} with some $c_N>0$ such that $\frac1N\log c_N\to 0$, 
 and let $M\subset V$ denote the set of leaves of $V$.
Note that for any $m$ with $\|m\|^2<N$, by definition,
\begin{equation}
\label{eq:1505-02}
\begin{aligned}
&\frac{\beta}{N}H_N(m)+\FF(m) - F_{N,\beta} 
\\
&=\frac1{ Nn_N}\log G_{N,\beta}^{\otimes n_N}\Big\{ \forall i< j\leq n_N:\ \big|R( \bs_i,\bs_j) - R(m,m)   \big| < \rho_N, \, \bs_i\in  {\rm Band}(m, \delta_N) \Big\}.
\end{aligned}
\end{equation}
Since this expression is non-positive, we only need to prove a lower bound for it for points in $V$.

Assuming that $\delta_N\geq \eta_N$ and that Point \eqref{enu:t12} of Corollary \ref{cor:embed} occurs, for any $m_k\in M$, \eqref{eq:1505-02} is bounded from below by
\begin{equation}
\label{eq:M}
	\frac1{ Nn_N}\log G_{N,\beta}^{\otimes n_N}\Big\{ \forall i< j\leq n_N:\ \big|R( \bs_i,\bs_j) - R(m_k,m_k)   \big| < \rho_N, \, \bs_i\in B_{k} \Big\},
\end{equation}
where $B_k$ is the subset associated to the leaf $m_k$.
Under the same assumption, for any $v\in V\setminus M$, \eqref{eq:1505-02} is bounded from below by
\begin{equation}
\label{eq:VnotM}
	\frac1{ Nn_N}\log G_{N,\beta}^{\otimes n_N}\Big\{ \forall i< j\leq n_N:\ \big|R( \bs_i,\bs_j) - R(v,v)   \big| < \rho_N, \, \bs_i\in B_i \Big\},
\end{equation}
where $B_i$ are the subsets that correspond to some leaves $m_i$, $i\leq n_N$, which we  choose such that $m_i\wedge m_j =v$ whenever $i\neq j$, provided that $n_N \leq d_N$.

Assuming that   $\rho_N\geq \epsilon_N$,  if Point \eqref{enu:t3} of Corollary \ref{cor:embed} occurs,
\begin{equation*}
	 G_{N,\beta}^{\otimes n_N}\Big\{ \forall i< j\leq n_N:\ \big|R( \bs_i,\bs_j) - R(v,v)   \big| < \rho_N \,\Big|\ 
	 \forall i\leq n_N:\bs_i\in B_i 
	 \Big\}\geq 1-n_N^2\epsilon_N,
\end{equation*}
and if Point \eqref{enu:t1} occurs,
\begin{equation*}
	G_{N,\beta}^{\otimes n_N}\Big\{ 	\forall i\leq n_N:\bs_i\in B_i 
	\Big\}\geq c_N^{n_N}.
\end{equation*}
Hence, also assuming that $\epsilon_Nn_N^2\to0$ for large $N$, 
\eqref{eq:VnotM} is bounded from below by $\frac1N\log c_N+\frac{1}{Nn_N}\log(1-n_N^2\epsilon_N))$ which tends to $0$. A similar bound holds for \eqref{eq:M}. 
This proves \eqref{eq:Eontree1}.

Assume towards contradiction that for some $t>0$, with probability that does not vanish as $N\to\infty$, there exists $v\in V$ such that
\[
\frac{H_N(v)}{N} < \max_{\|m\|=\|v\|}\frac{H_N(m)}{N} -\beta t.
\]
Let $m_\star$ be a point which attains the maximal energy above. Then, from  Proposition \ref{prop:concExpec} and \eqref{eq:Eontree1}, for an appropriate sequence $t_N\to0$ and with probability that does not vanish, we have that
\begin{equation*}
	\begin{aligned}
		\frac{\beta}{N}H_N(m_\star)+\FF(m_\star)&>\frac{\beta}{N}H_N(m_\star)+\E \FF(m_\star)-t_N
		\\&>\frac{\beta}{N}H_N(v)+\E \FF(v)+t-t_N
		\\
		&>\frac{\beta}{N}H_N(v)+ \FF(v)+t-2t_N
		 >  F_{N,\beta}+t-3t_N.
	\end{aligned}
\end{equation*}
Since for all $m$, $\frac{\beta}{N}H_N(m)+\FF(m)\leq  F_{N,\beta}$, we arrive at a contradiction. From this \eqref{eq:Eontree2} follows and the proof is completed.\qed

\subsection{Proof of Corollary \ref{cor:VTAPsol}}

Recall that from the inequality \eqref{eq:ineq}, the concentration of $F_{N,\beta}$ and of  $\FF(m)$ in Proposition \ref{prop:concExpec} and the uniform convergence  of $\E\FF(m)$ in Proposition \ref{prop:FFlim}, we concluded \eqref{eq:limsupTAP} and \eqref{eq:TAPcorrV}. Also recall that $F_{N,\beta}$ converges in probability to $\lim_{N\to\infty}\E F_{N,\beta}$. 

Denoting 
\[
f(m)=\frac{\beta}{N} H_N(m) +\frac{1}{2}\log(1-q)+F_{\beta}(q),
\]
it will be enough to show that for any $\epsilon>0$, $s\in(0,1)$ there exists $t>0$ such that with probability going to $1$,
\begin{equation}
\label{eq:gradfbd}
\sup_{m:\,\|m\|^2\leq Ns} \sup_{m':\,\|m-m'\|^2\leq Nt} \|\nabla f(m)- \nabla f(m')\|\leq \frac{\epsilon}{ \sqrt N},
\end{equation}
since this would imply that w.h.p. if there exists $v\in V$ such that $\|\nabla f(v)\|\geq 2\epsilon/\sqrt N$, then some point $m'$ within distance $\sqrt{Nt}$ from $v$ violates \eqref{eq:limsupTAP}, where we use the fact that w.h.p. uniformly in $v\in V$, $\|v\|^2/N\leq q_P+o(1)$, provided that $c_N$ (from Corollaries \ref{cor:embed} and \ref{cor:Eontree}) goes to $0$ slowly. 

Of course, it is enough to separately prove \eqref{eq:gradfbd} with $f(m)$ replaced by each of the summands in its definition. For the first one $\frac{\beta}{N} H_N(m)$, the bound as in \eqref{eq:gradfbd} follows since similarly to Lemma \ref{lem:Lip}, from \cite[Corollary C.2]{geometryMixed}, the Hessian $\nabla^2 H_N(m)$ has Lipschitz constant $C/N$ for some $C=C(\nu)$, w.r.t. the Frobenius norm. The second, logarithmic term is trivial to handle. Lastly, since $\nabla F_\beta(q) = \frac{2F'_\beta(q)}{N}\cdot m$, it  remains to show that $F'_\beta(q)$ is a continuous function of $q\in[0,1)$, which implies it is uniformly continuous on $[0,t]$ for any $t\in(0,1)$.

Let $H_{N,p}(\bs)$ denote the pure $p$-spin Hamiltonians, assumed to be independent for different $p$. Fix $q_0\in[0,1)$, recall that 
\[
F_\beta(q)=\lim_{N\to\infty}\frac1N \E  \log \int e^{\beta H_N^q(\bs)}d\bs,
\]
where $H_N^q(\bs)=\sum_p \alpha_p(q)H_{N,p}(\bs)$, and define $F_{\beta,q_0}(q)$ similarly, using 
the Hamiltonian $H_{N,q_0}^{q}(\bs):=\sum_p \big(\alpha_p(q_0)+\frac{d}{dq}\alpha_p(q_0)(q-q_0)\big)H_{N,p}(\bs)$ instead of $H_N^q(\bs)$.

Defining the interpolating free energy
\[
\varphi(t)=\lim_{N\to\infty}\frac1N \E  \log \int e^{\beta (\sqrt t H_N^q(\bs)+\sqrt{1-t} H_{N,q_0}^{q}(\bs))}d\bs,
\]
by a standard application of Gaussian integration by parts one can verify that for $q$ in a small neighborhood of $q_0$,
\begin{align*}
|F_\beta(q)&-F_{\beta,q_0}(q)|=|\varphi(1)-\varphi(0)|\\&\leq\sup_{\bs,\bs'\in \SN}\frac{\beta^2}{N}|\E H_N^q(\bs)H_N^q(\bs')-\E H_{N,q_0}^{q} (\bs) H_{N,q_0}^{q} (\bs')|\\
&=\beta^2\sup_{-1\leq t\leq 1}\Big| 
\sum_p \big(\alpha_p(q)\big)^2t^p
-\sum_p \big(\alpha_p(q_0)+\frac{d}{dq}\alpha_p(q_0)(q-q_0)\big)^2t^p\Big| \\
&\leq \beta^2 C_1 \
\sum_p \Big| \alpha_p(q) 
-\sum_p \alpha_p(q_0)+\frac{d}{dq}\alpha_p(q_0)(q-q_0)\Big| \leq \beta^2C_2(q-q_0)^2,
\end{align*}
where the last inequality follows from a Taylor expansion by a straightforward calculation using the exponential decay of $\gamma_p$, for appropriate constants $C_1,C_2>0$ depending on $(\gamma_p)_{p\geq2}$. 

From the quadratic upper bound above, $F'_\beta(q_0)= F'_{\beta,q_0}(q_0)$, provided they exist. From H\"{o}lder's inequality, $F_{\beta,q_0}(q)$ is convex in $q$ on a small neighborhood of $q_0$. Thus, using a similar argument to the proof of \cite[Theorem 1.2]{Talag}, we have that $F'_{\beta,q}(q)$, and thus $F'_\beta(q)$, exist (and are given by \eqref{qe:ddqFq}).  
E.g. by an interpolation argument as above, for any fixed $q_0$, as a function of $t$, $F_{\beta,q}(q+t)$ converges pointwise as $q\to q_0$. Since the latter functions are convex in $t$ on a neighborhood of $0$, we also have that $F'_{\beta,q}(q)$, and thus $F'_\beta(q)$, are continuous in $q$. \qed

\bibliographystyle{plain}
\bibliography{master}

\end{document}